\documentclass[leqno,10pt]{amsart}
\pdfoutput=1

\usepackage{lmodern,textcomp}
\usepackage{enumerate}
\usepackage{cite}
\usepackage{esint}

\usepackage{hyperref}

\usepackage{amsmath,amssymb,amsthm,mathtools,thmtools}
\usepackage[mathscr]{eucal}
\usepackage{bbm}
\usepackage{upgreek}
\usepackage{esint}
\usepackage{xparse}

\usepackage[capitalize,nameinlink,noabbrev]{cleveref}

\usepackage[a4paper,margin=1.15in]{geometry}
\allowdisplaybreaks

\theoremstyle{plain}
\newtheorem{theorem}{Theorem}[section]
\newtheorem{prop}{Proposition}[section]
\newtheorem{lemma}{Lemma}[section]
\newtheorem{cor}{Corollary}[section]
\newtheorem{defn}{Definition}[section]
\theoremstyle{remark}
\newtheorem{rem}{Remark}[section]
\numberwithin{equation}{section}

\DeclareMathOperator*{\tail}{{\rm Tail}} 

\DeclareMathAlphabet{\mathdutchcal}{U}{dutchcal}{m}{n}

\def\dz{\,{\rm d}z}
\def\dx{\,{\rm d}x}
\def\dv{\,{\rm d}v}
\def\dw{\,{\rm d}w}
\def\dt{\,{\rm d}t}
\def \d{\,{\rm d}}
\def\en{\mathbb{N}}
\def\er{\mathbb{R}}
\def\ern{\mathbb{R}^n}
\def\Om{\Omega}
\def\vs{\vspace{1mm}}

\newcommand{\Ac}{\mathcal{A}}

\newcommand{\Ec}{\mathcal{E}}

\newcommand{\Ic}{\mathcal{I}}
\newcommand{\Lc}{\mathcal{L}_v}
\newcommand{\Oc}{\mathcal{O}}
\newcommand{\Wc}{\mathcal{W}}

\newcommand{\Zs}{\mathsf{Z}}

\newcommand{\Xb}{\boldsymbol{X}}

\newcommand{\norma}[1]{{\left\|#1\right\|}}
\newcommand{\snr}[1]{\lvert #1\rvert}
\newcommand{\kk}{\kappa}
\newcommand{\eps}{\varepsilon}
\newcommand{\rr}{\rho}

\newcommand{\mean}[1]{-\hskip-1.075em\int_{#1}}

\def\Xint#1{\mathchoice
    {\XXint\displaystyle\textstyle{#1}}%
    {\XXint\textstyle\scriptstyle{#1}}%
    {\XXint\scriptstyle\scriptscriptstyle{#1}}%
    {\XXint\scriptscriptstyle\scriptscriptstyle{#1}}%
    \!\int}
\def\XXint#1#2#3{\setbox0=\hbox{$#1{#2#3}{\int}$}
    \vcenter{\hbox{$#2#3$}}\kern-0.5\wd0}
\def\bint{\Xint{-}}

\def\Xiint#1{\mathchoice
    {\XXiint\displaystyle\textstyle{#1}}%
    {\XXiint\textstyle\scriptstyle{#1}}%
    {\XXiint\scriptstyle\scriptscriptstyle{#1}}%
    {\XXiint\scriptscriptstyle\scriptscriptstyle{#1}}%
    \!\iint}
\def\XXiint#1#2#3{\setbox0=\hbox{$#1{#2#3}{\iint}$}
    \vcenter{\hbox{$#2#3$}}\kern-0.5\wd0}
\def\biint{\Xiint{-\!-}}

\makeatletter
\def\biiint{%
  \mathop{\mathpalette\biiint@{}}\nolimits\!\iiint}
\def\biiint@#1#2{%
  \begingroup
  \setbox0=\hbox{$\m@th#1\iiint$}%
  \dimen@=.81\wd0 
  \rlap{\hbox to\wd0{%
    \vphantom{\copy0}%
    \hfil$\m@th#1\vcenter{%
      \hrule width\dimen@ height.30pt
    }$\hfil
  }}%
  \endgroup}
\makeatother

\newcommand{\vertiii}[1]{{\left\vert\kern-0.25ex\left\vert\kern-0.25ex\left\vert #1 
		\right\vert\kern-0.25ex\right\vert\kern-0.25ex\right\vert}}
\makeatletter        
\newcommand{\linethrough}{\mathpalette\@thickbar}
\newcommand{\@thickbar}[2]{{#1\mkern0mu\vbox{
    \sbox\z@{$#1#2\mkern-0.5mu$}%
    \dimen@=\dimexpr\ht\tw@-\ht\z@+2\p@\relax 
    \hrule\@height0.5\p@ 
    \vskip\dimen@
    \box\z@}}}
\makeatother
\newcommand{\mathstrike}[1]{\ensuremath{\linethrough{#1}}}
\newcommand{\nra}[1]{\mathstrike{\lVert} #1 \rVert}

\begin{document}

\title[Kinetic integral equations]{De Giorgi-Nash-Moser theory for  kinetic equations \\ with nonlocal diffusions}


\author[F. Anceschi]{Francesca Anceschi}  
\address[Francesca Anceschi]{Dipartimento di scienza e alta tecnologia, Universit\`a degli Studi dell'Insubria, Via Valleggio 11, 22100 Como, Italy} \email{\url{francesca.anceschi@uninsubria.it}}

\author[G. Palatucci]{Giampiero Palatucci}  
\address[Giampiero Palatucci]{Dipartimento di Scienze Matematiche, Fisiche e Informatiche, Universit\`a di Parma, Campus - Parco Area delle Scienze 53/A, 43124 Parma, Italy}  \email{\url{giampiero.palatucci@unipr.it}}

\author[M. Piccinini]{Mirco Piccinini} 
\address[Mirco Piccinini]{Dipartimento di Matematica, Universit\`a di Pisa, L.go~B.~Pontecorvo~5, 56127 Pisa, Italy}
\email{\url{mirco.piccinini@dm.unipi.it}}

\makeatletter
\@namedef{subjclassname@2020}{\textup{2020} Mathematics Subject Classification}
\makeatother


\thanks{{\it Acknowledgement.} GP and MP are supported by PRIN 2022 PNRR Project SKYRDM, CUP\_D53D23018980001. GP is also supported by PRIN 2022 PNRR Project GEPSO, CUP\_D53D23005820006. The authors are supported by GNAMPA-INdAM Project ``Problemi variazionali: teoria cinetica, domini magnetici e non uniforme ellitticit\'a'' CUP\_E5324001950001. The authors are in debt with Moritz Ka{\ss}mann and Marvin~Weidner for their helpful and fundamental remarks on a preliminary version of the paper.}

\vspace{-3mm}

\begin{abstract} 
We extend the De Giorgi--Nash--Moser theory to nonlocal hypoelliptic equations arising in kinetic theory as linearized models for the non-cutoff Boltzmann equation.
Assuming that the nonlocal tail in velocity of weak solutions belongs locally to~$L^p_{t,x}$ for some
$p>N_{\tt d}/(2s)$, where~$N_{\tt d}$ is the drift homogeneous dimension and~$2s\in(0,2)$ is the order of diffusion, we prove a local~$L^2$-$L^\infty$ estimate, from which we also deduce a nonlocal strong Harnack inequality.
The tail summability threshold for the local estimate is optimal: for every~$1\leq p \leq N_{\tt d}/(2s)$, we construct a family of nonnegative, time-dependent solutions that rules out this estimate even for
the pure fractional Kolmogorov equation. 
The local boundedness estimate also allows for possibly unbounded source terms.
 \end{abstract}


\maketitle

\vspace{-3mm}
\setcounter{tocdepth}{1}
\tableofcontents 

\section{Introduction}

In this paper we establish quantitative De Giorgi--Nash--Moser ({\it DGNM}\,) estimates for nonlocal kinetic equations of Kolmogorov type.  We prove a local boundedness estimate and a strong Harnack inequality with a tail term, and show that the related tail summability threshold in the local estimate cannot be avoided. 

\smallskip

Nonlocal kinetic equations arise, in particular, as linearized models for the space non-homogeneous Boltzmann equation without cutoff. Therefore, their regularity theory must account for the interaction between diffusion in velocity and transport in space; see the introduction in~\cite{Mou18} and the references therein, alongside with Section~1 of~\cite{Sto19}.
Related Kolmogorov-type equations also occur in Mathematical Finance, for instance in models for Asian options; see~\cite{PasBook}.

\medskip
Let~$\Om$ be a bounded domain of~$\er^{1+2n}$, and let~$h$ be a possibly unbounded scalar field. Consider
\begin{equation}\label{problema}
	(\partial_t  + v\cdot \nabla_x) f = \Lc f + h  \quad \mbox{in}~\Omega\,,
\end{equation}
where the diffusion term~$\Lc$ is given by
\begin{equation}\label{operatore}
	\Lc f(t,x,v)
	:=\displaystyle  \ {p.~\!v.}\!\int_{\ern}\! \big(f(t,x,w)-f(t,x,v)\big)K(t,x,w,v)\dw.
\end{equation}
The kernel~$K:\er\times\er^n\times\er^{2n}\to[0,\infty)$
is measurable, symmetric in~$(v,w) \in \er^{2n}$ and satisfies the usual ellipticity/coercivity condition
\begin{equation}\label{def_kkk}
	\Lambda^{-1}|v-w|^{-n-2s} \leq \, K(t,x,v,w)\, \leq \Lambda|v-w|^{-n-2s}, \quad \forall\ v,w \in \er^n,~v\neq w\,,
\end{equation}
where~$\Lambda\geq 1$,~$s\in (0,1)$, and the condition above is assumed to hold for every $t$ and $x$. In what follows we often omit the $t$ and $x$ dependence in order to lighten the notation. As a prototype one can think of the usual fractional Laplacian operator~$-(-\Delta_v)^s$ with respect to the $v$-variables, namely when~$K(t,x,v,w) \approx |v-w|^{-n-2s}$; see, e.~\!g.,~\cite[Section~2]{DPV12} for further details.  

The integral in~\eqref{operatore} may be singular at the origin and must therefore be interpreted in the appropriate sense. For this, Equation~\eqref{problema} has to be understood through a natural weak formulation; see Section~\ref{sec_preliminaries} below.

\subsection{A brief recall of DGNM regularity theory}

Since the pioneering works of De Giorgi and Nash, together with Moser's subsequent contribution,  the {\it DGNM} theory has played a fundamental role in the study of uniformly {\it elliptic} and {\it parabolic} equations in divergence form. Therefore its extension to more general equations has been a central goal for several generations of mathematicians. For a long time, however, this theory remained essentially confined to equations whose diffusion acts in all directions of the phase space. This does not cover {\it kinetic equations}: as soon as the solution is not spatially homogeneous, the diffusion is coupled with a conservative Hamiltonian dynamics in position and velocity. The study of regularity properties for such equations goes back to Kolmogorov's short note~\cite{Kol34}, where the fundamental solution was explicitly computed in the plain Laplacian case. The strong regularizing properties of this kernel in the remaining variables, despite the lack of ellipticity of the equation, were a starting point for H\"ormander's seminal paper~\cite{Hor67}, where he gave general geometric conditions, based on commutator estimates, for such a regularization mechanism, called hypoellipitc, to hold. Nevertheless, due to their intrinsic difficulties, only recently the weak regularity theory for elliptic and parabolic equations found a complete counterpart in the kinetic setting only recently: see, e.~\!g.,~\cite{GIMV19} and~\cite{GM22,GI23}, where weak and strong Harnack inequalities are obtained together with H\"older regularity;~\cite{DH22} for a constructive proof of the weak Harnack inequality for rough kinetic equations. The $L^2$-$L^\infty$~estimate was first proved in~\cite{PP04}. We also mention~\cite{WZ11} for a preliminary H\"older regularity result based on an extension of Moser's original proof.

\vspace{2mm}

The picture becomes more delicate when dealing with general integro-differential operator, such as a fractional Laplacian with rough coefficients. Here, the {\it DGNM} theory has developed substantially over the last decades. In particular, after the breakthrough results by {Ka{\ss}mann}~\cite{Kas07,Kas11} on the validity of the classical Harnack inequality, a broad nonlocal theory was developed for integro-differential {\it elliptic} equations; see, e.~\!g.,~\cite{Kas09,DKP14,DKP16} and the references therein.

For nonlocal parabolic equations, despite the intrinsic dependance of the  time scale  on the order~$s$ of the diffusion, Harnack inequalities, H\"older continuity and $L^2$-$L^\infty$ estimates are proved in~\cite{KW23}, extending to the parabolic setting results from the elliptic theory in~\cite{DKP14,DKP16}.
In the latter works, the nonlocal tail quantifies the exterior contribution to local bounds and Harnack inequalities. Thus, in both the nonlocal elliptic and parabolic frameworks, as in the local ones, the corresponding {\it DGNM} theory is complete.

\medskip
For nonlocal kinetic equations such as~\eqref{problema}, much of the regularity theory has been developed in connection with the non-cutoff Boltzmann equation.
Imbert and Silvestre~\cite{IS20} proved a weak Harnack inequality for a broad class of kinetic integro-differential equations of the form~\eqref{problema}, under mild assumptions on the diffusion in velocity.
For the spatially inhomogeneous Boltzmann equation, this yields H\"older regularity in the conditional regime in which mass, energy and entropy are bounded from above and the mass is bounded away from vacuum.
Further regularity estimates in this regime were obtained in~\cite{IS22}; see also~\cite{Sto19,Loh22} for related results.
Boundedness estimates are available under additional assumptions relevant to kinetic models~\cite{Sil16,OS23,FRW24}, while other regularity results assume boundedness a priori.
Likewise, a quantitative control of the $L^\infty$-norm of solutions, namely an $L^2$-$L^\infty$ estimate, with a consequent strong Harnack inequality, was still unavailable.

\medskip
This state of affairs changed further with the counterexample of {Ka\ss mann and Weidner} in~\cite{KW24c}, where they constructed a sequence of solutions~$f_\eps: \er^{2n} \to [0,1]$ to
\begin{equation}\label{eq:stat-fp}
	v \cdot \nabla_x f + (-\Delta_v)^sf=0\,,
\end{equation}
for which
the ratio~$\sup f_\eps /\inf f_\eps$ blows up as~$\eps \to 0^+$; see~Theorem~1.1 there. This proves the failure of the Harnack inequality for~\eqref{eq:stat-fp}, and hence also for~\eqref{problema}, since the solutions~$\{f_\eps\}$ are time-independent.
A closer inspection also shows that such construction obstruct local~$L^2$-$L^\infty$ estimate for~\eqref{problema} when an additional tail-type error term is placed on the right-hand side with insufficient summability; see Formula~\eqref{def_tail_intro} below. This is in stark contrast with the parabolic and elliptic literature on fractional equations~(\!\!\cite{DKP16,KW23,KW24}) and appears a purely nonlocal kinetic  feature born from to the combination of the  diffusion with the anisotropic  drift. It is also worth observing that this  failure concerns quantitative estimates, not the qualitative smoothness of solutions to~\eqref{eq:stat-fp}; see~\cite{IS22}. Indeed, velocity averaging techniques~(\!\!\cite{Bou02}) transfer regularity from the $v$-variable to position, as in the case of purely local operators. See also the very recent paper~\cite{AN26} for related results in this direction.

\medskip

The counterexample in~\cite{KW24c} thus calls for a quantitative description of the exterior contribution, rather than a direct extension of the classical estimates. The aim of this paper is to identify such a contribution and to prove a~$L^2$-$L^\infty$ estimate with an arbitrarily small coefficient in front of the tail and use it to derive a nonlocal strong Harnack inequality.
The assumption on the tail is local $L^{p}_{t,x}$-summability for
$$
p \ >  \ 1 + \frac{n(1+2s)}{2s},
$$
and Theorem~\ref{thm_sharp} will show that the threshold cannot be lowered: for every~$p \leq 1 + n(1+2s)/(2s$ in the Lebesgue range, we construct a counterexample for the the simplest equation in this theory, i.~\!e. the fractional Kolmogorov equation.
 
This suggest that any broader regularity theory containing this model must account for the same threshold.
Theorems~\ref{thm_bdd} and~\ref{thm_strong} give the corresponding positive results for the measurable symmetric kernels considered here. The proofs distinguish the dependence of the tail on the drift variables from that of a general source on all kinetic variables.

\subsection{Main results} {Since the underlying geometry of equations~\eqref{problema} is determined by a homogeneous Lie group structure, to state our main results, which reflect this non-Euclidean background, we endow~$\er^{1+2n}$ with the Galilean transformation
\begin{equation} \label{def:action}
	z_{0} \circ z : = (t+t_{0},\, x+x_{0}+tv_{0}, \, v+v_{0}) \quad 
	\text{for any}~z_{0}:=(t_0,x_0,v_0), z:=(t,x,v) \in \mathbb{R}^{1+2n}\,,
\end{equation}
and the usual kinetic scaling~$\delta_r:  \er^{1+2n} \mapsto \er^{1+2n}$ defined by
\begin{equation} \label{def:dil} 
	\delta_r (z):=(r^{2s}t,\,r^{1+2s}x,\, r v) \quad \text{for any}~r > 0. 
\end{equation} 
So that the inverse of each element $z_0 \in \mathbb{R}^{1+2n}$ is defined and 
\begin{align*}
	z_{0}^{-1} \circ z = ( t- t_{0}, x - x_{0} - (t-t_{0}) v_{0}, v-v_{0}) \qquad \text{for any }~z \in \mathbb{R}^{1+2n}.
\end{align*}
For any~$r>0$, we denote by~${Q}_r$ a cylinder centered in the origin of radius~$r$; that is,
$$
{Q}_r \equiv {Q}_r({0}) \ :=  \  (-r^{2s},0]\times B_{r^{1+2s}}(0)\times B_r(0)\,.
$$ 
For every~$z_{0} \in \er^{1+2n}$ and for every~$r>0$, the {\it slanted} cylinder~${Q}_r(z_{0})$ is defined as follows,
\begin{eqnarray}
	\label{Q-classico}
	&& {Q}_r(z_{0})\!\!\! \ :=
	\big \{ z \in \er^{1+2n}: \, -r^{2s} < t -  t_{0} \leq  0 ,  \ | x - x_{0} - (t - t_{0}) v_{0}| <r^{1 + 2s},\ \snr{ v- v_{0} } < r \big\}. 
\end{eqnarray}
Moreover,  we also denote by
$$
U_r(z_0) \ := \
\left\{(t,x)\in\er^{1+n}:\,-r^{2s}<t-t_0\leq0,\ 
|x-x_0-(t-t_0)v_0|<r^{1+2s}\right\}\,,
$$
the projection of~$Q_r(z_0)$ onto the drift variables.

Furthermore, we denote with~$N_{\tt s}$ the {\it homogeneous dimension} related to~\eqref{def:dil}, i.e.
\begin{equation}\label{def:homo-dim}
	N_{\tt s}:=  n(2+2s)+2s,
\end{equation}	
which encodes the scaling properties of the underlying kinetic scalings, and letting~$\Om := (t_1,t_2) \times \Om_x \times \Om_v \subset \er^{1+2n}$ be a bounded domain in~$\er^{1+2n}$, we work under the following assumption on the source term in~\eqref{problema}: 
\begin{equation}\label{HP:source}
h\in L^{q}_{\rm loc}(\Om) \quad \text{with}~q>\frac{N_{\tt s}}{2s}.
\end{equation}
On the other hand, for the estimates involving the nonlocal tail we use the homogeneous dimension of the drift variables $(t,x)$, i.e.
\begin{equation}\label{def:drift-dim}
N_{\tt d} \ := \ n(1+2s) + 2s.
\end{equation}
The nonlocal tail was introduced in the nonlinear fractional elliptic setting~\cite{DKP14,DKP16} to quantify the influence of exterior values on local estimates.
For the kinetic equation, the ``{\it  nonlocal tail} of a function $f$ {centred in~$v_{0}\in\Omega_v\subset\er^n$ of  dif{f}usion radius~$r$''} is given by
\begin{equation}		 
	\label{def_tail_intro}
	\tail(f;B_r(v_{0}))  
	\, :=\,  r^{2s}\int_{ \ern \setminus B_r(v_{0})}\frac{\snr{f(v)}}{\snr{v_{0}-v}^{n+2s}} \dv.
\end{equation}
Thus, it is integrated in velocity and its summability is measured only in the drift variables, so that~$N_{\tt d}/(2s)$ is the scale-invariant threshold.

\medskip
Our first result is a~$\delta$-interpolative~$L^2$-$L^\infty$ estimate with tail for weak subsolutions to~\eqref{problema}, with possibly unbounded source terms~$h$.
The parameter~$0<\delta\leq 1$ makes the coefficient of the $L^{p}_{t,x}$-tail norm in~\eqref{kinetic_special} arbitrarily small, at the cost of increasing the coefficient of the local quadratic norm. This freedom will also be used in the proof of the strong Harnack inequality.
Our first main result reads as follows.\footnote{
Let us specify that, given any open set~$\Oc \subset \er^{1+2n}$, with positive Lebsegue measure $\snr{\Oc}>0$ we denote with~$\nra{g}_{L^p(\Oc)}$, its mean $L^p$-norm. With abuse of notation, we will adopt also this definition whenever~$\Oc \subset \er^{1+n}$ .
}

\begin{theorem}\label{thm_bdd}
Let~$\Omega:=(t_1,t_2)\times\Om_x\times\Om_v\subset\er^{1+2n}$ be a domain and~$s\in(0,1)$. Assume that~$f\in\Wc(\Om)$ is a weak subsolution to~\eqref{problema}, under~\eqref{HP:source}, such that, for every~$B_r(v_0) \Subset\Om_v$,
\begin{equation}\label{HP:tail}
\tail(f_+;B_r(v_0))\in L^{p}_{\rm loc}((t_1,t_2)\times\Omega_x) \quad \text{with}~p \ > \  \frac{N_{\tt d}}{2s}\,.
\end{equation}
Then, for every~$r>0$ such that~$Q_r(z_0)\Subset\Om$ and every~$\delta\in(0,1]$, it holds
\begin{equation}\label{kinetic_special}
\begin{aligned}
\operatorname*{ess\,sup}_{Q_\frac{r}{2}(z_0)}f
& \ \leq \  \frac{c}{\delta^\beta}
\nra{f_+}_{L^2(Q_r(z_0))}
+c\,r^{2s}\nra{h_+}_{L^{q}(Q_r(z_0))}\\
&\qquad+\delta\,\nra{\tail(f_+;B_\frac{r}{2}(v_0))}_{L^{p}(U_r(z_0))}\,,
\end{aligned}
\end{equation}
where~$c \equiv c (n,s,\Lambda,p,q,|v_0|)>0$, while~$\beta\equiv \beta(n,s,p,q)>0$.
\end{theorem}
The proof of Theorem~\ref{thm_bdd} uses a De Giorgi iteration in which the tail and the general source are treated in different spaces. We begin with a two-level Caccioppoli estimate: the test function is truncated at the upper level, while the exterior contribution is estimated at the lower one.
The gap between the levels allows us to relate the velocity superlevel sets to the local quadratic mass.

We then localize the truncation and compare it with a solution of the fractional Kolmogorov equation.
The cutoff terms and the terms involving the rough kernel are estimated in~$L^2_{t,x}H^{-s}_v$, without differentiating the kernel pointwise. The exterior contribution is bounded by a function of~$(t,x)$ multiplied by the indicator of a velocity superlevel set.
For each fixed pair of drift variables, Chebyshev's inequality controls the quadratic velocity norm of this indicator by the lower truncation and the level gap. This allow us to use mixed-norm estimates governed by the drift dimension~$N_{\tt d}$.

The general source~$h=h(t,x,v)$ is handled separately. Under the full-dimensional condition~$q>N_{\tt s}/(2s)$, the Kolomgorov potential of its localized contribution is bounded and can be absorbed by shifting the truncation level before applying the mixed estimates. The bound for this shift decreases with the measure of the superlevel set, which allows it to fit within the successive level gaps; the source contribution to the energy estimate is controlled in the same iteration.

To close the argument, we use two distinct gains of integrability for the fractional Kolmogorov semigroup.
The first improves integrability in the drift variables while retaining a quadratic norm in velocity, and controls the tail term in the energy estimate. The second improves integrability in all kinetic variables and produces a positive power of the superlevel measure in the quadratic estimate.
The strict inequality in the tail assumption is used here to make this power positive for the exterior term.
The coupled iteration then makes the local quadratic mass tend to zero, without imposing the full-dimensional source condition on the tail.

\medskip
The tail is already integrated in velocity and bounds the exterior source uniformly in that variable.
Compared with the general source condition~\eqref{HP:source}, its critical summability exponent is therefore lower by~$n/(2s)$ and is precisely the one in~\eqref{HP:tail}.

For each exponent outside the requested summability range, we are able to construct globally nonnegative, time-dependent solutions with uniformly bounded~$L^{p}_{t,x}$-tail norm for which estimate~\eqref{kinetic_special} fails.
Thus the local bound used to derive the strong Harnack inequality in Theorem~\ref{thm_strong} cannot be extended below this threshold. Our construction adapts ideas from~\cite{KW24c} and uses the kernel of~$(-\Delta_v)^s$, so the obstruction is not due to variable or rough coefficients.

\begin{theorem}\label{thm_sharp}
Let~$n\geq1$,~$s\in(0,1)$ and~$\mathcal Q :=(-2,1/2)\times B_2\times B_2$.
For any~$m \in \en$, there exists a family of nonnegative weak solutions~$f_m\in\Wc(\mathcal Q)$
to~\eqref{problema} in~$\mathcal Q$, with~$h=0$,
such that, for every exponent~$p$ satisfying
\[
1 \ \leq \  p \ \leq \ \frac{N_{\tt d}}{2s}\,,
\]
the tail of~$f_m$ satisfy~$\tail(f_m;B_{\frac{1}{2}})\in L^{p}(U_1)$
and the ratio
\[
\frac{f_m(0)}
{\nra{f_m}_{L^2(Q_1)}
+\nra{\tail(f_m;B_{\frac{1}{2}})}_{L^{p}(U_1)}}
\]
blows up as~$m  \to +\infty$.
\end{theorem}
This result is consistent with the stationary counterexample of Ka\ss mann and Weidner in~{\rm\cite{KW24c}}. In the stationary problem there is no time integration, and the relevant homogeneous dimension is~$n(1+2s)$; hence the corresponding obstruction is governed by the exponent~$n(1+2s)/(2s)$. In the time-dependent problem considered here, the tail is measured in~$L^{p}$ with respect to~$(t,x)$, and the time scale naturally contributes the additional homogeneous dimension~$2s$, thus giving~$N_{\tt d}=2s+n(1+2s)$. 
We prescribe nonnegative data outside the velocity domain and concentrate it in time and space at the kinetic scales.
Adding a correction obtained from the Dirichlet semigroup cancels the interior source produced by these data, so the resulting function solves the homogeneous equation.
A lower bound for the mass retained by the Dirichlet kernel gives an interior lower bound, while contraction estimates control the local quadratic norm and the scaling of the exterior data controls the tail.
Below the stated threshold, both norms decay faster than the interior value.

\vspace{2mm}

The local estimate in Theorem~\ref{thm_bdd} gives boundedness on interior cylinders without assuming that~$f$ is globally bounded.
For globally nonnegative solutions in the homogeneous case~$h=0$, we combine it with the weak Harnack inequality of Imbert and Silvestre~\cite{IS20}.
This gives the following strong Harnack inequality under the tail summability condition stated above.
\begin{theorem}\label{thm_strong}
Let~$s\in(0,1)$. 
Then, there exists~$\bar{R}>1$, depending only on~$n$ and~$s$, such that letting~$\widetilde{Q}_{\bar{R}}:= [-2,0]\times B_{\bar{R}^{1+2s}}\times B_{\bar{R}}$, if~$f\in \Wc(\widetilde{Q}_{\bar{R}})$ is a globally nonnegative weak solution to~\eqref{problema} with~$h=0$ in~$\widetilde{Q}_{\bar{R}}$, satisfying~\eqref{HP:tail} for any~$B_r(v_0) \subset B_{\bar{R}}$,  there exists~$r_0\in(0,1/2)$, depending only on~$n$ and~$s$, such that
\begin{equation}\label{strong}
\operatorname*{ess\,sup}_{Q^-_{r_0}}f
\ \leq \  c\operatorname*{ess\,inf}_{Q^+_{r_0}}f
+c\nra{\tail(f;B_{r_0})}_{L^{p}(U_{2r_0}(-1+(2r_0)^{2s},0,0))},
\end{equation}
where~$c\equiv c(n,s,p,\Lambda)>0$ and
\begin{equation}\label{def_slantedpm}
\begin{split}
Q^-_{r_0}:=(-1+(2r_0)^{2s},0,0)\circ Q_{r_0},\\
Q^+_{r_0}:=(-r_0^{2s},0]\times B_{r_0^{1+2s}}\times B_{r_0}.
\end{split}
\end{equation}
\end{theorem} 
\begin{figure}[htp!]
	\centering
	\includegraphics[scale=0.23]{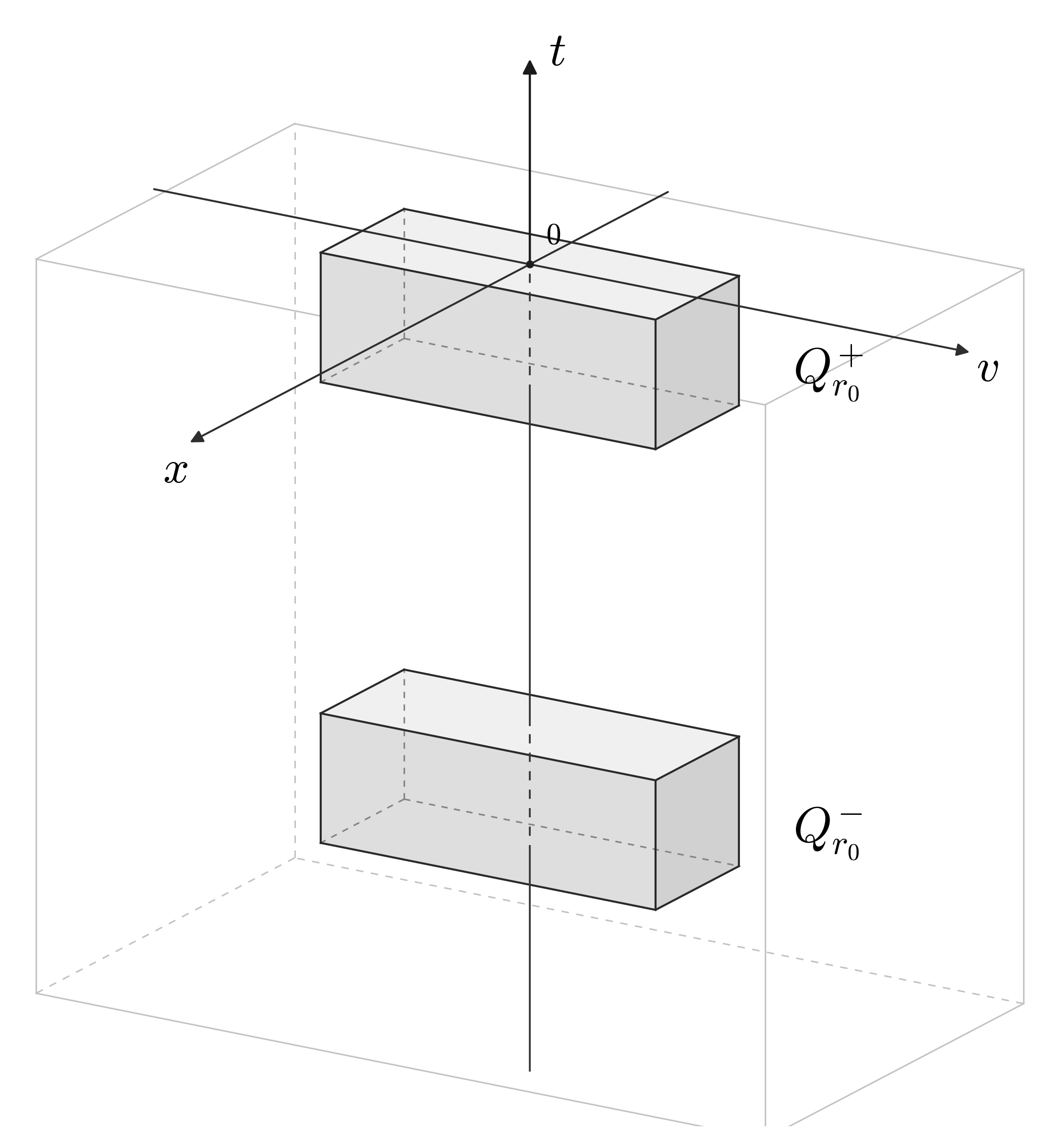}
	\caption{\small The geometry of the Harnack inequalities for kinetic equations.}
	\label{fig_slanted}
\end{figure}

To prove this estimate, we apply the local bound on small slanted cylinders covering the backward cylinder.
The tail at each centre is split into an annular contribution inside a fixed larger velocity ball and the contribution from outside that ball.
The annular term is controlled by a supremum on a larger cylinder.
After interpolation between the quadratic norm and the lower-power norm in the weak Harnack inequality, we choose the tail coefficient small enough to absorb the larger supremum in a radius iteration.
The weak Harnack inequality then bounds the remaining lower-power norm by the infimum on the forward cylinder.

Finally, notice that for fractional elliptic equations the tail~${\rm Tail}(f_-)$ measures the contribution of the exterior negative part. In the kinetic estimate, a tail term remains even for globally nonnegative solutions: it accounts for the interaction of the exterior values with transport in the interior.

\begin{rem}
The definition of weak solution to~\eqref{problema} does not require a finite~$L^{p}$-norm of the velocity tail. However, for nonnegative solutions, a uniform upper bound on the velocity mass plainly implies this requirement: for every fixed velocity ball, it bounds the tail uniformly in the drift variables.
Such mass bounds are among the assumptions used in models for the non-cutoff Boltzmann equation and related kinetic equations; see, for instance, {\rm\!\!\cite[Theorems~1.1-1.2]{Sil16},~\cite[Formula~(1.4)]{Mou18},~\cite[Formula~(1.3)]{GIMV19},~\cite[Formula~(1.3)]{IMS20},~\cite[Section~1.4-Assumption~(H)]{IS20a},~\cite[Section~1.3]{IS20},~\cite[Assumption~1.1]{IS22},~\cite[Formula~(1.2)]{OS23}, and~\cite[Formula~(1.9)]{FRW24}}.
\end{rem}

\subsection*{Outline of the paper} In Section~\ref{sec_preliminaries} we fix the notation, introduce the relevant function spaces and the weak formulation, and prove the mixed-norm estimates for the fractional Kolmogorov semigroup and its parametrix.
In Section~\ref{sec_gain} we prove the two-level Caccioppoli estimate and the gain of integrability in Theorem~\ref{thm:gain}, then close the coupled De Giorgi iteration leading to Theorem~\ref{thm_bdd}.
Section~\ref{sec_harnack} is devoted to the proof of Theorem~\ref{thm_strong}.
Finally, in Section~\ref{sec_psharp} we construct the counterexample that proves the sharpness of the tail summability threshold.

\medskip\section{Preliminaries}\label{sec_preliminaries}
In this section we fix the notation and briefly recall the necessary information on the underlying functional framework required to deal with~\eqref{problema}.

\subsection{Notation and basic tools}
We denote with~$c$ a positive universal constant greater than one, which may change from line to line. For the sake of readability, dependencies of the constants will be often omitted within the chains of estimates, therefore stated after the estimate.
 Relevant dependencies on parameters will be emphasized by using parentheses.
 
 For any~$\Oc\subset \er^n$ we denote with~$\mathbbm{1}_{\Oc}$ the indicator function of~$\Oc$.
As customary, for any~$r>0$ and any~$y_{0} \in \er^n$ we denote by~$
B_r(y_{0})$\,,
the open ball with radius~$r$ and center~$y_{0}$, when not important, or clear from the context, we shall omit denoting the center as follows:~$B_r \equiv B_r(y_0)$.

For any measurable function~$g$, we define  the positive 
and negative part of~$g$ as~$g_\pm(y):=\max\{\pm g(y),0\}$. 

For any function~$f$ smooth enough, we will denote with $\hat{f}$ its Fourier transform.
More precisely, for~$f\in\mathcal S(\er^m)$ we use the convention
\[
\widehat f(\xi) \ = \ \int_{\er^m}e^{-iy\cdot\xi}f(y)\d y
\quad \text{and} \quad 
f(y)\ = \ \frac1{(2\pi)^m}\int_{\er^m}e^{iy\cdot\xi}\widehat f(\xi)\d\xi.
\]

For any~$p \in [1,\infty)$ we define the weak Lebesgue space~$L^{p}_w(\Oc)$ as the space of measurable functions~$g: \Oc \to \er$ such that the following quasi-norm is finite
\[
\|g\|_{L^{p}_w(\Oc)} \ := \ 
\sup_{\lambda>0}\lambda\snr{\{z \in \Oc: |g(z)|>\lambda\}}^\frac{1}{p}.
\] 

We now recall some useful integral inequalities that will be used frequently throughout the section. We begin with the Minkowski's Inequality; see, e.~\!g.,~\cite[Proposition~1.3]{BCD11}.

Let~$(\mathbb{X}_1,\mu_1)$ and~$(\mathbb{X}_2,\mu_2)$ be two measure spaces and~$f$ being a nonnegative measurable function on~$\mathbb{X}_1\times \mathbb{X}_2$. For any~$1 \leq p \leq q \leq \infty$, we have
\begin{equation}\label{eq:mink-ineq}
 \left\|\|f(\cdot,x_2)\|_{L^p(\mathbb{X}_1,\mu_1)}\right\|_{L^q(\mathbb{X}_2,\mu_2)} \ \leq \ \left\|\|f(x_1,\cdot)\|_{L^q(\mathbb{X}_2,\mu_2)}\right\|_{L^p(\mathbb{X}_1,\mu_1)}.
\end{equation}

We also recall the Hardy--Littlewood--Sobolev inequality in~$\er^n$; see, e.~\!g.,~\cite[Theorem~1.7]{BCD11}.
Let~$\alpha \in (0,n)$ and~$p,q \in (1,\infty)$ such that
\[
\frac{1}{p}+\frac{\alpha}{n}  \ = \ 1+ \frac{1}{q}.
\]
Then, there exists~$c>0$ such that
\begin{equation}\label{eq:hls-ineq}
\||\cdot|^{-\alpha}\ast f \|_{L^q(\ern)} \ \leq \ c\,\|f\|_{L^p(\ern)}.
\end{equation}

For the sake of notation, throughout the paper we will denote with~$\|\cdot\|_{L^q_x}$,~$\|\cdot\|_{L^q_v}$ respectively, the $L^q$-norm in the spatial, or velocity variable. Similarly, we adopt the following notation for anisotropic Lebesgue norms
\[
\|g\|_{L^{q}_xL^{p}_v} \ : = \ \left(\int_{\ern}\left(\int_{\ern}|g(x,v)|^p\dv\right)^\frac{q}{p}\dx\right)^\frac{1}{q}.
\]

\medskip
Now we introduce some basic facts about fractional Sobolev spaces and the appropriate weak formulation of~\eqref{problema}. For further details, we refer the interested reader to~\cite{DPV12,AN26}.

 For~$s\in (0,1)$ we denote with~$H^s(\Oc)$ the classical fractional Sobolev space
\[
H^s(\Oc) \ := \ \left\{f \in L^2(\Oc) \; : \; \left[ f \right]_{H^s(\Oc)}  < +\infty \right\}\,,
\]
where the seminorm~$\left[ f \right]_{H^s(\Oc)} $ is the usual one defined via the Gagliardo kernel
\[
\left[ f \right]_{H^s(\Oc)}  \ :=  \ \left(\iint_{\Oc\times\Oc} \frac{\snr{f(v)-f(w)}^2}{\snr{v-w}^{n+2s}} \dv\dw\right)^\frac{1}{2}\,,
\]
and where we endowed~$H^s$ with the usual norm
\[ 
\norma{f}_{H^s(\Oc)} := \norma{f}_{L^2(\Oc)} + \left[ f \right]_{H^s(\Oc)} .
\]

A function~$ f$ belongs to $H_{\rm loc}^{s}(\Oc)$ if~$f\in H^s(\Oc')$ whenever~$\Oc' \Subset \Oc$. We will denote with~$H^{-s}(\ern)$ the dual of~$H^s(\ern)$ and denote with~$\langle \cdot,\cdot \rangle_{H^{-s},H^s}$ the usual duality pairing between~$H^{-s}$ and~$H^s$. Let us remark that, via Riesz--Fr\'echet's Representation Theorem, for any~$f \in H^{-s}(\ern)$ there exist two functions~$H_0$,~$H_1 \in L^2(\ern)$ such that
\begin{equation} \label{eq:repr}
f = H_1 + (-\Delta_v)^\frac{s}{2}H_{0}\quad \text{and} \quad \|H_{0}\|_{L^2(\ern)} + \|H_1\|_{L^2(\ern)} \approx \|f\|_{H^{-s}(\ern)}.
\end{equation}

For any~$f \in H^s(\ern)$ we define~$-\Lc f$ as an element of~$H^{-s}(\ern)$ that acts on~$\phi \in H^s(\ern)$ via
\[
\Ec(f,\phi):=\langle -\Lc f ,\phi \rangle_{H^{-s},H^s} = \frac{1}{2} \iint_{\ern\times\ern}(f(v)-f(w))(\phi(v)-\phi(w)) K(v,w)\dv\dw.
\]
Here and below the dependence of~$K$ on~$(t,x)$ is suppressed whenever these variables are fixed.

For~$s\in (0,1)$, consider the following  tail space
\[
L^1_{2s}(\ern):= \left\{g \in L^1_{\rm{loc}}(\ern)\, : \,  \int_{\ern}\frac{|g(v)|}{(1+\snr{v})^{n+2s}}\dv < \infty \right\},
\]
as used in~\cite{KKP16}; see Section~2~in~\cite{KKP17} for related properties.
 
Given~$\Om:=    (t_1,t_2)\times \Om_x \times \Om_v \subset \er^{1+2n}$,  
we denote by~$\Wc$ the natural functions space to which weak solutions to~\eqref{problema} belong to, and defined as
$$
 \begin{aligned}
 \Wc(\Omega) & :=  \Big\{ f \in L^2_{\rm loc}((t_1, t_2)\times \Om_x ;\, H^s_{\rm loc}(\Om_v))\cap L^1_{\rm loc}((t_1, t_2)\times \Om_x ;L^1_{2s}(\ern)) 
 	\\*
    & \qquad\qquad\qquad\qquad \qquad: (\partial_t +v \cdot \nabla_x)f \in L^2_{\rm loc} ((t_1,t_2)\times \Om_x;\,  H^{-s}(\ern )) \Big \}.
 \end{aligned} 
 $$

We are now in the position to recall the definition of weak sub- and supersolution.

\begin{defn}\label{def_weak_sol}
A function~$f \in \Wc(\Om)$ is a weak subsolution {\rm (}resp.,
supersolution{\rm)}
to~\eqref{problema} in~$\Om$ with
$h \in L^2(\Om)$ if
\[
\int_{t_1}^{t_2}\int_{\Om_x}\Ec(f,\phi)\dx\dt
+
\int_{t_1}^{t_2}\int_{\Om_x}
\langle \partial_t f+v\cdot\nabla_xf,\phi
\rangle_{H^{-s},H^s}\dx\dt
 \ \stackrel{\text{(resp., }\geq\text{)}}{\leq}  \
\iiint_{\Om}h\phi\dz\,,
\]
for any nonnegative~$\phi\in
L^2((t_1,t_2)\times\Om_x;H^s(\ern)) \cap L^\infty(\Om)$
such that~$
{\rm supp}\,\phi\Subset \Om$.

A function~$f\in\Wc(\Om)$ is a weak solution to~\eqref{problema}
if it is both a weak sub- and supersolution.
\end{defn}

\subsection{The fundamental solution}\label{sec_fondamentale}
As already remarked in the Introduction, in order to establish the gain of integrability in Theorem~\ref{thm:gain} we need to invoke the hypoelliptic nature of equation~\eqref{problema}; see, e.~\!g.,~\cite{HPZ21}. Indeed, we shall rely on the regularizing properties of the fundamental solution of the fractional Kolmogorov equation, whose expression -- with pole in the origin -- is given by
\begin{equation}\label{eq:fundamental}
		P_t(x,v) := 
	\begin{cases}
		\displaystyle \frac{1}{t^{n+\frac{n}{s}}} \mathdutchcal{P}\left(\frac{x}{t^{1+\frac{1}{2s}}}, \frac{v}{t^\frac{1}{2s}}\right) \quad & \text{if}~t>0\,,\\*[0.7ex]
		0 \quad & \text{if}~t \leq 0\,,
	\end{cases}
\end{equation}
where the kernel~$\mathdutchcal{P}$ is defined in Fourier variables by
\[
\widehat{\mathdutchcal{P}}(\xi,\eta) := e^{-\int_0^1 \snr{\eta + \tau\xi}^{2s}\d\tau}.
\]

Formula~\eqref{eq:fundamental} and the basic properties below are contained in~\cite[Proposition~2.1]{IS20}, after translating its Fourier convention. Sharp two-sided estimates in arbitrary dimension are proved in~\cite[Theorem~1.1]{HZ24}, while the one-dimensional case is also treated by Fourier methods in~\cite{Gru24}.
    
The following proposition summarizes the main properties of the fundamental solution~$P(\cdot)$ we are going to use throughout the manuscript.
	\begin{prop}[Proposition~2.1 in~\cite{IS20}]\label{prop:fund}
    The function~$P$ and $\mathdutchcal{P}$ have the following properties:
    \begin{enumerate}[(i)]
    \item the function $\mathdutchcal{P}$ is $C^\infty$ and decays polynomially at infinity. Moreover, $\mathdutchcal{P}$ and all its derivatives are integrable in $\er^{2n}$;
    \item both the functions $\mathdutchcal{P}$ and $P$ are nonnegative;
    \item for any $t>0$,~$\|P_t\|_{L^1(\er^{2n})}=1$;
    \item for every~$f_{0} \in L^2(\er^{2n})$ and every~$h \in L^2((0,\infty)\times\er^{2n})$ the function
    \begin{eqnarray*}
    	f(t,x,v) & = & \iint_{\ern\times\ern} P_t(x-y-tw,v-w)f_{0}(y,w)\dw\d y \\
    	&& + \int_0^t \iint_{\ern\times\ern} P_{t-\tau}(x-y-(t-\tau)w,v-w)h(\tau,y,w) \dw \d y \d \tau\,,
    \end{eqnarray*}
    is the unique solution to the Cauchy problem
   $$
    \begin{cases}
     \partial_t f + v \cdot \nabla_x f + (-\Delta_v)^sf = h & \quad \text{in}~(0,\infty) \times \er^{2n},\\
     f(0,\cdot) = f_{0} & \quad \text{in}~\{0\} \times \er^{2n}.
    \end{cases}
    $$
    \end{enumerate}
	\end{prop}

    Associated with the fundamental solution~$P_t$ in~\eqref{eq:fundamental} we define with~$\mathscr{S}_t$ the fractional Kolmogorov semigroup, which, for sufficiently regular function~$f$, is given by
    \begin{equation}\label{sharp:m1:eq_semigroup}
     \mathscr{S}_t f(\tau,x,v) \ := \ \iint_{\ern\times\ern}P_t(x-y-tw,v-w)f(\tau,y,w)\,\d y\,\dw\,,
     \end{equation}
     and
      we also let
    \begin{equation}\label{sharp:m1:eq_potential}
     \mathscr{P}f(t,x,v) \ := \ \int_{-\infty}^t\big[\mathscr{S}_{t-\tau}f(\tau,\cdot,\cdot)\big](x,v)\,\mathrm d\tau.
     \end{equation}

  \begin{rem}\label{rem:P-conv}
  We remark that~\eqref{sharp:m1:eq_potential} can be written as~$\mathscr{P}f = f\ast_{\rm kin}P$ where we defined the group convolution according to the group law in~\eqref{def:action} on~$(\er^{1+2n},\circ)$ as
  \[
 (f\ast_{\rm kin}g)(z) \ := \ \iiint_{\er^{1+2n}}f(z')g((z')^{-1}\circ z)\dz'\,,
 \]
 for a couple of sufficiently regular functions~$f$,~$g$ on~$(\er^{1+2n},\circ)$.
  \end{rem}

In the subsequent part of this section we will prove some estimate for the functional~$\mathscr{P}$ in~\eqref{sharp:m1:eq_potential}.
We begin recalling the following estimate for the fundamental solution~$P_t$.
\begin{lemma}\label{sharp:m1:lem_kernel_scale}
For every~$t>0$ and $q_v,q_x \in [1,\infty]$, there exists a constant~$c \equiv c(n,s,q_x,q_v)>0$ such that
\begin{equation}\label{sharp:m1:eq_kernel_mixed}
\|P_t\|_{L^{q_v}_vL^{q_x}_x} \ \leq \ c\,t^{-\frac{n(1+2s)}{2s}(1-\frac{1}{q_x})-\frac{n}{2s}(1-\frac{1}{q_v})}.
\end{equation}
\end{lemma}

\begin{proof}
By~\cite[Theorem~1.1 and Lemma~1.5]{HZ24}, with~$\alpha=2s$ and~$\beta=n+2s$, the profile satisfies
\[
\int_{\ern}\mathdutchcal P(x,v)\,\mathrm dx\leq C(1+|v|)^{-\beta},
\qquad
\sup_x\mathdutchcal P(x,v)\leq C(1+|v|)^{-1-\beta}.
\]
Interpolation in~$x$ gives
\[
\|\mathdutchcal P(\cdot,v)\|_{L^{q_x}_x}
\leq C(1+|v|)^{-1-\beta+1/q_x}.
\]
The right-hand side belongs to~$L^{q_v}_v$ for every~$q_v\in[1,\infty]$, since~$\beta=n+2s>n$. Finally, the scaling in~\eqref{eq:fundamental}, with the changes~$x=t^{1+1/(2s)}X$ and~$v=t^{1/(2s)}V$, gives
\[
\|P_t\|_{L^{q_v}_vL^{q_x}_x}
=t^{-n-n/s+\frac{n(1+2s)}{2sq_x}+\frac{n}{2sq_v}}
\|\mathdutchcal P\|_{L^{q_v}_vL^{q_x}_x},
\]
which is exactly~\eqref{sharp:m1:eq_kernel_mixed}.
\end{proof}

The following heat kernel estimates for the kinetic semigroup~$\mathscr{S}_t$ hold true.

\begin{prop}
\label{sharp:m1:prop_fixed}
Let $t>0$. Then, for any~$1\leq p\leq2\leq q\leq\infty$, there exists a constant~$c\equiv c(n,s,p,q)>0$ such that 
\begin{equation}
\label{sharp:m1:eq_fixed_vector}
\|\mathscr{S}_t f\|_{L^q_xL^2_v} \ \leq \ c\, t^{-\frac{n(1+2s)}{2s}(\frac{1}{p}-\frac{1}{q})}\|f\|_{L^p_xL^2_v}\,,
\end{equation}
and
\begin{equation}
\label{sharp:m1:eq_fixed_scalar}
\|\mathscr{S}_t f\|_{L^q_{x,v}} \ \leq \ c\,t^{-\frac{n(1+2s)}{2s}\left(\frac{1}{p}-\frac{1}{q}\right)
-
\frac{n}{2s}\left(\frac{1}{2}-\frac{1}{q}\right)}
\|f\|_{L^p_xL^2_v}.
\end{equation}
\end{prop}

\begin{proof}
Let us fix~$1\leq p_x,p_v,q_x,q_v,\ell_x,\ell_v\leq\infty$ such that
\begin{equation}
\label{sharp:m1:eq_young_rel}
1+\frac{1}{q_x}=\frac{1}{p_x}+\frac{1}{\ell_x},
\quad \text{and} \quad 
1+\frac{1}{q_v}=\frac{1}{p_v}+\frac{1}{\ell_v}.
\end{equation}
Then, for any fixed~$v$, applying Young's Inequality in~$x$ gives
\[
\|(K*G)(\cdot,v)\|_{L^{q_x}_x}
\leq
\int_{\ern}\|K(\cdot,v-w)\|_{L^{\ell_x}_x}\|G(\cdot,w)\|_{L^{p_x}_x}\,\dw\,,
\]
so that a final application of Young's Inequality in~$v$ yields
\begin{equation}
\label{eq_mixed_young}
\|K*G\|_{L^{q_v}_vL^{q_x}_x}
\leq
\|K\|_{L^{\ell_v}_vL^{\ell_x}_x}
\|G\|_{L^{p_v}_vL^{p_x}_x}.
\end{equation}

Moreover, let us note that, for~$1\leq p\leq2\leq q\leq\infty$, the integral Minkowski's Inequality gives
\begin{equation}
\label{eq_minkowski_orders}
\|G\|_{L^2_vL^p_x}
\leq\|G\|_{L^p_xL^2_v}
\quad \text{and} \quad
\|G\|_{L^q_xL^2_v}
\leq\|G\|_{L^2_vL^q_x}.
\end{equation}

Letting~$\mathscr{T}_tf(y,w) := f(y-tw,w)$, by changing~$y$ into~$y+tw$ in~\eqref{sharp:m1:eq_semigroup}, we obtain
\begin{equation}
\label{sharp:m1:eq_ordinary_conv}
\mathscr{S}_t f=P_t*\mathscr{T}_tf.
\end{equation}

In order to prove~\eqref{sharp:m1:eq_fixed_vector}, choose~$\ell_v=1$ and~$\ell_x$ such that
\[
1-\frac{1}{\ell_x}=\frac{1}{p}-\frac{1}{q}.
\]

Combining Lemma~\ref{sharp:m1:lem_kernel_scale} with~\eqref{eq_mixed_young} and~\eqref{eq_minkowski_orders}, we get
\[
\|\mathscr{S}_tf\|_{L^q_xL^2_v} \ \leq \ \|\mathscr{S}_tf\|_{L^2_vL^q_x}  \ \leq \ \|P_t\|_{L^1_vL^{\ell_x}_x}
\|\mathscr{T}_tf\|_{L^2_vL^p_x}\ \leq \ c\, t^{-\frac{n(1+2s)}{2s}(\frac{1}{p}-\frac{1}{q})}
\|f\|_{L^p_xL^2_v},
\]
which proves~\eqref{sharp:m1:eq_fixed_vector}.

For~\eqref{sharp:m1:eq_fixed_scalar}, we proceed in a similar fashion by choosing~$\ell_x$ and~$\ell_v$ such that
\[
1-\frac{1}{\ell_x}=\frac{1}{p}-\frac{1}{q}
\quad \text{and} \quad 1-\frac{1}{\ell_v} = \frac{1}{2}-\frac{1}{q}\,,
\]
 so that Lemma~\ref{sharp:m1:lem_kernel_scale},~\eqref{eq_mixed_young} and~\eqref{eq_minkowski_orders} give
\[
\|\mathscr{S}_tf\|_{L^q_{x,v}} \ \leq \ \|P_t\|_{L^{\ell_v}_vL^{\ell_x}_x}\|\mathscr{T}_tf\|_{L^2_vL^p_x} \ \leq \ 
c\,t^{-\frac{n(1+2s)}{2s}\left(\frac{1}{p}-\frac{1}{q}\right)
-
\frac{n}{2s}\left(\frac{1}{2}-\frac{1}{q}\right)}\|f\|_{L^p_xL^2_v}.
\]
\end{proof}

Next, we prove the estimate which will be applied to the error between the rough operator and $(-\Delta_v)^s$.

\begin{lemma}
\label{sharp:m1:lem_derivative}
For every $t>0$ and $h\in L^2(\er^{2n})$,
\begin{equation}
\label{sharp:m1:eq_derivative_L2}
\|\mathscr{S}_t(-\Delta_v)^\frac{s}{2}h\|_{L^2_{x,v}}
\ \leq  \ \frac{c\,\|h\|_{L^2_{x,v}}}{\sqrt{t}}\,,
\end{equation}
where~$c \equiv c(n,s)>0$.
\end{lemma}

\begin{proof}
With no loss of generality, upon reasoning by density, take~$h$ in the Schwartz class.
By applying the Fourier transform, the semigroup has the representation
\begin{equation}\label{sharp:m1:eq_fourier_semigroup}
\widehat{\mathscr{S}_th}(\xi,\eta)
=
\exp\left(-\int_0^t|\eta+\rho\xi|^{2s}\,\mathrm d\rho\right)
\widehat h(\xi,\eta+t\xi).
\end{equation}

We shall use the elementary coercivity estimate
\begin{equation}
\label{sharp:m1:eq_coercivity}
\int_0^1|\eta+\rho\xi|^{2s}\,\mathrm d\rho
\geq c_s\left(|\eta|^{2s}+|\xi|^{2s}\right).
\end{equation}

By Plancherel's theorem and~\eqref{sharp:m1:eq_fourier_semigroup}, followed by the change of variables $\eta \mapsto \zeta=\eta+t\xi$, we have
\begin{align*}
\|\mathscr{S}_t(-\Delta_v)^\frac{s}{2}h\|_2^2
&= \frac{1}{(2\pi)^{2n}}
\iint_{\ern\times\ern}
\exp\left(-2\int_0^t|\zeta-r\xi|^{2s}\,\mathrm dr\right)
|\zeta|^{2s}|\widehat h(\xi,\zeta)|^2
\,\mathrm d\xi\,\mathrm d\zeta.
\end{align*}
By scaling~\eqref{sharp:m1:eq_coercivity},
$$
\int_0^t|\zeta-r\xi|^{2s}\,\mathrm dr
\geq c_st\left(|\zeta|^{2s}+t^{2s}|\xi|^{2s}\right)
\geq c_st|\zeta|^{2s}.
$$
Therefore, since~$\er_+ \ni r \mapsto r^{2s} e^{-2c_st r^{2s}}$ is bounded by~$c/t$, for some~$c=c(s)>0$, we obtain
~\eqref{sharp:m1:eq_derivative_L2} 
\end{proof}

Now, we introduce a set of exponents, whose definition is based on~$N_{\rm d}$ and~$N_s$ introduced in~\eqref{def:drift-dim} and~\eqref{def:homo-dim}, respectively. First of all, set
\begin{equation}\label{sharp:eq_drift_exponents}
q_{\rm d}:=\frac{2N_{\rm d}}{N_{\rm d}-2s},
\qquad
a_{\rm d}:=\frac{2N_{\rm d}}{N_{\rm d}+2s}.
\end{equation}
Notice that
\begin{equation}\label{sharp:m1:eq_drift_rel}
\frac1{a_{\rm d}}-\frac1{q_{\rm d}}=\frac{2s}{N_{\rm d}},
\qquad
\frac12-\frac1{q_{\rm d}}=\frac{s}{N_{\rm d}}.
\end{equation}
Recall that~$p_*$ is defined in~\eqref{HP:tail}.
For $p\in[p_*,\infty)$, we define
\begin{equation}\label{sharp:m1:eq_ap_qp}
a_p:=\frac{2p}{p+2},
\qquad
\frac1{q_p}:=\frac12-\frac{2s-N_{\rm d}/p}{N_s}.
\end{equation}
Equivalently,
\begin{equation}
\label{sharp:m1:eq_qp_formula}
q_p=\frac{2N_s}{N_s-4s+2N_{\rm d}/p}.
\end{equation}
Then $1<a_p<2\leq q_p<\infty$, and $q_p>2$ if and only if $p>p_*$.

\begin{prop}\label{sharp:m1:prop_fixed_negative}
For every~$t>0$ and~$h\in L^2(\er^{2n})$,
\begin{equation}
\label{sharp:m1:eq_fixed_negative}
\|\mathscr{S}_t(-\Delta_v)^\frac{s}{2}h\|_{L^{q_{\rm d}}_xL^2_v} \ \leq \ \frac{c\,\|h\|_{L^2_{x,v}}}{t^\frac{N_{\tt d}-s}{N_{\tt d}}}\,,
\end{equation}
where~$N_{\tt d}$ is defined in~\eqref{def:drift-dim} and the constant~$c\equiv c(n,s)>0$.
\end{prop}

\begin{proof}
Using the semigroup property, we write
$$
\mathscr{S}_t(-\Delta_v)^\frac{s}{2}h
=
\mathscr{S}_{t/2}\bigl(\mathscr{S}_{t/2}(-\Delta_v)^\frac{s}{2}h\bigr).
$$
We apply Proposition~\ref{sharp:m1:prop_fixed} with $p=2$ and $q=q_{\rm d}$, and then Lemma~\ref{sharp:m1:lem_derivative}. Since $1/2-1/q_{\rm d}=s/N_{\rm d}$, we obtain
\begin{align*}
\|\mathscr{S}_t(-\Delta_v)^\frac{s}{2}h\|_{L^{q_{\rm d}}_xL^2_v}
&\leq
ct^{-\frac{n(1+2s)}{2s}(\frac12-\frac1{q_{\rm d}})}
\|\mathscr{S}_{t/2}(-\Delta_v)^\frac{s}{2}h\|_2
\\
&\leq
ct^{-\frac{n(1+2s)}{2N_{\rm d}}-\frac12}\|h\|_2
=ct^{-(1-s/N_{\rm d})}\|h\|_2.
\end{align*}
\end{proof}

Now, our main result reads as follows.
\begin{theorem}
\label{sharp:thm_mixed_parametrix}
The following estimates hold.
\begin{enumerate}[(i)]
\item For every $F\in L^{a_{\rm d}}(\er^{1+n};L^2(\ern))$,
\begin{equation}
\label{sharp:m1:eq_K1}
\|\mathscr{P}F\|_{L^{q_{\rm d}}(\er^{1+n};L^2(\ern))}
\leq c\|F\|_{L^{a_{\rm d}}(\er^{1+n};L^2(\ern))}.
\end{equation}
\item For every $H\in L^2(\er^{1+2n})$,
\begin{equation}
\label{sharp:m1:eq_K2}
\|\mathscr P((-\Delta_v)^\frac{s}{2}H)\|_{L^{q_{\rm d}}(\er^{1+n};L^2(\ern))}
\leq c\|H\|_{L^2(\er^{1+2n})}.
\end{equation}
\item Let $p\in[p_*,\infty)$. For every $F\in L^{a_p}(\er^{1+n};L^2(\ern))$,
\begin{equation}
\label{sharp:m1:eq_K3}
\|\mathscr P F\|_{L^{q_p}(\er^{1+2n})}
\leq c_p\|F\|_{L^{a_p}(\er^{1+n};L^2(\ern))}.
\end{equation}
\end{enumerate}
Here, the constant $c$ depends only on $n$ and $s$, while $c_p$ depends only on $n$, $s$, and $p$.
\end{theorem}
\begin{proof}
We first take smooth compactly supported sources. 

\noindent
{\it (i)} Applying Proposition~\ref{sharp:m1:prop_fixed} with $p=a_{\rm d}$ and $q=q_{\rm d}$, we obtain
\begin{align*}
\|\mathscr PF(t,\cdot,\cdot)\|_{L^{q_{\rm d}}_xL^2_v}
&\leq
c\int_{-\infty}^t
(t-\tau)^{-\frac{n(1+2s)}{2s}(\frac1{a_{\rm d}}-\frac1{q_{\rm d}})}
\|F(\tau,\cdot,\cdot)\|_{L^{a_{\rm d}}_xL^2_v}
\,\mathrm d\tau
\\
&=
c\int_{-\infty}^t
(t-\tau)^{\frac{2s}{N_{\rm d}}-1}
\|F(\tau,\cdot,\cdot)\|_{L^{a_{\rm d}}_xL^2_v}
\,\mathrm d\tau,
\end{align*}
where we used~\eqref{sharp:m1:eq_drift_rel}.

Then, taking the $L^{q_d}$-norm in time on both sides of the inequality and applying \eqref{eq:hls-ineq}, with $\alpha=1 - \frac{2s}{N_{\rm d}}$ and $q=q_d$ gives~\eqref{sharp:m1:eq_K1}.

\noindent
{\it (ii)} For every $t>0$, Proposition~\ref{sharp:m1:prop_fixed_negative} yields
$$
\|\mathscr P((-\Delta_v)^\frac{s}{2}H)(t,\cdot,\cdot)\|_{L^{q_{\rm d}}_xL^2_v}
\leq
c\int_{-\infty}^t
(t-\tau)^{s/N_{\rm d}-1}
\|H(\tau,\cdot,\cdot)\|_2\,\mathrm d\tau.
$$
Since $1/q_{\rm d}=1/2-s/N_{\rm d}$, $0<s/N_{\rm d}<1$ and~$1<2<q_{\rm d}<\infty$, taking the $L^{q_d}$-norm in time on both sides of the inequality and applying \eqref{eq:hls-ineq} 
with $\alpha=1-s/N_{\rm d}$ and $q=q_d$ proves~\eqref{sharp:m1:eq_K2}.

\noindent
{\it (iii)} To prove~\eqref{sharp:m1:eq_K3} we set
\begin{equation}
\label{sharp:m1:eq_Ap_Bp}
A_p:=\frac1{a_p}-\frac1{q_p},
\qquad
B_p:=\frac12-\frac1{q_p}.
\end{equation}
By~\eqref{sharp:m1:eq_ap_qp}, $0<A_p<1$ and $B_p\geq0$, with $B_p>0$ if and only if $p>p_*$ (with $p^*$ defined in ~\eqref{HP:tail}). 
More explicitly,
\[
B_p=\frac{2s-N_{\rm d}/p}{N_s},
\qquad
A_p=\frac{2s}{N_s}+\frac{n}{N_sp}.
\]
Thus~$B_p\geq0$, while~$0<A_p\leq A_{p_*}=2s/N_{\rm d}<1$; in particular,~$a_p<q_p$. We also have
\begin{equation}
\label{sharp:m1:eq_scaling_identity}
N_{\rm d} A_p+nB_p=2s.
\end{equation}

Also, identity~\eqref{sharp:m1:eq_scaling_identity} gives
\begin{equation}
\label{sharp:m1:eq_beta_identity}
\frac{n(1+2s)}{2s}A_p+\frac n{2s}B_p
=1-A_p.
\end{equation}
Applying~\eqref{sharp:m1:eq_fixed_scalar} with $p=a_p$ and $q=q_p$, we obtain
\begin{align*}
\|\mathscr PF(t,\cdot,\cdot)\|_{L^{q_p}_{x,v}}
&\leq
c_p\int_{-\infty}^t
(t-\tau)^{A_p-1}
\|F(\tau,\cdot,\cdot)\|_{L^{a_p}_xL^2_v}
\,\mathrm d\tau.
\end{align*}
Since $1/q_p=1/a_p-A_p$, taking the $L^{q_p}$ norm on both sides of the previous inequality, we apply \eqref{eq:hls-ineq} with $\alpha = 1-A_p$
and $q=q_p$. This proves~\eqref{sharp:m1:eq_K3}.
The three estimates extend to the stated source spaces by density.
\end{proof}

The following consequence will be used for the $L^2_v+H^{-s}_v$ decomposition of the rough-kernel error.

\begin{cor}
\label{sharp:m1:cor_rough}
Let $H_0,H_1\in L^2(\er^{1+2n})$, and assume that the projection of $\operatorname{supp}H_1$ onto the $(t,x)$ variables is contained in a measurable set $E\subset\er^{1+n}$ with $|E|<\infty$. Then
\begin{equation}
\label{sharp:m1:eq_rough_pair}
\|\mathscr P(H_1+(-\Delta_v)^\frac{s}{2}H_0)\|_{L^{q_{\rm d}}_{t,x}L^2_v}
\leq c\left(|E|^{s/N_{\rm d}}\|H_1\|_{L^2_{t,x,v}}+\|H_0\|_{L^2_{t,x,v}}\right).
\end{equation}
No support assumption in the velocity variable is required.
\end{cor}

\begin{proof}
Since $1/a_{\rm d}-1/2=s/N_{\rm d}$, H\"older's Inequality in $(t,x)$ gives
$$
\|H_1\|_{L^{a_{\rm d}}_{t,x}L^2_v}
\leq |E|^{s/N_{\rm d}}\|H_1\|_2.
$$
We apply Theorem~\ref{sharp:thm_mixed_parametrix}~\textup{(i)} to $H_1$, Theorem~\ref{sharp:thm_mixed_parametrix}~\textup{(ii)} to $H_0$, and use the triangle inequality.
\end{proof}

\begin{cor}
\label{sharp:m1:cor_slab}
Let $I\subset\er$ be an interval. If a source in Theorem~\ref{sharp:thm_mixed_parametrix} is supported in $I\times\er^{2n}$, then the corresponding estimate holds with the same constant after restricting the output to $I\times\er^{2n}$, or to any measurable subset thereof.
\end{cor}

\begin{proof}
Extend the source by zero outside $I\times\er^{2n}$, apply Theorem~\ref{sharp:thm_mixed_parametrix}, and restrict the output.
\end{proof}

\begin{cor}
\label{sharp:m1:cor_local}
Let $R>0$ and $z_0\in\er^{1+2n}$. Extend every source below by zero outside $Q_R(z_0)$.
\begin{enumerate}[(i)]
\item If~$F$ is supported in~$Q_R(z_0)$, then
\[
	\left(\biint_{U_R(z_0)}\nra{\mathscr P F}_{L^2(B_R(v_0))}^{q_{\tt d}}\dx\dt\right)^\frac{1}{q_{\tt d}} \ 
	\leq \  cR^{2s}\left(\biint_{U_R(z_0)}\nra{F}_{L^2(B_R(v_0))}^{a_{\tt d }}\dx\dt\right)^\frac{1}{a_{\tt d}} .
\]
\item If~$H$ is supported in~$Q_R(z_0)$, then
\begin{align}
\label{sharp:m1:eq_local_K2}
	\left(\biint_{U_R(z_0)}\nra{\mathscr P((-\Delta_v)^\frac{s}{2}H}_{L^2(B_R(v_0))}^{q_{\tt d}}\dx\dt\right)^\frac{1}{q_{\tt d}} \leq cR^s  \left(\biint_{U_R(z_0)}\nra{ H}_{L^2(B_R(v_0))}^{2}\dx\dt\right)^\frac{1}{2} 
\end{align}
\item Let $p\in[p_*,\infty)$. If $F$ is supported in $Q_R(z_0)$, then
\begin{align}
\label{sharp:m1:eq_local_K3}	
	\nra{\mathscr P F}_{L^{q_{p}}} \leq c_pR^{2s} \left(\biint_{U_R(z_0)}\nra{F}_{L^2(B_R(v_0))}^{a_p}\dx\dt\right)^\frac{1}{a_p}.
\end{align}
\end{enumerate}
The constants are independent of $R$ and $z_0$.
\end{cor}

\begin{cor}\label{sharp:m1:cor_rough_local}
Let $R>0$ and $z_0\in\er^{1+2n}$. Assume that $H_0,H_1\in L^2(\er^{1+2n})$ vanish for $(t,x)\notin U_R(z_0)$, with no support assumption in velocity. Then
\begin{align}
\label{sharp:m1:eq_rough_local}
&\left(\biint_{U_R(z_0)} \| \mathscr P(H_1+(-\Delta_v)^\frac{s}{2}H_0) \|_{L^2(\mathbb{R}^{n})}^{q_d} \dx\dt\right)^\frac{1}{q_d}  \\ \nonumber
&\leq c \left( R^{2s}  \left(\biint_{U_R(z_0)} \| H_1 \|_{L^2(\mathbb{R}^{n})}^{2} \dx\dt\right)^\frac{1}{2}  
+ R^s \left(\biint_{U_R(z_0)} \| H_0 \|_{L^2(\mathbb{R}^{n})}^{2} \dx\dt\right)^\frac{1}{2}    \right).
\end{align}
\end{cor}

\begin{proof}[Proof of Corollary~\ref{sharp:m1:cor_rough_local}]
Apply Corollary~\ref{sharp:m1:cor_rough} with~$E=U_R(z_0)$ and restrict the output to~$U_R(z_0)\times\ern_v$. Since~$|U_R|=c_nR^{N_{\rm d}}$ and~$1/2-1/q_{\rm d}=s/N_{\rm d}$,
\[
|U_R|^{-1/q_{\rm d}}|U_R|^{s/N_{\rm d}}\|H_1\|_2
=c R^{2s}|U_R|^{-1/2}\|H_1\|_2,
\]
while
\[
|U_R|^{-1/q_{\rm d}}\|H_0\|_2
=c R^s|U_R|^{-1/2}\|H_0\|_2.
\]
These are precisely the two terms in~\eqref{sharp:m1:eq_rough_local}.
\end{proof}

\begin{lemma}\label{lemma:fundamental}
Let~$1<a<N_s/(2s)$ and define~$q$ by
\begin{equation}\label{sharp:eq_full_hls_relation}
\frac1q=\frac1a-\frac{2s}{N_s}.
\end{equation}
Then
\begin{equation}\label{sharp:eq_full_hls}
\|\mathscr P F\|_{L^q(\er^{1+2n})}
\leq c\|F\|_{L^a(\er^{1+2n})},
\end{equation}
where~$c\equiv c(n,s,a)>0$.
\end{lemma}

\begin{proof}
The kernel~$P$ is homogeneous of degree~$-(N_s-2s)$ with respect to the dilations in~\eqref{def:dil}, namely
$$
P(\delta_r (z))=r^{-(N_s-2s)}P(z).
$$
We first notice that~$|\{P>1\}|<+\infty$. Indeed, Proposition~\ref{prop:fund}~\textup{(i)} gives~$\|\mathdutchcal P\|_{L^\infty}<+\infty$, and hence~$P(t,x,v)>1$ implies that~$t$ belongs to a bounded interval. For every such~$t$, Proposition~\ref{prop:fund}~\textup{(ii)} and Chebyshev's Inequality give
$$
|\{(x,v):P(t,x,v)>1\}|\leq\int_{\er^{2n}}P(t,x,v)\,\d x\,\d v=1.
$$
Consequently, by homogeneity and~$|\delta_r(E)|=r^{N_s}|E|$,
$$
|\{P>\lambda\}|
=\lambda^{-\frac{N_s}{N_s-2s}}|\{P>1\}|
\qquad\text{for every }\lambda>0.
$$
Thus
$$
P\in L^{\vartheta}_w(\er^{1+2n})
\quad \text{with}~\vartheta:=\frac{N_s}{N_s-2s}.
$$
 The weak Young inequality on homogeneous groups, see~\cite[Proposition~1.10]{Fol75}, yields
$$
\|F\ast_{\rm kin}P\|_{L^q}
\leq c\|P\|_{L^{\vartheta}_w}\|F\|_{L^a}
$$
whenever~$1+1/q=1/a+1/\vartheta$. This relation is precisely~\eqref{sharp:eq_full_hls_relation}, and the conclusion follows.
\end{proof}

\begin{lemma}[Proposition 2.2 in~\cite{IS20}]\label{lemma:kolmogorov-young}
Let~$I\subset\er$ be an interval of length at most one and set~$q_0:=2N_s/(N_s-2s)$.
If~$H_0,H_1\in L^2(I\times\er^{2n})$, extended by zero outside~$I\times\er^{2n}$, then
\[
\begin{aligned}
&\|\mathscr{P}H_1\|_{L^{q_0}(I\times\er^{2n})}
+
\|\mathscr{P}(( -\Delta_v)^\frac{s}{2}H_0)\|_{L^{q_0}(I\times\er^{2n})}
\\
&\qquad\leq
c\left(\|H_1\|_{L^2(I\times\er^{2n})}
+\|H_0\|_{L^2(I\times\er^{2n})}\right).
\end{aligned}
\]
where the constant~$c\equiv c(n,s)>0$.
\end{lemma}

\subsection{Auxiliary results}
 We recall the weak Harnack inequality in the form needed below. It is possible to verify that the assumptions of Theorem~1.6 in~\cite{IS20} are satisfied in our framework. So that, the statement adapted to our case reads as follows.
 
   \begin{theorem}\label{thm:weak} Let~$s \in (0,1)$.Then, there exist~$r_{0}\in(0,1)$,~$\bar{R} >1$,~$\zeta>0$ and~$c>0$ such that, letting~$\widetilde{Q}_{\bar{R}} := [-2,0]\times B_{\bar{R}^{1+2s}} \times B_{\bar{R}}$, if~$f \in \Wc(\widetilde{Q}_{\bar{R}})$ is a globally nonnegative weak supersolution to
      \[
   (\partial_t + v\cdot \nabla_x) f = \Lc f \quad {in}~\widetilde{Q}_{\bar{R}}\,,
      \]
   then,
   \[
       \|f\|_{L^\zeta(Q^-_{r_0})} \leq c\, \operatorname*{ess\,inf}_{Q^+_{r_{0}}}f\,,
   \]
   where~$c$ and~$\zeta$ depend only on~$n$,~$s$ and~$\Lambda$, while~$r_{0}$ and~$\bar{R}$ depend only on~$n$ and~$s$
    and
    \begin{equation*}
		\begin{split}
			& {Q}^-_{r_{0}} := (-1,-1+r_{0}^{2s}] \times B_{r_{0}^{1+2s}}\times B_{r_{0}}\\*[0.5ex] \quad \text{and} & \quad  {Q}^+_{r_{0}} := (-r_{0}^{2s},0] \times B_{r_{0}^{1+2s}}\times B_{r_{0}}.
		\end{split}
	\end{equation*}
   \end{theorem}

\medskip
We recall a covering property of the slanted cylinders in~\eqref{Q-classico}.  For a similar result in the classic kinetic framework we refer to~\cite[Lemma~4.2]{PP04}; see also~\cite{AP20}.

\begin{lemma}\label{lemma:cov}
There exist two universal constants~$c_* \equiv c_*(s) \in (0,1)$ and~$\beta \equiv \beta(s) \geq 1$ such that, for any~$1/2 \leq \rho < r \leq 1$ and any~$z_{0} \in \er^{1+2n}$, it holds
\begin{equation}\label{eq:inclusion}
{Q}_{(c_*(r - \rho))^\beta}(z_1) \subset {Q}_{r}(z_{0}) \qquad \forall z_1 \in {Q}_{\rho}(z_{0}).
\end{equation}
\end{lemma}
   \begin{proof}
Set~$d:=r-\rho$,~$c_*(s):=\frac12(1\wedge2s)$ and~$\beta(s):=1\vee(2s)^{-1}$, and put~$\delta:=(c_*d)^\beta$. Then~$\delta\leq d$ and
\[
\delta^{2s}\leq
\begin{cases}
d^{2s},&s\geq1/2,\\
sd,&s<1/2.
\end{cases}
\]
On the other hand,
\[
r^{2s}-\rho^{2s}\geq
\begin{cases}
d^{2s},&s\geq1/2,\\
2sd,&s<1/2,
\end{cases}
\]
so the time and velocity components of every~$z_2\in Q_\delta(z_1)$, with~$z_1\in Q_\rho(z_0)$, lie in the corresponding components of~$Q_r(z_0)$. For the spatial component, the group law gives
\[
|x_2-x_0-(t_2-t_0)v_0|
\leq \delta^{1+2s}+\rho^{1+2s}+\delta^{2s}\rho
\leq \rho^{1+2s}+\delta^{2s}r.
\]
If~$s\geq1/2$, then~$\delta^{2s}r\leq d^{2s}r\leq r^{1+2s}-\rho^{1+2s}$. If~$s<1/2$, then~$\delta^{2s}r=sd r\leq r^{1+2s}-\rho^{1+2s}$, because~$\rho\geq1/2$. This proves~\eqref{eq:inclusion}.
 \end{proof}

    We conclude by stating a classical iteration argument which will turn out to be useful in establishing our main results.

   \begin{lemma}[Lemma~4.11 in~\cite{Coz17}]\label{lemma_giusti}
Let~$\Uppsi:[{\rho}, {r}] \rightarrow [0,+\infty)$ be a 
bounded function,~$\eps \in (0,1)$,~$A_1,$ $A_2$, $A_3$,~$\beta_1$,~$\beta_2 \geq 0$ and~${\rho}, {r} > 0$. Assume that
\[
\Uppsi(\sigma') \leq \eps \Uppsi(\sigma) + {{A}_1}{(\sigma-\sigma')^{-\beta_1}}+ {{A}_2}{(\sigma-\sigma')^{-\beta_2}} + A_3,
\]
holds whenever~${\rho} \leq  \sigma' < \sigma \leq  {r}.$ Then, 
\[
\Uppsi({\rho}) \leq  c{A_1}{({r} - {\rho})^{-\beta_1}} +  c{A_2}{({r} - {\rho})^{-\beta_2}} + cA_3,
\]
where~$c\equiv c(\eps, \beta_1,\beta_2) > 0.$
\end{lemma}

\section{Local boundedness estimates}\label{sec_gain}

This section is devoted to the proof of local boundedness for subsolutions to~\eqref{problema}, as stated in~Theorem~\ref{thm_bdd}. We begin collecting some fundamental observations.

\begin{rem}\label{rem:blow-up}
We first work on cylinders centered at the origin, with radii between~$[1/2,1]$. Indeed,  see, e.~\!g.,~\cite[Lemma~5.1]{Sto19} the function~$\tilde{f}(z):=f(z_{0}\circ \delta_r(z))$ satisfies~\eqref{problema} with adjusted source term and kernel. In particular the ellipticity class of~$K$ does not change. However,  a dependence on the norm~$|v_0|$ will appear in the constant multiplying~$\nra{f}_{L^2(Q_r)}$ in~\eqref{kinetic_special} as consequence of the sensibility of the transport operator~$\partial_t +v\cdot \nabla_x$ to the Galilean translation in~\eqref{def:action}. 
\end{rem}

\begin{rem}\label{rem:exponent}
We just recall here the three exponents used in the forthcoming anisotropy estimates:
\begin{equation}\label{eq:exponent-crt}
p_*:=1+\frac{n(1+2s)}{2s},
\quad
q_{\tt d}
=
\frac{2N_{\tt d}}{N_{\tt d}-2s},
\quad \text{and} \quad 
a_{\tt d}
=
\frac{2N_{\tt d}}
{N_{\tt d}+2s}\,,
\end{equation}
where~$N_{\tt d}$ is the homogeneous dimension in the drift variables defined in~\eqref{def:drift-dim}. 
In particular, the following relations hold true:
\begin{equation}\label{eq:exponent-identities}
\frac{1}{q_{\tt d}} = \frac{1}{2}-\frac{s}{N_{\tt d}},
\quad \text{and} \quad \frac{1}{a_{\tt d}} = \frac{1}{2}+\frac{s}{N_{\tt d}}.
\end{equation}
Moreover, given a general exponent~$p > p_\ast$ we let
\begin{equation}\label{eq:bar-p-T}
    \bar{p} := \min\{p,2p_\ast\}
\end{equation}

\end{rem}

Let us also record that, for~$0\leq\ell$, since~$(f-\ell)_+\leq f_+$,~$B_\frac{1}{2}\subset B_r$ and~$U_r\subset U_1$ whenever~$1/2\leq r\leq1$, it holds~$\tail((f-\ell)_+;B_r)\leq 4\,\tail(f_+;B_\frac{1}{2})$ pointwise and, by~$\bar{p} \leq p$ and the monotonicity of normalized norms,~$\nra{\tail((f-\ell)_+;B_r)}_{L^{\bar{p}}(U_r)}\leq c(n,s)\nra{\tail(f_+;B_\frac{1}{2})}_{L^{p}(U_1)}$. This will be used in the rest of the this section  
 to pass from the tail quantities produced by the proofs to those appearing in the statements.

\subsection{The two-level energy estimate}
Firstly, we prove a suitable (two-level) Caccioppoli-type estimates with tail for subsolutions to~\eqref{problema}.
The upper level~$\kk$ is used in the test function, while the lower level~$\ell<\kk$ is used only to estimate the exterior interaction.

\begin{lemma}\label{lemma:energy-est}
	Let~$Q_1 \subset \er^{1+2n}$  and~$s \in (0,1)$.  Assume that~$f \in \Wc(Q_1)$ is a weak  subsolution to~\eqref{problema} in~${Q}_1$. Let~$p_\ast$ and~$q_{\tt d}$ be given in~\eqref{eq:exponent-crt}. 
	Then, for any~$ \frac{1}{2} \leq {\rho} < {r} \leq 1$, any~$\beta >0$, any~$p>p_\ast$ and any~$0 \leq \ell<\kk$, it holds 
	\begin{equation}\label{caccioppoli-def}
    \begin{aligned}
     &\biiint_{Q_\rho}(f-\kk)_+^2\d z 
     +\biiint_{Q_\rho}\int_{B_\rho}
      \frac{\left((f(t,x,v)-\kk)_+-(f(t,x,w)-\kk)_+\right)^2}{|v-w|^{n+2s}}\dw\dz \\
      &
\quad\quad \ \leq \ \frac{\beta}{(r-\rho)^{2(1+2s)}}\left(\biint_{U_r}\nra{(f-\ell)_+}_{L^2(B_r)}^{q_{\tt d}}\dx\dt\right)^\frac{2}{q_{\tt d}} +\frac{c}{(r-\rho)^{2(1+2s)}}\biiint_{Q_r}(f-\ell)_+^2\dz \\
     &\quad\qquad 
       +c\left(\frac{\nra{\tail(f_+;B_\frac{1}{2})}_{L^{{p}}(U_1)}}{(r-\rho)^{n+2s}(\kk-\ell)}\right)^{\frac{{p}}{{p}-p_*} }\biiint_{Q_r}(f-\ell)_+^2\dz\\
    &\quad\qquad  +c\biiint_{Q_{\frac{\rho+r}{2}}}h_+(f-\kk)_+\dz\,,
\end{aligned}
\end{equation}
	where~$c \equiv c(n,s,\Lambda,\beta,p)>0$.
\end{lemma}
\begin{proof} 
	For the sake of readability we denote with~$f_{\kk} := (f-\kk)_+$ and with~$f_\ell := (f-\ell)_+$. Let $\sigma := (r+\rho)/2$ and consider a cut-off function~$\phi \in C^\infty_c(B_{\sigma^{1+2s}}\times B_\sigma)$ such that~$0 \leq \phi \leq {1}$,~$\phi =1 $ on~$B_{\rho^{1+2s}}\times B_\rho$ and~$|\nabla_v \phi| + |v\cdot \nabla_x\phi| \leq c/(r-\rho)^{1+2s}$, for some universal constant~$c>0$.
	Consider in the weak formulation the test function~$\phi^2 f_\kk$ so that, up to mollification (see for instance~\cite[Section~2]{GM22}), for a.~\!e.~$t \in (-{{r}^{2s}},0]$, we have that
	\begin{equation}\label{eq:cacc-est1}
    \begin{aligned}
		\iint_{B_{{r}^{1+2s}}\times B_{{r}}} h (\phi^2f_\kk) \dx\dv\  & \ \geq \ \iint_{B_{{{r}}^{1+2s}}\times B_{{r}}} (\partial_t f+v\cdot \nabla_x f)( \phi^2f_\kk )\dx\dv\\
		& \quad \ +\int_{B_{{{r}}^{1+2s}}} \Ec(f,\phi^2f_\kk)\dx \  =: \   I_1 + I_2.
    \end{aligned}
	\end{equation}
	
	We start by considering~$I_1$. Using the fact  that~$\partial_t\phi=0$, and that, by~\cite[Formula (A8)]{IS20}, 
	\begin{equation}\label{eq:transport-truncation}
		(\partial_t + v\cdot \nabla_x)\,f_\kk=
		(\partial_t + v\cdot \nabla_x)\,(f-\kk)\mathbbm{1}_{\{f>\kk\}}\,,
    \end{equation}
    we have that
	\begin{equation}\label{mancava} 
    \begin{aligned}
		I_1 & \ \geq  \ \frac{1}{2}\frac{\d}{\dt} \iint_{B_{{{r}}^{1+2s}}\times B_{{r}}} \snr{(\phi f_\kk){(t,x,v)}}^2  \dx\dv - \iint_{B_{{{r}}^{1+2s}}\times B_{{r}}}f_\kk^2\snr{v\cdot \nabla_x \phi}\dx\dv \\
		& \ \geq \  \frac{1}{2}\frac{\d}{\dt} \iint_{B_{{{r}}^{1+2s}}\times B_{{r}}} \snr{ (\phi f_\kk){(t,x,v)}}^2  \dx\dv - \frac{c}{({r}-{\rho})^{1+2s}}\iint_{B_{{{r}}^{1+2s}}\times B_{{r}}} f_\kk^2\dx\dv.
    \end{aligned}
	\end{equation}
	
	Now we focus on the term~$I_2$. By symmetry of the involved kernel~$K$, we can simply restrict on the subcase when~$f(t,x,v) \geq f(t,x,w)$ \big(up to exchange the roles of~$v$ and~$w$\big), and thus we simply use that
	\[
    \begin{aligned}
			&   \big(f(t,x,v)-f(t,x,w)\big)\big((\phi^2f_\kk)(t,x,v)-(\phi^2f_\kk)(t,x,w)\big)\\*
			& \quad \ = \   \big((f(t,x,v)-\kk)-(f(t,x,w)-\kk)\big)\big((\phi^2f_\kk)(t,x,v)-(\phi^2f_\kk)(t,x,w)\big)\\
			& \quad \geq \
			\begin{cases}
				\big(f_\kk(t,x,v)-f_\kk(t,x,w)\big)\\
                \qquad\times \big((\phi^2f_\kk)(t,x,v)-(\phi^2f_\kk)(t,x,w)\big)  & \text{if}~f(t,x,v) \geq f(t,x,w) > \kk,\\*
				\big (f(t,x,v) - \kk\big)(\phi^2f_\kk)(t,x,v) & \text{if}~f(t,x,v) >\kk \geq f(t,x,w),\\*
				0 & \text{if}~\kk \geq f(t,x,v) \geq f(t,x,w),
			\end{cases}     \\*
			& \quad \geq \ 
			\big((\phi f_\kk)(t,x,v)-(\phi f_\kk)(t,x,w)\big) ^2 -f_\kk(t,x,v)f_\kk(t,x,w)\big(\phi(x,v)-\phi(x,w)\big)^2\,,     
    \end{aligned}
	\]
	which, in turn, yields the bound
	\begin{equation}\label{eq:energy-0}
		\begin{aligned}
		I_2
		& \ \geq  \   \int_{B_{{r}^{1+2s}}} \Ec(\phi f_\kk,\phi f_\kk)\dx \\
		&  \quad \ -  \Lambda\int_{B_{{r}^{1+2s}}}\iint_{B_r \times B_r}f_\kk(t,x,v)f_\kk(t,x,w)\frac{(\phi(x,v)-\phi(x,w))^2}{\snr{v-w}^{n+2s}}\dw\dv\dx\\
		& \quad \ -  \Lambda\int_{B_{{r}^{1+2s}}}\iint_{B_r \times (\ern \setminus B_r)}\frac{(\phi^2f_\kk)(t,x,v)f_\kk(t,x,w)}{\snr{v-w}^{n+2s}}\dw\dv\dx\\
		& \ \geq \ \int_{B_{{r}^{1+2s}}} \Ec(\phi f_\kk,\phi f_\kk)\dx -  \frac{c}{({r}-{\rho})^{2(1+2s)}}\iint_{B_{r^{1+2s}}\times B_r}f_\kk^2 \dv \dx\\
		& \quad \ -  c\int_{B_{{r}^{1+2s}}}\iint_{B_r \times (\ern \setminus B_r)}\frac{(\phi^2 f_\kk)(t,x,v)f_\kk(t,x,w)}{\snr{v-w}^{n+2s}}\dw\dv\dx\,,
        \end{aligned}
	\end{equation}
where in the last line we have used that, 
 up to exchanging the roles of~$v$ and~$w$, we can assume that~$f_\kk(t,x,v) \geq f_\kk(t,x,w)$, so that
	\begin{equation}\label{eq:singular-int}
	\begin{aligned}
        & \iint_{B_r \times B_r}f_\kk(t,x,v)f_\kk(t,x,w)\frac{(\phi(x,v)-\phi(x,w))^2}{\snr{v-w}^{n+2s}} \dw\dv \\*
		&\qquad \qquad \qquad \qquad\qquad \leq \  \frac{c}{({r}-{\rho})^{2(1+2s)}}\int_{B_{{r}}}f_\kk^2(t,x,v) \left(\int_{B_{2{r}}(v)}\frac{\dw}{|v-w|^{n-2(1-s)}}\right)\dv\\*
		&\qquad \qquad \qquad \qquad\qquad\leq \ \frac{c}{({r}-{\rho})^{2(1+2s)}}\int_{B_{{r}}}f_\kk^2 \dv.
    \end{aligned}
	\end{equation}
	
	Combining~\eqref{mancava} and~\eqref{eq:energy-0} with~\eqref{eq:cacc-est1} yields
    \[
    \begin{aligned}
   & \frac{1}{2}\frac{\d}{\dt} \iint_{B_{r^{1+2s}}\times B_r} \snr{(\phi f_\kk)(t)}^2\dx\dv + \int_{B_r^{1+2s}} \Ec(\phi f_\kk,\phi f_\kk) \dx \\*
   &  \qquad\quad \leq \ c\iint_{B_{r^{1+2s}}}\iint_{B_r \times (\ern \setminus B_r)} \frac{(\phi^2f_\kk)(t,x,v)f_\kk(t,x,w)}{\snr{v-w}^{n+2s}} \dw\dv\dx \\*
	&\qquad\qquad \ + \frac{c}{({r}-{\rho})^{2(1+2s)}}\iint_{B_{r^{1+2s}}\times B_r} f_\kk^2\dv\dx  + c\,\iint_{B_{\sigma^{1+2s}}\times B_{\sigma}} h_+ f_\kk\dx\dv.\notag
    \end{aligned}
	\]
	
	After integration in time and manipulations precisely as in~\cite[Lemma 2.2]{Sto19}, via the ellipticity of the kernel in~$Q_\rho$, and since~$\ell<\kk$, we eventually arrive at
	\begin{equation}\label{eq:cacc-final-1}
    \begin{aligned}
	& \biiint_{Q_\rho} f_\kk^2\dz   +\biiint_{Q_\rho}\int_{B_\rho}        \frac{\left(f_\kk(t,x,v)-f_\kk(t,x,w)\right)^2}{|v-w|^{n+2s}}\dw\dz\\*
	& \qquad\quad  \leq  \ 
		 c\biiint_{Q_r}\int_{\ern \setminus B_r} \frac{(\phi^2f_\kk)(t,x,v)f_\kk(t,x,w)}{\snr{v-w}^{n+2s}} \dw\d z \\*
	& \qquad\qquad + \frac{c}{({r}-{\rho})^{2(1+2s)}}\biiint_{Q_r}f_\ell^2\dz+ c\,\biiint_{Q_{\frac{r + \rho}{2}}} h_+ f_\kk\dz.
    \end{aligned}
	\end{equation}
	
	In order to conclude, it remains to estimate the nonlocal contribution and the contribution given by the source term. Let us first note that, for any~$v \in B_r \cap \textup{supp}(\phi) \subset B_{(r+\rho)/2}$ and any~$w \in \ern \setminus B_r$, it  holds
	\[
	\frac{|w|}{|v-w|}\, \leq \,1 + \frac{|v|}{||w|-|v||} \,\leq\, 1 + \frac{r+\rho}{r-\rho}\, \leq \, \frac{c}{r-\rho}.
	\]
    Moreover, on~$\{f>\kk\}$ we also have that
    \[
    0 \ \leq \ f-\kk \leq \frac{(f-\ell)^2}{\kk-\ell} \quad \forall \ell< \kk.
    \]
    Therefore, we can estimate the tail term on the right-hand side of~\eqref{eq:cacc-final-1} as follows
	\[
    \begin{aligned}
		& \biiint_{Q_r}\int_{\ern \setminus B_r} \frac{(\phi^2f_\kk)(t,x,v)f_\kk(t,x,w)}{\snr{v-w}^{n+2s}} \dw\dz\\
        & \qquad\quad \leq  \ \frac{c}{(r-\rho)^{n+2s}} \biiint_{Q_r} f_\kk(t,x,v) \tail(f_\kk; B_r)\dz \\
		& \qquad\quad \leq \ \frac{c}{(r-\rho)^{n+2s}(\kk-\ell)} \biint_{U_r} \tail(f_\kk; B_r) \nra{f_\ell}_{L^2(B_r)}^2\dx\dt.
    \end{aligned}
	\]

    Now, considering~${p}> p_\ast$, with~$p_\ast$ as in~\eqref{eq:exponent-crt}, we apply H\"older's Inequality in~$(t,x)$  with exponents~${p}$ and~${p}'$, obtaining
    \begin{equation}\label{sharp:cacc:tail-holder}
     \biint_{U_r}\tail(f_\ell;B_r)\nra{f_\ell}_{L^2(B_r)}^2\dx\dt
   \ \leq \ \nra{\tail(f_\ell;B_r)}_{L^{{p}}(U_r)}\left(\,\biint_{U_r}\nra{f_\ell}_{L^2(B_r)}^{2{p}'}\dx\dt\right)^\frac{1}{{p}'}.
   \end{equation}

 Define now
  \begin{equation}\label{def:theta-t}
  \theta \ := \ \frac{p_*}{{p}} \ =  \ \frac{N_{\tt d}}{2s{p}} \in (0,1).
  \end{equation}

  Using~\eqref{eq:exponent-identities}, we compute
 \begin{equation}\label{eq:interpolation}
 \frac{1}{2{p}'} \ = \ \frac{1}{2}-\frac{1}{2{p}} \ = \ \frac{\theta}{q_{\tt d}} +\frac{1-\theta}{2}.
 \end{equation}
   Then, by~\eqref{eq:interpolation} and the standard inequality~$\|g\|_{L^{2{p}'}} \leq \|g\|_{L^{q_{\tt d}}}^{\theta}\|g\|_{L^2}^{1-\theta}$, we obtain
   \begin{equation}\label{sharp:tail-int}
   \left(\biint_{U_r}\nra{f_\ell}_{L^2(B_r)}^{2{p}'}\dx\dt\right)^\frac{1}{{p}'} \ \leq \ \left(\biint_{U_r}
   \nra{f_\ell}_{L^2(B_r)}^{q_{\tt d}} \dx\dt\right)^\frac{2\theta}{q_{\tt d}}\left(\biiint_{Q_r}f_\ell^2\dz\right)^{1-\theta}.
   \end{equation}

   Combining~\eqref{sharp:cacc:tail-holder} and~\eqref{sharp:tail-int}, upon applying Young's Inequality with exponent~$1/\theta$ and~$1/(1-\theta)$, for any~$\beta>0$, we obtain
   \begin{eqnarray}\label{eq:sharp-tail-final}
  && \frac{c}{(r-\rho)^{n+2s}(\kk-\ell)}\nra{\tail(f_\ell;B_r)}_{L^{{p}}(U_r)}\left(\biint_{U_r}\nra{f_\ell}_{L^2(B_r)}^{q_{\tt d}}\dx\dt\right)^\frac{2\theta}{q_{\tt d}}\left(\biiint_{Q_r}f_\ell^2\dz\right)^{1-\theta}\notag\\
  &&\qquad\qquad \qquad \leq \ \beta\left(\biint_{U_r}\nra{f_\ell}_{L^2(B_r)}^{q_{\tt d}}\dx\dt\right)^\frac{2}{q_{\tt d}}\ + c\left(\frac{\nra{\tail(f_\ell;B_r)}_{L^{{p}}(U_r)}}{(r-\rho)^{n+2s}(\kk-\ell)}\right)^\frac{1}{1-\theta}\biiint_{Q_r}f_\ell^2\dz\,,
  \end{eqnarray}
  where~$c \equiv c(n,s,\Lambda,\beta,p)>0$.
  
  Finally, noting that by~\eqref{def:theta-t} and~\eqref{def:drift-dim} we have
  \[ 
  \frac{1}{1-\theta} \ = \ \frac{{p}}{{p}-p_*}, 
    \]
  we can conclude the proof of~\eqref{caccioppoli-def} by combining~\eqref{eq:sharp-tail-final} into~\eqref{eq:cacc-final-1}. 
\end{proof}

\subsection{Hypoelliptic gain of integrability in~$(t,x)$}
This section is devoted to the proof of the local gain used in the De Giorgi iteration argument in Theorem~\ref{thm_bdd}. In order to prove forthcoming Theorem~\ref{thm:gain}, we localize the upper truncation~$(f-\kk)_+$ and apply a Kato inequality to the resulting equation. The local commutator is treated as an~$L^2_{t,x}H^{-s}_v$ source term, whereas the exterior contribution retains the form of the tail~$\tail((f-\ell)_+;B_r)(t,x)$ multiplied by the velocity superlevel set. We then apply the mixed-norm estimates for the  fractional Kolmogorov fundamental solution and estimate the contribution generated by~$h$ in~$L^\infty$ so that it is absorbed into the truncation level by a shift~$\lambda_{h,\kk}(r)$; combining all these estimates will finally yield the estimate below.
\begin{theorem}[Local gain in integrability]\label{thm:gain}
Let~$Q_1\subset\er^{1+2n}$, let~$s\in(0,1)$. Let~$p > p^*$ and~$q> N_s/(2s)$, with~$p_\ast$ given in~\eqref{eq:exponent-crt} and~$N_s$ given in~\eqref{def:homo-dim}.
For any~$r>0$ and any~$\kk \in \er$ define 
\begin{equation}\label{def:lamb-h-k}
\lambda_{h,\kk}(r) \ := \ 
 c_0\nra{h_+}_{L^{q}(Q_1)}
\left(\frac{|Q_r \cap \{f>\kk\}|}{|Q_r|}\right)^{\frac{1}{\theta}-\frac{1}{q}},
\end{equation}
for some~$c_0\equiv c_0(n,s,\Lambda,q)>0$ and~$\theta := \frac{1}{2}(q+N_s/(2s))$.

Assume that~$f \in \Wc(Q_1)$ is a weak  subsolution to~\eqref{problema} in ${Q}_1$. Then,
 for any~$ \frac{1}{2} \leq {\rho} < {r} \leq 1$, any~$\beta>0$ and any~$0 \leq \ell<\kk$ it holds
\begin{equation}\label{sharp:gain_scalar_source}
\begin{aligned}
&\left(\biint_{U_\rho}\nra{(f-\kk-\lambda_{h,\kk}(r))_+}_{L^2(B_\rho)}^{q_{\rm d}}\dx\dt\right)^\frac{2}{q_{\rm d}}\\
&\qquad\quad\leq \frac{c}{(r-\rho)^{2(1+2s)}}
\biiint_{Q_r}(f-\kk)_+^2\dz + \beta\left(\biint_{U_r}\nra{(f-\ell)_+}_{L^2(B_r)}^{q_{\rm d}}\dx\dt\right)^\frac{2}{q_{\rm d}}\\
&\qquad \qquad
 +\frac{c}{(r-\rr)^{2(1+2s)}}\biiint_{Q_r}\int_{B_r}
\frac{\left((f(t,x,v)-\kk)_+-(f(t,x,w)-\kk)_+\right)^2}{|v-w|^{n+2s}}\dw\dz\\
&\qquad\qquad +c\,\left(\frac{\nra{\tail(f_+;B_\frac{1}{2})}_{L^{{p}}(U_1)}}{(r-\rho)^{n+2s}(\kk-\ell)}\right)^\frac{\bar{p}}{\bar{p}-p_\ast}\biiint_{Q_r}(f-\ell)_+^2\dz\,,
\end{aligned}
\end{equation}
where~$q_{\tt d}$ is given in~\eqref{eq:exponent-crt}, $\overline p$ is defined in \eqref{eq:bar-p-T}, and~$c \equiv c(n,s, \Lambda, p, q,\beta)>0$.
\end{theorem}

\begin{proof} 
Fix~$ 1/2 \leq \rho < r \leq 1$ such that~$Q_{r} \subset \Om$, and, for~$\sigma:= (r+\rho)/2$,
define
two cut-off functions:~$\phi \in C^\infty_c(B_{\sigma^{1+2s}}\times B_{\sigma})$ such that~$0 \leq \phi \leq 1$,~$\phi \equiv 1$ on~$B_{\rho^{1+2s}}\times B_\rho$,
 \[
 \snr{\nabla_v \phi} \ \leq \ \frac{c}{r-\rho} \quad \text{and}  \quad \snr{v\cdot \nabla_x \phi} \ \leq  \ \frac{c}{(r-\rho)^{1+2s}}\,, 
 \]
 for some universal constant~$c>0$.

  For any~$\eps>0$, let~$\eta_\eps \in C^\infty_c ((-\sigma^{2s},0])$ be such that~$0 \leq \eta_\eps \leq 1$,~$\eta_\eps \equiv 1$ on~$[-\rr^{2s},-\eps]$. Moreover, once denoted with~$\mu_\eps(t):= \mathbbm{1}_{[-\eps,0]}(t) \eta_\eps'(t)$ we get that~$\mu_\eps\leq 0$ and
  \[
  \sup_{\eps>0} \| \eta_\eps' -\mu_\eps\|_{L^\infty}< \frac{c}{(r-\rho)^{1+2s}}< \infty.
  \]

  With this bit of notation we divide the proof into several steps.

  \medskip
  {\it Step 1 - The localized equation.} Consider the function~$g:=\eta_\eps\phi f_\kk$, where~$f_\kk :=(f-\kk)_+$. By weak chain-rule, in the sense of distribution, such an auxiliary function~$g$ satisfies
  \[
  \begin{aligned}
  (\partial_t+v\cdot\nabla_x)g & \ = \ f_\kk(\partial_t+v\cdot\nabla_x)(\eta_\eps\phi) +\eta_\eps\phi(\partial_t+v\cdot\nabla_x)f_\kk\\
  & \ = \ f_\kk(\partial_t+v\cdot\nabla_x)(\eta_\eps\phi)+\eta_\eps\phi\mathbbm{1}_{Q_{\sigma}\cap \{f>\kk\}}(\partial_t+v\cdot\nabla_x)(f-\kk)\\
 & \ \leq \ f_\kk \bigl((\eta_\eps'-\mu_\eps)\phi+\eta_\eps v\cdot\nabla_x\phi\bigr)
 +\eta_\eps\phi\mathbbm{1}_{Q_{\sigma}\cap\{f>\kk\}}\Lc(f-\kk) +\eta_\eps\phi h_+\mathbbm{1}_{Q_{\sigma}\cap \{f>\kk\}}\,,
\end{aligned}
\]
 where we have used~\eqref{eq:transport-truncation}, the sign of~$\mu_\eps$ and the fact that~$f$ is a weak subsolution of~\eqref{problema}.

 Now, we employ the integration by parts formula for the diffusion operator~$\Lc$:
 \begin{equation}\label{eq:product-rule}
 \Lc \big((f-\kk)\phi\big) \ = \ \phi\Lc(f-\kappa) + (f-\kk)\Lc\phi+\Ic(f-\kappa,\phi),
 \end{equation}
 where~$\Ic(\cdot,\cdot)$ is a reminder term given by
 \begin{equation}\label{eq:kato}
 \Ic(f-\kk,\phi)(v) \ := \ \int_{\ern}\bigl(f(v)-f(w)\bigr)\bigl(\phi(v)-\phi(w)\bigr)K(t,x,v,w)\dw\,,
 \end{equation}
  together with the Kato inequality (see, e.~\!g.,~\cite{CS18})
  \[
  \mathbbm{1}_{\{\Upphi>0\}}\Lc\Upphi \ \leq \ \Lc\Upphi_+.
  \]
  Then, since~$\eta_\eps$ depends only on time,~\eqref{eq:product-rule} and~\eqref{eq:kato} give
  \[
  \begin{aligned}
   \eta_\eps\phi\mathbbm{1}_{\{f>\kk\}}\Lc(f-\kk) & \ = \
  \mathbbm{1}_{\{f>\kk\}}\Lc[\eta_\eps\phi(f-\kk)]-\eta_\eps(f-\kk)\Lc\phi\mathbbm{1}_{\{f>\kk\}} -\eta_\eps\Ic(f-\kk,\phi)\mathbbm{1}_{\{f>\kk\}}\\
  & \ \leq \ \Lc g +\eta_\eps\mathbbm{1}_{\{f>\kk\}}(v)
  \int_{\ern}f_\kk(t,x,w)\bigl(\phi(x,v)-\phi(x,w)\bigr)K(t,x,v,w)\dw.
  \end{aligned}
  \]
  {
  Indeed, on~$\{f>\kk\}$, in view of ~\eqref{eq:product-rule} and~\eqref{eq:kato}, we have~$\Lc\Upphi_+-\Lc\Upphi=\eta_\eps\int_{\ern}\phi(w)(f(w)-\kk)_-K\dw$ with~$\Upphi:=\eta_\eps\phi(f-\kk)$, and it exactly compensates the replacement of~$f(w)-\kk$ by~$f_\kk(w)$ in the commutator, the difference between the two sides being~$\eta_\eps\phi(v)\int_{\ern}(f(w)-\kk)_-K\dw\geq0$, and outside~$\{f>\kk\}$, the left-hand side is zero, while the right-hand one is~$\Lc g(v) \geq 0$, being~$g(v)=0$ there. 
  }
  
  We now split the last integral into~$B_r$ and~$\er^n\setminus B_r$.  Adding~$(-\Delta_v)^sg$ to both sides, we obtain
  \begin{equation}\label{eq:gain:localized-eq}
  \begin{aligned}
   (\partial_t+v\cdot\nabla_x)g+(-\Delta_v)^sg
   &  \ \leq \ f_\kk\bigl((\eta_\eps'-\mu_\eps)\phi+\eta_\eps v\cdot\nabla_x\phi\bigr)+(\Lc+(-\Delta_v)^s)g\\
   & \quad + \eta_\eps\mathbbm{1}_{B_r}(v)\int_{B_r}f_\kk(t,x,w)(\phi(x,v)-\phi(x,w))K(t,x,v,w)\dw\\
   &\quad + F_T + F_h\,,
  \end{aligned}
  \end{equation}
  where in the last line we denote with
  \[
  \begin{aligned}
  F_T & \ := \ \eta_\eps\phi\mathbbm{1}_{\{f>\kk\}}(v)\int_{\er^n\setminus B_r}f_\kk(t,x,w)K(t,x,v,w)\dw\,,\\
  \text{and} \quad F_h & \ := \ \eta_\eps\phi h_+\mathbbm{1}_{\{f>\kk\}}.
  \end{aligned}
  \]
  Let us justify the presence of the indicator~$\mathbbm{1}_{B_r}(v)$ in the local term of~\eqref{eq:gain:localized-eq}.
  For~$v\notin B_r$ one has~$\phi(x,v)=0$, so the discarded local contribution equals~$-\eta_\eps\mathbbm{1}_{\{f>\kk\}}\int_{B_r}f_\kk(w)\phi(x,w)K\dw\leq0$; on~$\{f\leq\kk\}$, where the transport term of~$g$ vanishes by~\eqref{eq:transport-truncation} and~$g(v)=0$, inequality~\eqref{eq:gain:localized-eq} holds directly since 
  $$\Lc g(v)+\eta_\eps\mathbbm{1}_{B_r}(v)\int_{B_r}f_\kk(w)(\phi(v)-\phi(w))K\dw = \eta_\eps\phi(x,v)\int_{B_r}f_\kk(w)K\dw+\int_{\ern\setminus B_r}g(w)K\dw\geq0.$$
  
  \smallskip

  We remark that the local integral in~$B_r(v)$ appearing in~\eqref{eq:gain:localized-eq} is understood in the duality sense. Indeed, as in~\cite[Lemma~4.11]{IS20}, by pairing with~$\xi\in H^s(\ern)$, truncating the singularity and symmetrizing in~$(v,w)$, it is represented by
  \[
  \begin{aligned}
  \frac{\eta_\eps}{2}\iint_{B_r \times B_r} & \Big[(f_\kk(t,x,v)\bigl(\xi(v)-\xi(w)\bigr)+\xi(v)\bigl(f_\kk(t,x,w)-f_\kk(t,x,v)\bigr)\Big]
  \\
  &\hspace{2cm}\times\bigl(\phi(x,v)-\phi(x,w)\bigr)K(t,x,v,w)\dw\dv.
  \end{aligned}
  \]
  Hence, H\"older's Inequality and the Lipschitz bound on~$\phi$  give
  \begin{equation}\label{eq:singular-duality}
  \begin{aligned}
   & \left|\left\langle
   \eta_\eps\mathbbm{1}_{B_r}
    \int_{B_r}f_\kk(w)(\phi(\cdot)-\phi(w))K(\cdot,w)\dw,\xi\right\rangle_{H^{-s},H^s}\right|\\
   &\qquad\qquad \ \leq \ \frac{c}{r-\rho}\left(\|(f_\kk(t,x,\cdot)\|_{L^2(B_r)}
  + [(f_\kk(t,x,\cdot)]_{H^s(B_r)}\right)\|\xi\|_{H^s(\er^n)}.
   \end{aligned}
  \end{equation}
  By a similar duality argument the remaining nonlocal operators satisfy
  \begin{equation}\label{eq:perturbation-duality}
   \| (\Lc+(-\Delta_v)^s)g\|_{H^{-s}(\er^n)} \ \leq \ c\|g\|_{H^s(\er^n)}.
  \end{equation}

   All in all, we proved that the sum of the first three terms on the right-hand side of~\eqref{eq:gain:localized-eq} belongs to~$L^2(U_r;H^{-s}(\ern))$ and can be written as~$H_1+(-\Delta_v)^\frac{s}{2}H_0$, where~$H_0,H_1\in L^2(U_r\times\er^n)$ vanish outside~$U_r$ in the~$(t,x)$-variables. The preceding estimates~\eqref{eq:singular-duality} and~\eqref{eq:perturbation-duality} give
   \begin{equation}\label{sharp:gain:R-bound-1}
   \begin{aligned}
    & \frac{1}{|Q_r|}\left(\iiint_{U_r\times\ern}|H_0|^2\dz  +\iiint_{U_r\times\ern}|H_1|^2\dz\right)\\
    & \quad \leq \   \frac{c}{(r-\rr)^{2(1+2s)}}\biiint_{Q_r}f_\kk^2\d z 
     +\frac{c}{(r-\rr)^{2(1+2s)}}\biiint_{Q_r}\int_{B_r}
      \frac{\left(f_\kk(t,x,v)-f_\kk(t,x,w)\right)^2}{|v-w|^{n+2s}}\dw\dz\\
    & \qquad  +\frac{1}{|B_r|}\biint_{U_r} \iint_{\ern\times\ern}\frac{(g(t,x,v)-g(t,x,w))^2}{|v-w|^{n+2s}}\dw\dz.
    \end{aligned} 
    \end{equation}

    Now, we estimate the fractional norm of~$g$ in the right-hand side of~\eqref{sharp:gain:R-bound-1}. For this, we split 
    \[
    \begin{aligned}
            \iint_{\ern\times\ern}\frac{(g(t,x,v)-g(t,x,w))^2}{|v-w|^{n+2s}}\dw\dv & \ = \   \iint_{B_r \times B_r}\frac{(g(t,x,v)-g(t,x,w))^2}{|v-w|^{n+2s}}\dw\dv \\
            & \qquad + 2   \iint_{B_r \times (\ern\setminus B_r)}\frac{g^2(t,x,v)}{|v-w|^{n+2s}}\dw\dv .
    \end{aligned}
    \]
    Moreover since,~$g=\eta_\eps \phi f_\kk$, if~$v\in\operatorname{supp}(\phi)\subset B_{\sigma}$ and~$w \in \ern\setminus B_r$, then
    \begin{equation}\label{eq:kernel-est}
     \frac{|w|}{|v-w|} \ \leq \ 1 + \frac{|v|}{|w|-|v|} \ \leq \ 1+\frac{c r}{r-\rho}  \ \leq \ \frac{c}{r-\rho}\,,
    \end{equation}
   since~$1/2 \leq \rho < r \leq 1$.

 Thus, since~$g(t,x,\cdot)$ is supported in~$B_\sigma$ and~$|v-w| \geq r-\sigma = (r-\rho)/2$ for every~$v \in B_\sigma$ and~$w \in \ern\setminus B_r$, we get
   $$
    \iint_{B_r \times (\ern\setminus B_r)}\frac{g^2(t,x,v)}{|v-w|^{n+2s}}\dw\dv  \ \leq \ \frac{c}{(r-\rho)^{2s}}\int_{B_r}g^2\dv,
   $$
  
   On the other hand, by the identity
   \[
   f_\kk(v)\phi(v) - f_\kk(w)\phi(w) = f_\kk(v) \big(\phi(v)-\phi(w)\big) + \phi(w)(f_\kk(v)-f_\kk(w))
   \]
   we obtain
   \[
   \begin{aligned}
       \iint_{B_r \times B_r}\frac{(g(t,x,v)-g(t,x,w))^2}{|v-w|^{n+2s}}\dw\dv & \ \leq \ c\iint_{B_r \times B_r}\frac{f_\kk^2(t,x,v)(\phi(x,v)-\phi(x,w))^2}{|v-w|^{n+2s}}\dw\dv\\
       & \quad + c\iint_{B_r \times B_r}\frac{(f_\kk(t,x,v)-f_\kk(t,x,w))^2}{|v-w|^{n+2s}}\dw\dv\\
       & \ \leq \ \frac{c}{(r-\rho)^{2(1+2s)}}\int_{B_r}f_\kk^2(t,x,v)\dv  \\
       & \quad + c\iint_{B_r \times B_r}\frac{(f_\kk(t,x,v)-f_\kk(t,x,w))^2}{|v-w|^{n+2s}}\dw\dv\,,
   \end{aligned}
   \]
   upon proceeding as in~\eqref{eq:singular-int}. Therefore,
   \begin{equation}\label{eq:g-control}
   \begin{aligned}
       &\frac{1}{|B_r|}\biint_{U_r} \iint_{\ern\times\ern}\frac{(g(t,x,v)-g(t,x,w))^2}{|v-w|^{n+2s}}\dw\dz \\
       &\quad \leq \ \frac{c}{(r-\rho)^{2(1+2s)}}\biiint_{Q_r}f_\kk^2\dz  
        + c\biiint_{Q_r}\int_{B_r}\frac{(f_\kk(t,x,v)-f_\kk(t,x,w))^2}{|v-w|^{n+2s}}\dw\dz
   \end{aligned}
   \end{equation}
  All in all, by combining~\eqref{eq:g-control} with~\eqref{sharp:gain:R-bound-1} we finally arrive at
 \begin{equation}\label{sharp:gain:R-bound}
   \begin{aligned}
    & \frac{1}{|Q_r|}\left(\iiint_{U_r\times\ern}|H_0|^2\dz  +\iiint_{U_r\times\ern}|H_1|^2\dz\right)\\
    & \quad \leq \   \frac{c}{(r-\rr)^{2(1+2s)}}\biiint_{Q_r}f_\kk^2\d z 
     +\frac{c}{(r-\rr)^{2(1+2s)}}\biiint_{Q_r}\int_{B_r}
      \frac{\left(f_\kk(t,x,v)-f_\kk(t,x,w)\right)^2}{|v-w|^{n+2s}}\dw\dz.
    \end{aligned} 
    \end{equation}
    
  \medskip
  { \it Step 2 - Tail estimate.} Now, we establish an estimate on the tail term~$F_T$.
  For any~$\ell < \kk$, letting~$f_\ell := (f-\ell)_+$, by~\eqref{eq:kernel-est}, we can recenter the kernel into the origin so that, by ellipticity of~$K$, we get
   \begin{equation}\label{sharp:gain:FT-pointwise}
   0 \ \leq  \ F_T(t,x,v) \ \leq \ \frac{c}{(r-\rho)^{n+2s}}\tail(f_\ell;B_r)(t,x)\mathbbm{1}_{\{f(t,x,v)>\kk\}}.
   \end{equation}
  
    Since~$(f-\ell)_+ \geq \kk-\ell$ on~$\{f>\kk\}$, for a.~\!e.~$(t,x)\in U_r$, we obtain the following Chebychev's Inequality
    \begin{equation}\label{sharp:gain:cheb-fibre}
    \left(\mean{B_r}\mathbbm{1}_{\{f(t,x,v)>\kk\}}\dv\right)^\frac{1}{2} \ \leq \ \frac{\nra{f_\ell}_{L^2(B_r)}}{\kk-\ell}.
    \end{equation}

Define now
  \[
  \lambda \ := \ \frac{2p_*}{\bar{p}}-1\in[0,1).
  \]

   By~\eqref{eq:exponent-crt} and~\eqref{eq:exponent-identities}, 
   {we introduce the parameter $\bar m$, for which it holds}
   \[
   \frac{1}{\bar{m}} \ := \ \frac{1}{a_{\rm d}}-\frac{1}{\bar{p}} \ = \ \frac{\lambda}{q_{\rm d}}+\frac{1-\lambda}{2}
   \]

  Applying, in order, H\"older's Inequality with~$1/a_{\tt d} = 1 /\bar{m}+1/\bar{p}$  and the interpolation inequality~$\|\psi\|_{L^{\bar{m}}} \leq \|\psi\|_{L^{q_{\tt d}}}^{\lambda}\|\psi\|_{L^2}^{1-\lambda}$, we have
 \[
\begin{aligned}
  & r^{2s}\left(\biint_{U_r}\nra{F_T}_{L^2(B_r)}^{a_{\rm d}}\dx\dt\right)^\frac{1}{a_{\rm d}}\\
  & \quad \  \leq \ \frac{c}{(r-\rho)^{n+2s}(\kk-\ell)}\left(\biint_{U_r}\tail(f_\ell,B_r)^{a_{\rm d}}\nra{f_\ell}_{L^2(B_r)}^{a_{\rm d}}\dx\dt\right)^\frac{1}{a_{\rm d} }\\
   & \quad \  \leq \ \frac{c\,\nra{\tail(f_\ell;B_r)}_{L^{\bar{p}}(U_r)}}{(r-\rho)^{n+2s}(\kk-\ell)}\left(\biint_{U_r}\nra{f_\ell}_{L^2(B_r)}^{\bar{m}}\dx\dt\right)^\frac{1}{\bar{m}}\\
   & \quad \leq \ \frac{c\,\nra{\tail(f_\ell;B_r)}_{L^{\bar{p}}(U_r)}}{(r-\rho)^{n+2s}(\kk-\ell)}\left(\biint_{U_r}\nra{f_\ell}_{L^2(B_r)}^{q_{\tt d}}\dx\dt\right)^\frac{\lambda}{q_{\tt d}}
    \left(\biiint_{Q_r}f_\ell^2\dz\right)^\frac{1-\lambda}{2}.
   \end{aligned}
   \]
 
 Lastly, we apply Young's Inequality with~$1/\lambda$ and~$1/(1-\lambda)$, for every~$\beta>0$, in the previous inequality, to get
  \begin{equation}\label{sharp:gain:FT-mixed-input}
  \begin{aligned}
   & r^{4s}\left(\biint_{U_r}\nra{F_T}_{L^2(B_r)}^{a_{\rm d}}\dx\dt\right)^\frac{2}{a_{\rm d}}\\
   & \quad \leq \
   \beta\left(\biint_{U_r}\nra{f_\ell}_{L^2(B_r)}^{q_{\tt d}}\dx\dt\right)^\frac{2}{q_{\tt d}} + \left(\frac{c\,\nra{\tail(f_\ell;B_r)}_{L^{\bar{p}}(U_r)}}{(r-\rho)^{n+2s}(\kk-\ell)}\right)^\frac{\bar{p}}{{\bar{p}-p_\ast}}\biiint_{Q_r}f_\ell^2\dz\,,
  \end{aligned}
  \end{equation}
   where~$c \equiv c(n,s,\Lambda,\beta, p)>0$ and
    we have also used that $\frac{2}{1-\lambda} \ = \ \frac{\bar{p}}{\bar{p}-p_*}$.

    \medskip
   {\it Step 3 - Comparison and potential estimates.} We now apply the higher integrability properties given by the kinetic semigroup~\eqref{sharp:m1:eq_potential}. Indeed, upon extending the right-hand side to zero outside the support of the cut-off  in~\eqref{eq:gain:localized-eq}, we consider the solution~$G$ to the following problem
   \begin{equation}\label{eq:cauchy-kin}
	\begin{cases}
		(\partial_t+v\cdot\nabla_x)G + (-\Delta_v)^sG =  H_1 + (-\Delta_v)^\frac{s}{2}H_{0} + F_T + F_h\quad &\text{in}~(-\sigma^{2s},\infty)\times \er^{2n},\\*[0.5ex]
		G(-\sigma^{2s},x,v) =0 \quad & \text{on}~\{-\sigma^{2s}\}\times\er^{2n},
	\end{cases}
	\end{equation}
    where~$H_1$,~$H_0$ are the representation via duality introduced in~\eqref{eq:singular-duality} and~\eqref{eq:perturbation-duality}. 

Note that such a solution~$G$ does exist via convolution with the fundamental solution in~\eqref{eq:gain:localized-eq}. Moreover, once we define
    \begin{equation}\label{def:G}
    G_{\rm loc} \ := \ \mathscr P\bigl(H_1+(-\Delta_v)^\frac{s}{2}H_0\bigr)\,, \quad G_T \ := \ \mathscr P F_T \quad \text{and} \quad G_h \ := \ \mathscr P F_h\,, 
    \end{equation}
   then the maximum principle for the Cauchy problem~\eqref{eq:cauchy-kin} -- see, e.~\!g.,~\cite[Lemma~A.12]{IS20} -- applied to the function~${g}-G$ gives
    \begin{equation}\label{sharp:gain:comparison}
     0 \ \leq \ g \ \leq \ G =  G_{\rm loc} + G_T + G_h.
    \end{equation}
Hence, applying Corollary~\ref{sharp:m1:cor_rough_local}, 
     together with~\eqref{sharp:gain:R-bound}, implies
     \begin{equation}\label{sharp:gain:Gloc}
     \begin{aligned}
     \left(\biint_{U_r}\nra{G_{\rm loc}}_{L^2(B_r)}^{q_{\rm d}}\dx\dt\right)^\frac{2}{q_{\rm d}}
      &  \ \leq \ \frac{c}{(r-\rho)^{2(1+2s)}}\biiint_{Q_r}f_\kk^2\dz\\
      & \quad + \frac{c}{(r-\rho)^{2(1+2s)}}\biiint_{Q_r}\int_{B_r}\frac{\left(f_\kk(t,x,v)-f_\kk(t,x,w)\right)^2}{|v-w|^{n+2s}}\dw\dz.
      \end{aligned}
      \end{equation}
      
     As for the term involving the tail, using first Corollary~\ref{sharp:m1:cor_local}-{(i)}  combined with~\eqref{sharp:gain:FT-mixed-input}, we get
     \begin{equation}
     \begin{aligned}\label{sharp:gain:GT-mixed}
     &\left(\biint_{U_r}\nra{G_T}_{L^2(B_r)}^{q_{\rm d}}\dx\dt\right)^\frac{2}{q_{\rm d}}\\
     &\quad \leq \ \beta\left(\biint_{U_r}\nra{f_\ell}_{L^2(B_r)}^{q_{\rm d}}\dx\dt\right)^\frac{2}{q_{\rm d}}+\left(\frac{c\,\nra{\tail(f_\ell;B_r)}_{L^{\bar{p}}(U_r)}}{(r-\rho)^{n+2s}(\kk-\ell)}\right)^\frac{\bar{p}}{{\bar{p}-p_\ast}}\biiint_{Q_r}f_\ell^2\dz.
     \end{aligned}
     \end{equation}
Now, we are left to estimate the source term in~$G_h$. For this, we set
\[
\theta \ := \ \frac{1}{2}\left(q+\frac{N_s}{2s}\right),
\qquad
\frac{N_s}{2s}< \theta < q,
\qquad
\frac{1}{\theta}-\frac1{q}>0.
\]
H\"older's Inequality in the definition of~$\mathscr{P}F_h$ gives
\[
|\mathscr PF_h(z)|
\ \leq \
\left(
\int_0^{r^{2s}}
\|P_t\|_{L^{\theta'}_{x,v}(\er^{2n})}^{\theta'}\,\d t
\right)^\frac{1}{\theta'}
\|F_h\|_{L^\theta(Q_r)}.
\]
{Since, by the scaling relation~\eqref{eq:fundamental} and Proposition~\ref{prop:fund}~{(i)},} we have that
\[
\|P_t\|_{L^{\theta'}(\er^{2n})}
\ \leq \ c \,t^{-n(1+\frac{1}{s})(1-\frac{1}{\theta'})},
\]
we obtain
\[
\int_0^{r^{2s}}
\|P_t\|_{L^{\theta'}_{x,v}(\er^{2n})}^{\theta'}\,\d t \ \leq \ c \int_0^{r^{2s}}t^{-\frac{n}{\theta-1}(1+\frac{1}{s})}\dt  < +\infty,
\]
{where we also used the fact that~$\theta>N_s/(2s)$.}

Therefore, since~$F_h =\eta_\eps\phi h_+\mathbbm{1}_{\{f>k\}}$, applying H\"older's Inequality once more gives
\begin{equation}\label{sharp:gain:Gh-infinity}
\|G_h\|_{L^\infty}
\ \leq  \ c\nra{h_+}_{L^{q}(Q_r)}
\left(\frac{|Q_r \cap \{f>k\}|}{|Q_r|}\right)^{\frac{1}{\theta}-\frac1{q}}.
\end{equation}

Now, we choose~$c_0$ in the definition of~$\lambda_{h,\kk}(r)$ in~\eqref{def:lamb-h-k} so that the right-hand side of~\eqref{sharp:gain:Gh-infinity} is at most~$\lambda_{h,\kk}(r)$. Then, since~$g \equiv \eta_\eps(f-k)_+$ on~$Q_\rho$, the comparison in~\eqref{sharp:gain:comparison} gives
\begin{equation}\label{eq:shifted-comp}
\bigl(\eta_\eps(f-\kk)_+-\lambda_{h,\kk}(r)\bigr)_+
\leq |G_{\rm loc}|+G_T
\quad\hbox{on}~Q_\rho\,,
\end{equation}
{and we note that~$\bigl(\eta_\eps(f-\kk)_+-\lambda_{h,\kk}(r)\bigr)_+ \to (f-\kk-\lambda_{h,\kk}(r))_+$ a.~\!e. as~$\eps$ goes to~$0$.}

We first combine~\eqref{sharp:gain:Gloc} with~\eqref{sharp:gain:GT-mixed}, and then we let~$\eps\downarrow0$, via Fatou's Lemma, since the constants above are independent of~$\eps$. This yields the desired~\eqref{sharp:gain_scalar_source}.
\end{proof}

\subsection{Proof of the~$L^2$-$L^\infty$ estimate}
This section is devoted to the proof of Theorem~\ref{thm_bdd}. We begin establishing two auxiliary results that will be its building blocks. In order to make the forthcoming computations less heavy, we introduce some further notation that will be also adopted throughout the proof of Theorem~\ref{thm_bdd}.

First of all, for any~$r>0$ and any~$\kk \in \er$ we define
\begin{equation}\label{def:energy-j}
\Ec_{f,\kk}(r) \ := \ \biiint_{Q_r}(f-\kk)_+^2\dz
+\biiint_{Q_r}\int_{B_r}\frac{\left((f(t,x,v)-\kk)_+-(f(t,x,w)-\kk)_+\right)^2}{|v-w|^{n+2s}}\dw\dz\,,
\end{equation}
and
\begin{equation}\label{def:level-k}
\Ac^+(\kk,z_0,r) \ := \ \frac{|Q_r(z_0) \cap \{ f > \kk\}|}{|Q_r|} \quad \text{with} \quad \Ac^+(\kk,0,r) = \Ac^+(\kk,r).
\end{equation}

The first auxiliary result is a nonlinear recurrence estimate that will be coupled with~\eqref{sharp:gain_scalar_source} in order to build the proper iteration scheme that will lead to~\eqref{kinetic_special}.
\begin{lemma}\label{lemma:recurrence}
Let~$Q_1\subset\er^{1+2n}$ and~$p^*$ be defined in~\eqref{eq:exponent-crt},~$p > p_\ast$,~$\bar p$ defined in~\eqref{eq:bar-p-T} and~$q_{\bar{p}}$ being the exponent associated to $\bar p$ in~\eqref{sharp:m1:eq_qp_formula}.
For any~$r>0$ and~$\kk \in \er$ let~$\lambda_{h,\kk}(r)$ be defined in~\eqref{def:lamb-h-k}. 
Assume that~$f \in \Wc(Q_1)$ is a weak  subsolution to~\eqref{problema} in ${Q}_1$. Then, for any~$ \frac{1}{2} \leq {\rho} < {r} \leq 1$ and any~$0 \leq \ell<\kk$ it holds
   \begin{equation}\label{eq:old-recurrence}
   \begin{aligned}
   & \biiint_{Q_\rr}(f-\kk-\lambda_{h,\kk}(r))_+^{2}\dz\\
   & \qquad\quad  \ \leq \ \frac{c\,\Ec_{f,\kk}(r)}{(r-\rho)^{2(n+2s)}}\left(\Ac^+(\kk,r)\right)^\frac{2s}{N_s}\\
   & \qquad\qquad +\ c\,\left(\frac{\nra{\tail(f_+;B_\frac{1}{2})}_{L^{p}(U_1)}}{(\kk-\ell)(r-\rho)^{n+2s}}\right)^2\left(\Ac^+(\kk,r)\right)^{1-\frac{2}{q_{\bar{p}}}}\biiint_{Q_r}(f-\ell)_+^2 \dz\,,
   \end{aligned}
   \end{equation}
  where~$N_s$ is given in~\eqref{def:homo-dim}
  and~$c \equiv c (n,s,\Lambda,p)>0$.
  \end{lemma}
  
  \begin{proof}
  The initial part of the proof follows the localization and comparison properties established in Theorem~\ref{thm:gain}; see respectively Step~1 and~Step~3 there.

For any~$r>0$ and~$\kk\in \er$, let~$\lambda_{h,\kk}(r)$ be given in~\eqref{def:lamb-h-k}. We recall
that, by~\eqref{eq:shifted-comp}, we have that 
\begin{equation}\label{sharp:bdd:comparison-recurrence}
(\eta_\eps f_{\kk}-\lambda_{h,\kk}(r))_+ \ \leq \ |G_{\rm loc}|+G_T \quad \text{on}~Q_\rho\,,
\end{equation}
where~$G_{\rm loc}$ and~$G_T$ are defined in~\eqref{def:G}, and~$\eta_\eps$ is the time cut-off introduced in the proof of Theorem~\ref{thm:gain}.

Let~$\bar{p}$ be given in~\eqref{eq:bar-p-T},~$q_{\bar{p}}$ being the one in in~\eqref{sharp:m1:eq_qp_formula} associated with~$\bar{p}$ and~$q_0 = 2N_s/(N_s-2s)$. By H\"older's Inequality, the estimate in~\eqref{sharp:bdd:comparison-recurrence} yields
\begin{equation}\label{shapr:bdd-recurrence-y}
\begin{aligned}
\biiint_{Q_\rho} (\eta_\eps f_{\kk}-\lambda_{h,\kk}(r))_+ ^2\dz & \ \leq \ c\,\left(\biiint_{Q_r}|G_{\rm loc}|^{q_0}\dz\right)^\frac{2}{q_0}\left(\Ac^+(\kk,r)\right)^\frac{2s}{N_s}\\
& \quad + c\,\left(\biiint_{Q_r}|G_{T}|^{q_{\bar{p}}}\dz\right)^\frac{2}{q_{\bar{p}}}\left(\Ac^+(\kk,r)\right)^{1-\frac{2}{q_{\bar{p}}}}\,,
\end{aligned}
\end{equation}
since, by~\eqref{sharp:m1:eq_ap_qp},~$q_{\bar{p}} >2$ being~$\bar{p}$ in~\eqref{eq:bar-p-T} greater than~$p_\ast$, and where~$\Ac^+(\kk,r)$ is given in~\eqref{def:level-k}.

We now estimate separately all the different contributions in the right-hand side of~\eqref{shapr:bdd-recurrence-y}.

\medskip
\noindent
We begin estimating the quantity involving~$G_{\rm loc}$. For this, we apply Lemma~\ref{lemma:kolmogorov-young} and~\eqref{sharp:gain:R-bound}, so that, recalling the definition of~$\mathcal{E}_{f,\kk}(r)$ in~\eqref{def:energy-j}, we get
\begin{equation}\label{eq:source-Gloc}
 \left(\biiint_{Q_r}|G_{\rm loc}|^{q_0}\dz\right)^\frac{2}{q_0}  \ \leq \ \frac{c\,\mathcal{E}_{f,\kk}(r)} {(r-\rho)^{2(1+2s)}}
\end{equation}

 Now, we give an estimate for the tail of~$f_+$ that will control the term~$G_T$ in~\eqref{shapr:bdd-recurrence-y}. Let~$\bar{p}$ be given in~\eqref{eq:bar-p-T} and let~$q_{\bar{p}}$ and~$a_{\bar{p}}$ be the associated exponents in~\eqref{sharp:m1:eq_ap_qp}. Since~$a_{\bar{p}} = 2\bar{p}/(\bar{p}+2)$ we have that
  \[
  \frac{1}{a_{\bar{p}}}   
  \ = \ \frac{1}{2} +\frac{1}{\bar{p}}.
  \]
  Hence, by H\"older's Inequality,~\eqref{sharp:gain:FT-pointwise} and~\eqref{sharp:gain:cheb-fibre} we get
  \begin{equation}\label{sharp:gain-tail-old}
   \begin{aligned} 
   & r^{2s}\left(\biint_{U_r}\nra{F_T}^{a_{\bar{p}}}_{L^2(B_r)}\dx\dt\right)^\frac{1}{a_{\bar{p}}} \\
   &\qquad \quad  \leq \ \frac{c}{(r-\rho)^{n+2s}(\kk-\ell)}\left(\biint_{U_r}\tail(f_\ell;B_r)^{a_{\bar{p}}} \nra{f_\ell}_{L^2(B_r)}^{a_{\bar{p}}}\dx\dt\right)^\frac{1}{a_{\bar{p}}}\\*[0.8ex]
   &\qquad \quad \leq \ \frac{c\nra{\tail(f_\ell;B_r)}_{L^{\bar{p}}(U_r)}}{(r-\rho)^{n+2s}(\kk-\ell)}\left(\biiint_{Q_r}f_\ell^2 \dz\right)^\frac{1}{2}.
  \end{aligned}
  \end{equation}
  
  Thus, Corollary~\ref{sharp:m1:cor_local}~{(\it iii)} applied with~$\bar{p}$ -- which by~\eqref{eq:bar-p-T} is admissible being greater than~$p_\ast$ in~\eqref{eq:exponent-crt} --  combined with~\eqref{sharp:gain-tail-old}  gives
     \begin{equation}\label{sharp:gain:GT-scalar}
      \left(\biiint_{Q_r}|G_T|^{q_{\bar{p}}}\dz\right)^\frac{2}{q_{\bar{p}}}
      \ \leq \ c\,\left(\frac{\nra{\tail(f_+;B_\frac{1}{2})}_{L^{{p}}(U_1)}}{(r-\rho)^{n+2s}(\kk-\ell)}\right)^2\biiint_{Q_r}f_\ell^2\dz.
      \end{equation}

All in all, combining~\eqref{eq:source-Gloc} and~\eqref{sharp:gain:GT-scalar} into~\eqref{shapr:bdd-recurrence-y}, and letting~$\eps \to 0^+$, Fatou's Lemma yields the desired estimate in~\eqref{eq:old-recurrence}.
 \end{proof}

Now, we prove a precise estimate on the source contribution appearing in the energy estimate~\eqref{caccioppoli-def}. This will be a fundamental part when building the iteration procedure needed for the forthcoming proof of Theorem~\ref{thm_bdd}.

\begin{prop}\label{prop:source-contribution}
Let~$Q_1\subset\er^{1+2n}$. Let~$p > p_\ast$, with~$p_\ast$ given in~\eqref{eq:exponent-crt},~$q_{{p}}$ being the parameter in~\eqref{sharp:m1:eq_qp_formula} associated with~${p}$ and~$q > N_s/(2s)$, with~$N_s$ given in~\eqref{def:homo-dim}. Assume that~$f \in \Wc(Q_1)$ is a weak  subsolution to~\eqref{problema} in ${Q}_1$. 
    Then, for any~$ \frac{1}{2} \leq {\rho} < {r} \leq 1$, any~$\eps \in(0,1]$ and any~$0 \leq \ell<\kk$, it holds
    \begin{equation}\label{eq:source-contribution}
  \begin{aligned}
   & \biiint_{Q_\rho}h_+(f-\kk)_+ \dz \\
   & \quad  \ \leq \  \eps\Ec_{f,\kk}(r) + \frac{c\nra{h_+}_{L^{q}(Q_1)}^2}{\eps(r- \rho)^{2(n+2s)}}\left(\biiint_{Q_r}\frac{(f-\ell)_+^2}{(\kk-\ell)^2}\dz\right)^{1+\frac{2s}{N_s}-\frac{2}{q}}\\
   & \qquad + \frac{c\,\nra{\tail(f_+;B_\frac{1}{2})}_{L^{{p}}(U_1)}\nra{h_+}_{L^{q}(Q_1)}}{(r-\rho)^{2(n+2s)}}\left(\biiint_{Q_r}\frac{(f-\ell)_+^2}{(\kk-\ell)^2}\dz\right)^{\frac{3}{2}-\frac{1}{q}-\frac{1}{q_{p}}}\,,
\end{aligned}
\end{equation}
where~$c\equiv c(n,s,\Lambda,p) >0$.
\end{prop}

\begin{proof}
In view of~\eqref{sharp:gain:comparison}, we have that 
\begin{equation}\label{sharp:bdd:unshifted-comparison}
\eta_\eps f_{\kk} \ \leq \ |G_{\rm loc}|+G_T+G_h \quad \text{on}~Q_\rho\,,
\end{equation}
where~$G_{\rm loc}$,~$G_T$ and~$G_h$ are defined in~\eqref{def:G} and~$\eta_\eps$ is the time cut-off defined in the proof of Theorem~\ref{thm:gain}.
Below, we estimate separately all the different term in the right hand-side of \eqref{sharp:bdd:unshifted-comparison}.

\vspace{1mm}
\noindent
By applying H\"older's Inequality three times, firstly with $q$ and $q/(q-1)$, then with~$2$ and~$2$ and lastly with $N_s(q-1)/(q(N_s-2s))$ and $N_s(q-1)/(2sq-N_s)$, since~$2q/(q-1) < q_0:= 2N_s/(N_s-2s)$ being~$q > N_s/(2s)$, we obtain
\begin{equation}\label{eq:source-loc1}
\begin{aligned}
\biiint_{Q_{\rho}}h_+\mathbbm{1}_{\{f>\kk\}}|G_{\rm loc}|\dz & \  \leq \ c\nra{h_+}_{L^{q}(Q_1)}\left(\Ac^+(\kk,r)\right)^{\frac{1}{2}-\frac{1}{2q}}\left(\biiint_{Q_r}\mathbbm{1}_{\{f>\kk\}}|G_{\rm loc}|^\frac{2q}{q-1}\dz\right)^\frac{q-1}{2q}\\
& \leq \ c\nra{h_+}_{L^{q}(Q_1)}\left(\Ac^+(\kk,r)\right)^{\frac{1}{2}+\frac{s}{N_s}-\frac{1}{q}} \left(\biiint_{Q_r}|G_{\rm loc}|^{q_0}\dz\right)^\frac{1}{q_0}\,,
\end{aligned}
\end{equation}
where~$\Ac(\kk,r)$ is given in~\eqref{def:level-k}.

Thus, by combining~\eqref{eq:source-Gloc} and~\eqref{eq:source-loc1} with Young's Inequality, and recalling the definition of~$\Ec_{f,\kk}(\cdot)$ in~\eqref{def:energy-j}, we obtain that, for any~$\eps>0$, it holds
\begin{equation}\label{eq:source-loc}
\begin{aligned}
& \biiint_{Q_{\rho}}h_+\mathbbm{1}_{\{f>\kk\}}|G_{\rm loc}|\dz \\
& \quad  \ \leq \ \eps\Ec_{f,\kk}(r) + \frac{c\nra{h_+}^2_{L^{q}(Q_1)}}{\eps(r-\rho)^{2(1+2s)}}\left(\biiint_{Q_r}\frac{f_\ell^2}{(\kk-\ell)^2}\dz\right)^{1+\frac{2s}{N_s}-\frac{2}{q}}\,,
\end{aligned}
\end{equation}
for some~$c \equiv c(n,s,\Lambda)>0$ and where we have used that
 \begin{equation}\label{eq:source-level}
 1< \frac{f_\ell^2}{(\kk-\ell)^2} \quad \text{on}~\{f >\kk\},\,\,\text{for}~\ell < \kk.
 \end{equation}

Similarly, by Lemma~\ref{lemma:fundamental} with~$\displaystyle a_0  :=  \frac{2N_s}{N_s+2s}$, 
so that~$\displaystyle \frac{1}{a_0} -\frac{2s}{N_s}  =  \frac{1}{q_0}$,
recalling that~$F_h = \eta_\eps\phi\, h_+ \mathbbm{1}_{\{f>\kk\}}$, we obtain that
\begin{equation}\label{eq:source-Gh}
\left(\biiint_{Q_r}|G_h|^{q_0}\dz\right)^\frac{1}{q_0}  \ \leq \ c \left(\biiint_{Q_r}|F_h|^{a_0}\right)^\frac{1}{a_0} \ \leq \ c\left(\biiint_{Q_1}h_+^{q}\dz\right)^\frac{1}{q}\left(\Ac^+(\kk,r)\right)^{\frac{1}{2}+\frac{s}{N_s}-\frac{1}{q}}\,,
\end{equation}
where we have applied H\"older's Inequality with exponents~$q/a_0$ and~$q/(q-a_0)$, which are admissible since~$q > N_s/(2s)$; in the display above~$c\equiv c(n,s)>0$.

Thus, by applying again H\"older's Inequality, first with~$a_0$ and~$q_0$  and then with~$q/a_0$ and~$q/(q-a_0)$, by~\eqref{eq:source-Gh} and~\eqref{eq:source-level} we get
\begin{equation}\label{eq:source-h}
\begin{aligned}
\biiint_{Q_{\rho}}h_+\mathbbm{1}_{\{f>\kk\}}|G_h|\dz & \  \leq \ c\left(\biiint_{Q_r}h_+^{a_0}\mathbbm{1}_{\{f>\kk\}}\dz\right)^\frac{1}{a_0}\left(\biiint_{Q_r}|G_h|^{q_0}\dz\right)^\frac{1}{q_0}\\
& \ \leq \ c\left(\biiint_{Q_1}h_+^{q}\dz\right)^\frac{2}{q}\left(\biiint_{Q_r}\frac{f_\ell^2}{(\kk-\ell)^2}\dz\right)^{1+\frac{2s}{N_s}-\frac{2}{q}}\,,
\end{aligned}
\end{equation}
for some~$c \equiv c(n,s)>0$.

\smallskip

We are left with estimating the remaining tail term. For this, since~${p} > p_\ast$, we apply~Corollary~\ref{sharp:m1:cor_local}~{\it (iii)} obtaining
\[
\left(\biiint_{Q_r}|G_T|^{q_{{p}}}\dz\right)^\frac{1}{q_{p}}  \ \leq \ c\left(\biint_{U_r}\nra{F_T}_{L^2(B_r)}^{a_{{p}}}\dx\dt\right)^\frac{1}{a_{{p}}}
\]
where, recalling~\eqref{sharp:m1:eq_ap_qp}, we have~$a_{{p}} = 2{p}/({p}+2)$, so that~$\displaystyle \frac1{a_{{p}}} \ = \ \frac{1}{{p}}+\frac{1}{2}$.
  Thus,~\eqref{sharp:gain:FT-pointwise}, H\"older's Inequality in~$(t,x)$ and~\eqref{sharp:gain:cheb-fibre} yields
  \[
  \begin{aligned}
   & \left(\biint_{U_r}\nra{F_T}_{L^2(B_r)}^{a_{{p}}}\dx\dt\right)^\frac{1}{a_{{p}}}\\
   &\qquad  \quad \ \leq \ \frac{c}{(r-\rho)^{n+2s}}\left(\biint_{U_r}\tail(f_\ell,B_r)^{a_{{p}}}\left(\mean{B_r}\mathbbm{1}_{\{f(t,x,v)>\kk\}}\dv\right)^\frac{a_{{p}}}{2}\dx\dt\right)^\frac{1}{a_{{p}}}\\
   & \qquad \quad \ \leq \ \frac{c}{(r-\rho)^{n+2s}(\kk-\ell)}\left(\biint_{U_r}\tail(f_\ell,B_r)^{a_{{p}}}\nra{f_\ell}_{L^2(B_r)}^{a_{{p}}}\dx\dt\right)^\frac{1}{a_{{p}}}\\*[0.8ex]
  &\qquad \quad \leq \ \frac{c\,\nra{\tail(f_\ell;B_r)}_{L^{{p}}(U_r)}}{(r-\rho)^{n+2s}}\left(\biiint_{Q_r}\frac{f_\ell^2}{(\kk-\ell)^2}\dz\right)^\frac{1}{2}\,,
  \end{aligned}
  \]
 where~$c \equiv c (n,s)>0$.
  
  All in all, we have obtained that
  \begin{equation}\label{eq:source-GT}  
      \left(\biiint_{Q_r}|G_T|^{q_{{p}}}\dz\right)^\frac{1}{q_{p}}  \ \leq \ \frac{c\,\nra{\tail(f_+;B_\frac{1}{2})}_{L^{{p}}(U_1)}}{(r-\rho)^{n+2s}}\left(\biiint_{Q_r}\frac{f_\ell^2}{(\kk-\ell)^2}\dz\right)^\frac{1}{2}.
  \end{equation}

  Now, note that, by~\eqref{sharp:m1:eq_ap_qp} and the definition of~$q_0$, we have
  \[
  \frac{1}{q_{{p}}}-\frac{1}{2} 
  \leq 0\,,
  \]
  since~${p} > p_\ast = N_{\tt d}/(2s) $.
   
   Therefore, since~$q > N_s/(2s)>2$, its H\"older conjugate~$q/(q-1)<2 $, by applying H\"older's Inequality twice, firstly with~$q$ and~$q/(q-1)$ and then with~$q_{p}(q-1)/q$ and~$q_{p}(q-1)/(q_p(q-1)-q)$, we obtain
   \begin{eqnarray}\label{eq:source-tail1}
             \biiint_{Q_r}h_+\mathbbm{1}_{\{f>\kk\}}|G_T| \dz && \ \leq \  \left(\biiint_{Q_r}h_+^{q}\dz\right)^\frac{1}{q}\left(\biiint_{Q_r}\mathbbm{1}_{\{f>\kk\}}|G_T|^\frac{q}{q-1}\dz\right)^\frac{q-1}{q}\notag\\
             && \ \leq \ \left(\biiint_{Q_r}h_+^{q}\dz\right)^\frac{1}{q}\left(\frac{|Q_r\cap \{f>\kk\}|}{|Q_r|}\right)^{1-\frac{1}{q}-\frac{1}{q_{p}}}
             \left(\biiint_{Q_r}|G_T|^{q_{p}}\dz\right)^\frac{1}{q_{p}}.
   \end{eqnarray}

   Finally, by~\eqref{eq:source-tail1},~\eqref{eq:source-GT} and~\eqref{eq:source-level} we get
    \begin{equation}\label{eq:source-tail}
       \begin{aligned}
&             \biiint_{Q_r}h_+\mathbbm{1}_{\{f>\kk\}}|G_T| \dz 
             \\
             & \qquad \quad \leq \ \frac{c\,\nra{\tail(f_+;B_\frac{1}{2})}_{L^{{p}}(U_1)}\nra{h_+}_{L^{q}(Q_1)}}{(r-\rho)^{n+2s}}\left(\biiint_{Q_r}\frac{f_\ell^2}{(\kk-\ell)^2}\dz\right)^{\frac{3}{2}-\frac{1}{q}-\frac{1}{q_{p}}}.
       \end{aligned}
   \end{equation}

  The desired estimate~\eqref{eq:source-contribution} will then follow by combining~\eqref{eq:source-loc},~\eqref{eq:source-h} and~\eqref{eq:source-tail}  after letting~$\eps\downarrow0$, using Fatou's Lemma, since all the constants in the above estimates are independent of~$\eps$.
\end{proof} 
\medskip
We are now in the position to prove of Theorem~\ref{thm_bdd}.

\begin{proof}[Proof of Theorem~\ref{thm_bdd}]
For every~$\kk\in\er$, let~$f_\kk:=(f-\kk)_+$, and define for any~$j \in \en$ the following sequences of radii and truncation levels
\[
r_j \ := \ \frac{1}{2}(1+2^{-j}), \quad \kk_j \ := \ \kk(1-2^{-j}) \quad \widehat{\kk}_j:=\frac{\kk_j+\kk_{j+1}}{2} \quad \text{and} \quad \widehat{r}_j:=\frac{r_j+r_{j+1}}{2}.
\]
In this way we have that~$r_j\downarrow 1/2$,~$\kk_j\uparrow \kk$, and
\begin{equation}\label{sharp:bdd:level-gaps}
\widehat{\kk}_j-\kk_j\ = \ \kk_{j+1}-\widehat{\kk}_j \ = \ 2^{-j-2}\kk \quad \text{and} \quad 
r_j-r_{j+1}\ = \ 2^{-j-2}.
\end{equation}

Set now the following quantities,
\begin{equation}\label{sharp:bdd:xy}
\mathcal{X}_j \ := \ \left(\biint_{U_{r_j}}\frac{\nra{f_{\kk_j}}_{L^2(B_{r_j})}^{q_{\rm d}}}{\kk^{q_{\rm d}}}\dx\dt\right)^\frac{2}{q_{\rm d}}
\quad \text{and} \quad
\mathcal{Y}_j \ := \ \biiint_{Q_{r_j}}\frac{f_{\kk_j}^2}{\kk^2}\dz.
\end{equation}
We now divide the proof into several steps.

\medskip
{\it Step 1 - Energy  estimate at the intermediate level.}  Let us assume that for a fixed~$\delta\in(0,1]$, $\kk$  satisfies the following lower bound
\begin{equation}\label{sharp:bdd:K-tail-source}
 \kk \ \geq \ \delta\nra{\tail(f_+;B_\frac{1}{2})}_{L^{p}(U_1)} 
+ c_0\nra{h_+}_{L^{q}(Q_1)},
\end{equation}
where~$c_0>1$ is a fixed constant that will be chosen later on.

Now, recalling the definition of~$\Ec_{f,\hat{\kk}_j}(\cdot)$ in~\eqref{def:energy-j}. We claim that, for every~$\beta\in(0,1)$, 
there are~$c_\beta\geq1$ and~$b>1$,with~$b$ independent of~$\beta$,
 such that
\begin{equation}\label{sharp:bdd:nested-energy}
\frac{\Ec_{f,\hat{\kk}_j}(\widehat{r}_j)}
{(\widehat{r}_j-r_{j+1})^{2(1+2s)}}
\ \leq \ \beta\kk^2\mathcal{X}_j+\frac{c_{\beta}\,b^j\kk^2}{\delta^\frac{\bar{p}}{\bar{p}-p_\ast}}\left(\mathcal{Y}_j+\mathcal{Y}_j^{1+\frac{2s}{N_s}-\frac{2}{q}} +\mathcal{Y}_j^{\frac{3}{2}-\frac{1}{q}-\frac{1}{q_{p}}}\right)\,,
\end{equation}
where~$q_{{p}}$ is the one in in~\eqref{sharp:m1:eq_qp_formula} associated with~${p}$, and~$\bar{p}$ is given by~\eqref{eq:exponent-crt}.

Indeed, fix
\[
\widehat{r}_j \ \leq \ \rho \ < \ R \ \leq  \ r_j \quad \text{and} \quad \rho_1 \ := \ \frac{\rho+R}{2}\,,
\]
and apply Lemma~\ref{lemma:energy-est} with~$\widehat{\kk}_j>\kk_j$ and~$\rho < \rho_1$. 

Now, since~$\bar{p} \leq p$ by~\eqref{eq:bar-p-T}, applying H\"older's Inequality, together with~\eqref{sharp:bdd:level-gaps} and~\eqref{sharp:bdd:K-tail-source}, yields
\begin{equation}\label{eq:tail-trivial}
\frac{\nra{\tail(f_{\kk_j};B_{\rho_1})}_{L^{\bar{p}}(U_{\rho_1})}}
{\widehat{\kk}_j-\kk_j} \ \leq  \ \frac{c\, 2^j}{\delta}.
\end{equation}

Choosing the parameter~$\beta :=  c_1\tilde{\beta}(R-\rho)^{2(1+2s)}$, where~$\tilde{\beta}>0$ and~$c_1>0$ is  sufficiently small, the first term in~\eqref{caccioppoli-def} is then bounded by~$\tilde{\beta}\kk^2\mathcal{X}_j$. Moreover, by~\eqref{eq:tail-trivial} we obtain that there exists~$m = m(n,s,\bar{p},p_\ast)>1$ and~$b>1$ such that
\begin{equation}\label{sharp:bdd:energy-source-pre}
\Ec_{f,\hat{\kk}_j}(\rho) \ \leq \ \tilde{\beta}\kk^2 \mathcal{X}_j
+\frac{c_{\tilde{\beta}}\,b^j\kk^2\mathcal{Y}_j}{(R-\rho)^m \delta^\frac{\bar{p}}{\bar{p}-p_\ast}} + c\biiint_{Q_{\rho_1}}h_+f_{\hat{k}_j}\dz.
\end{equation}

We now estimate the source term on the right-hand side in~\eqref{sharp:bdd:energy-source-pre} by relying on Proposition~\ref{prop:source-contribution}. Indeed, we apply~\eqref{eq:source-contribution} with~$1/2 \leq \rho_1 < R \leq 1$,~$\kk_j < \widehat{\kk}_j$ and~$\eps>0$. This combined with~\eqref{sharp:bdd:xy},~\eqref{sharp:bdd:K-tail-source} and~\eqref{def:energy-j} yields
\begin{equation}\label{sharp:bdd:source-final}
\biiint_{Q_{\rho_1}}h_+f_{\hat{\kk}_j}\dz  \ \leq \ \varepsilon \Ec_{f,\hat{\kk}_j}(R)
+\frac{c\,\kk^2 \tilde{b}^j}{\eps\delta(R-\rho)^{2(n+2s)}}\left(\mathcal{Y}_j^{1+\frac{2s}{N_s}-\frac{2}{q}} +\mathcal{Y}_j^{\frac{3}{2}-\frac{1}{q}-\frac{1}{q_{p}}}\right)\,,
\end{equation}
where~$c \equiv c(n,s,\Lambda)>0$,~$\tilde{b}>1$ and~$q_{p}$ is the one given in~\eqref{sharp:m1:eq_qp_formula}.

Combining the inequality in the display above with~\eqref{sharp:bdd:energy-source-pre}, upon enlarging the constants~$c>0$,~$m>1$ and~$b>1$ if necessary, yields that
\begin{equation}\label{sharp:bdd:hole-pre}
\Ec_{f,\hat{\kk}_j}(\rho) \ \leq \ \varepsilon \Ec_{f,\hat{\kk}_j}(R)+\tilde{\beta}\kk^2 \mathcal{X}_j
+\frac{c_{\tilde{\beta}}\,b^j\kk^2}{(R-\rho)^m \delta^\frac{\bar{p}}{\bar{p}-p_\ast}}\left(\mathcal{Y}_j+\mathcal{Y}_j^{1+\frac{2s}{N_s}-\frac{2}{q}} +\mathcal{Y}_j^{\frac{3}{2}-\frac{1}{q}-\frac{1}{q_{p}}}\right).
\end{equation}

Applying now Lemma~\ref{lemma_giusti} to~\eqref{sharp:bdd:hole-pre} with~$\Psi(\cdot):= \Ec_{f,\hat{\kk}_j}(\cdot)$,~$\widehat{r}_j$,~$r_{j}$,~$\beta_1:=m$,
\[
A_1 \ := \ \frac{c_{\tilde{\beta}}\,b^j\kk^2}{\delta^\frac{\bar{p}}{\bar{p}-p_\ast}}\left(\mathcal{Y}_j+\mathcal{Y}_j^{1+\frac{2s}{N_s}-\frac{2}{q}} +\mathcal{Y}_j^{\frac{3}{2}-\frac{1}{q}-\frac{1}{q_{p}}}\right), \quad A_2 \ := \ 0 \quad \text{and} \quad A_3 \ := \ \tilde{\beta}\kk^2\mathcal{X}_j,
\]
gives
\[
\Ec_{f,\hat{\kk}_j}(\widehat{r}_j) \ \leq \ c\,\tilde{\beta}\kk^2 \mathcal{X}_j
+\frac{c_{\tilde{\beta}}\,b^j\kk^2}{\delta^\frac{\bar{p}}{\bar{p}-p_\ast}}\left(\mathcal{Y}_j+\mathcal{Y}_j^{1+\frac{2s}{N_s}-\frac{2}{q}} +\mathcal{Y}_j^{\frac{3}{2}-\frac{1}{q}-\frac{1}{q_{p}}}\right)\,,
\]
for larger~$b>1$, since~$(r_j-\widehat{r}_j)^{-m} = 2^{m(j+3)}$, for some universal~$c>0$.

Choosing now~$\tilde{\beta} := \beta(\widehat{r}_j-r_{j+1})^{2(1+2s)}/c$ for any~$\beta>0$ yields  the desired estimate in~\eqref{sharp:bdd:nested-energy}, since the inverse powers of~$\tilde{\beta}$ produced by Young's Inequality are fixed powers of~$2^j$ and are therefore absorbed into~$c_{\beta}b^j$.

\medskip
{\it Step 2 - The coupled nonlinear recurrences.} We now build the double recurrence. First of all, let us begin noticing that, recalling the definition of~$\Ac^+(\hat{\kk}_j,\hat{r}_j)$ in~\eqref{def:level-k}, by~\eqref{sharp:bdd:level-gaps} and~\eqref{sharp:bdd:xy}, we have
\begin{equation}\label{sharp:bdd:cheb}
\Ac^+(\hat{\kk}_j,\hat{r}_j)  \ \leq \  \biiint_{Q_{\hat{r}_j}}\frac{f_{{\kk}_j}^2}{(\widehat{\kk}_j-\kk_j)^2}\dz
 \ \leq \ \frac{2^{2(j+2)}}{|Q_{\hat{r}_j}|} \int_{Q_{r_j}}\frac{f_{\kk_j}^2}{\kk^2}\dz  \ \leq \ c2^{2(j+2)}\mathcal{Y}_j\,,
\end{equation}
where~$c \equiv c(n,s)>0$ since~$|Q_{r_j}||Q_{\hat{r}_j}|^{-1} \leq c(n,s)$, for any~$j \in \en$.

Therefore, by~\eqref{def:lamb-h-k}, we obtain that the condition
\begin{equation}\label{sharp:bdd:source-gap}
  \frac{c_0\,\nra{h_+}_{L^{q}(Q_1)}b^j \mathcal{Y}_j^{\frac{1}{\theta}-\frac{1}{q}} }{\kk}  \ \leq \ 1 \quad \forall j \in \en\,,
\end{equation}
where~$N_s/(2s) < \theta< q$ and~$c_0(n,s,\Lambda,q) > 1$ is sufficiently large, and~\eqref{sharp:bdd:cheb} implies
\begin{equation}\label{sharp:bdd:lambda-cond}
    \lambda_{h,\hat{\kk}_j} (\widehat{r}_j) \ \leq \ \kk_{j+1} -\widehat{\kk}_j \quad \forall j \in \en\,,
\end{equation}
upon enlarging~$c_0$, if necessary.

Assume now that~\eqref{sharp:bdd:source-gap} and
\begin{equation}\label{eq:con-Yj}
\mathcal{Y}_j \ \leq \ 1 \quad j \in \en\,,
\end{equation}
do hold. Then, by applying Theorem~\ref{thm:gain} with~$r_{j+1} < \widehat{r}_j$,~$\kk_j<\widehat{\kk}_j$ we get that~\eqref{sharp:bdd:lambda-cond} and~\eqref{sharp:bdd:nested-energy} yield
\begin{equation}\label{sharp:bdd:x-recurrence}
\mathcal{X}_{j+1} \ \leq \ \beta\mathcal{X}_j+\frac{c_{\beta}\,b^j}{\delta^\frac{\bar{p}}{\bar{p}-p_\ast}}\mathcal{Y}_j^\alpha,
\end{equation}
with
\[
\alpha \ := \ \min\left\{1,1+\frac{2s}{N_s}-\frac{2}{q},\frac{3}{2}-\frac{1}{q}-\frac{1}{q_{p}}\right\} \in (0,1].
\]

On the other hand, by applying  Lemma~\ref{lemma:recurrence} with~$r_{j+2} < \widehat{r}_{j+1}$,~$\kk_{j+1}<\widehat{\kk}_{j+1}$ we get that a combination of H\"older's Inequality, a similar estimate as in~\eqref{sharp:bdd:cheb},~\eqref{sharp:bdd:lambda-cond} and~\eqref{sharp:bdd:nested-energy}, gives
\begin{equation}\label{sharp:bdd:y-recurrence}
	\mathcal{Y}_{j+2}  \ \leq \  \frac{c\,b^j}{\delta^\frac{\bar{p}}{\bar{p}-p_\ast}}
	\left[\bigl(\mathcal{X}_j+\mathcal{Y}_j^\alpha\bigr)\mathcal{Y}_{j+1}^{\frac{2s}{N_s}}\
	+\mathcal{Y}_j\mathcal{Y}_{j+1}^{1-\frac{2}{q_{\bar{p}}}}\right]\,,
\end{equation}
upon sufficiently enlarging~$c>0$ and~$b>1$, where~$q_{\bar{p}}$ is the one in in~\eqref{sharp:m1:eq_qp_formula} associated with~$\bar{p}$ in~\eqref{eq:bar-p-T}.

\medskip
{\it Step 3 Iteration:} Let us begin noticing that in~\eqref{def:lamb-h-k}, since~$\theta \in (N_s/(2s), q)$,~$1/\theta-1/q>0$. Also, since~$q > N_s/(2s)$, by~\eqref{HP:source}, we have that~$1+4s/N_s-2/q>1$. Lastly, by~\eqref{sharp:m1:eq_ap_qp} and  
by~\eqref{HP:tail}, we have that
\[
\frac{1}{2}-\frac{1}{q_{p}} \ =  \ \frac{2s}{N_s} - \frac{N_{\tt d}}{N_s p} 
\ > \ 0\,,
\]
so that, by~\eqref{HP:source}, we have
\[
\frac{3}{2}-\frac{1}{q}-\frac{1}{q_{p}} +\frac{2s}{N_s} 
\ > \ 1.
\]

All in all, we obtain that
\[
\alpha + \frac{2s}{N_s} \ > \ 1.
\]

Now, letting~$b>1$ be the one appearing in~\eqref{sharp:bdd:x-recurrence} and~\eqref{sharp:bdd:y-recurrence}, we choose~$\omega\in(0,1)$ sufficiently small so that
\begin{equation}\label{sharp:bdd:omega-choice}
{b\left(\omega^{\alpha}+ \omega^{1-\frac{2}{q_{\bar{p}}}} + \omega^{\frac{1}{\theta}-\frac{1}{q}}\right) + b^2 \omega^{\alpha+\frac{2s}{N_s}-1} \ \leq \ \frac{1}{4},}
\end{equation}
and set
\begin{equation}\label{eq:useful-def}
\sigma_0 \ := \ b \omega^\alpha \in (0,1)\,, \quad \beta \ := \ \frac{\sigma_0}{2} \quad \text{and} \quad M \ := \ \max \left\{1,\frac{2\max\{c,c_\beta\}}{\sigma_0\delta^\frac{\bar{p}}{\bar{p}-p_\ast}}\right\}\,,
\end{equation}
where~$c_\beta>0$ is the constant appearing in~\eqref{sharp:bdd:x-recurrence} (and it is fixed by the previous choice of~$\beta$), while~$c>0$ is the one appearing in~\eqref{sharp:bdd:y-recurrence}.

Finally choose~$\eps\in(0,1)$ sufficiently small such that
\begin{equation}\label{sharp:bdd:epsilon-choice}
\frac{2\max\{c,c_\beta\}}{\delta^\frac{\bar{p}}{\bar{p}-p_\ast}}(M+1)
\eps^{\alpha+\frac{2s}{N_s}-1}\omega^{\frac{2s}{N_s}} \ \leq \ \frac{\omega^2}{2} \quad \text{and} \quad 
{
\frac{2\max\{c,c_\beta\}}{\delta^\frac{\bar{p}}{\bar{p}-p_\ast}}\eps^{1-\frac{2}{q_{\bar{p}}}}\omega^{1-\frac{2}{q_{\bar{p}}}}
\ \leq  \ \frac{\omega^2}{2}.
}
\end{equation}

Assuming that the initial conditions below are satisfied:
\begin{equation}\label{sharp:bdd:initial-decay}
\begin{cases}
\mathcal{Y}_1 \ \leq \ \eps\omega,\\
\mathcal{Y}_2 \ \leq \ \eps\omega^2,\\
\mathcal{X}_1 \ \leq  \ M\eps^{\alpha}\sigma_0\,.
\end{cases}
\end{equation}
Then, for every~$j \in \en$, it holds
\begin{equation}\label{sharp:bdd:decay-induction}
\mathcal{Y}_j  \ \leq \ \eps\omega^j \quad \text{and} \quad 
\mathcal{X}_j \ \leq \  M\eps^{\alpha}\sigma_0^j.
\end{equation}

Indeed, let us assume that for~$j \in \en$, it holds
\[
\begin{cases}
\mathcal{Y}_j \ \leq \ \eps\omega^j,\\
\mathcal{Y}_{j+1} \ \leq \ \eps\omega^{j+1},\\
\mathcal{X}_j \ \leq  \ M\eps^{\alpha}\sigma_0^j.
\end{cases}
\]

Then,~\eqref{sharp:bdd:source-gap} (assuming~\eqref{sharp:bdd:K-tail-source} holds with~$c_0>1$ as in~\eqref{sharp:bdd:source-gap}) and~\eqref{eq:con-Yj} are verified. Moreover, letting~$c_0$ being as in~\eqref{sharp:bdd:source-gap}, if~$c_0\nra{h_+}_{L^{q}(Q_1)} \leq \kk$, then, by~\eqref{sharp:bdd:omega-choice} since~$\eps \in (0,1)$, we get
\[
\frac{c_0\nra{h_+}_{L^{q}(Q_1)}\,b^j\mathcal{Y}_j^{\frac{1}{\theta}-\frac{1}{q}}}{\kk} \ \leq \ \eps^{\frac{1}{\theta}-\frac{1}{q}}(b\omega^{\frac{1}{\theta}-\frac{1}{q}})^j \ \leq \ 1,
\]
so~\eqref{sharp:bdd:source-gap} and both recurrences~\eqref{sharp:bdd:x-recurrence} and~\eqref{sharp:bdd:y-recurrence} hold true. Moreover, since~$(b\omega^{\alpha})^j=\sigma_0^j$,~\eqref{sharp:bdd:x-recurrence} yields
\[
\mathcal{X}_{j+1} \ \leq  \ \beta M\eps^{\alpha}\sigma_0^j
+\frac{M\eps^\alpha\sigma_0^{j+1}}{2}
\ \leq \ M\varepsilon^{\alpha}\sigma_0^{j+1}\,,
\]
since~$\beta = \sigma_0/2$.

For~$\mathcal{Y}_{j+2}$ we obtain via~\eqref{sharp:bdd:y-recurrence},~\eqref{sharp:bdd:omega-choice} and~\eqref{sharp:bdd:epsilon-choice} that
\[
\begin{aligned}
\mathcal{Y}_{j+2} & \ \leq  \ \frac{\max\{c,c_\beta\}}{\delta^\frac{\bar{p}}{\bar{p}-p_\ast}}\left\{(M+1)\eps^{\alpha+\frac{2s}{N_s}}
\bigl(b^2\omega^{\alpha+\frac{2s}{N_s}}\bigr)^j\omega^{\frac{2s}{N_s}} 
{+\eps^{2-\frac{2}{q_{\bar{p}}}}
\bigl(b\omega^{2-\frac{2}{q_{\bar{p}}}}\bigr)^j\omega^{1-\frac{2}{q_{\bar{p}}}}
}\right\}\\
& \ \leq \ \frac{\eps\omega^{j+2}}{2}\bigl(b^2\omega^{\alpha+\frac{2s}{N_s}-1}\bigr)^j + \frac{\eps\omega^{j+2}}{2}\bigl(b\omega^{1-\frac{2}{q_{\bar{p}}}}\bigr)^j
\\ 
& \ \leq \ \eps \omega^{j+2}.
\end{aligned}
\]

This proves~\eqref{sharp:bdd:decay-induction}, and in particular
\begin{equation}\label{sharp:bdd:xy-decay}
\mathcal{X}_j+\mathcal{Y}_j \to 0 \quad \text{as}~j \to +\infty.
\end{equation}

We are thus left to show that initial conditions in~\eqref{sharp:bdd:initial-decay} are satisfied in order to start the iteration.  For this,  recall that
\[
\kk_0 \ = \ 0, \quad \widehat{\kk}_0 \ = \ \frac{\kk}{4} \quad \text{and} \quad \kk_1 \ = \  \frac{\kk}{2}.
\]
Repeating now the proof of Lemma~\ref{lemma:energy-est} with radii~$\widehat{r}_0$ and~$1$ and level~$0<\widehat{\kk}_0$, and estimating the exterior interaction before using the two-level inequality, gives
\begin{equation}\label{sharp:bdd:initial-energy}
\mathcal{\Ec}_{f,\hat{\kk_0}}(\widehat{r}_0)
\ \leq \ c\,\nra{f_+}_{L^2(Q_1)}^2
+c\,\left(\nra{\tail(f_+;B_\frac{1}{2})}_{L^{p}(U_1)} + \nra{h_+}_{L^{q}(Q_1)}\right)\nra{f_+}_{L^2(Q_1)}\,,
\end{equation}
by the inclusions~$L^2\hookrightarrow L^{p'}$ and~$L^2\hookrightarrow L^{q'}$, since~$p>p_*>2$ and~$q>N_s/(2s)>2$.

Take now the source term~$G_{h,0}$ in~\eqref{sharp:gain:comparison} at level~$\widehat{\kk}_0$. By~\eqref{sharp:gain:Gh-infinity} and the Chebychev inequality
\begin{equation}\label{sharp:bdd:initial-alpha}
\Ac^+(\widehat{\kk}_0,\widehat{r}_0) \ \leq \ \frac{c\,\nra{f_+}_{L^2(Q_1)}^2}{\kk^2}\,,
\end{equation}
we obtain that
\begin{equation}\label{sharp:bdd:initial-source-gap}
\frac{\|G_{h,0}\|_{L^\infty}}{\kk}
\ \leq \ \frac{c_0\,\nra{h_+}_{L^{q}(Q_1)}}{\kk}
\left(\frac{\nra{f_+}_{L^2(Q_1)}}{\kk}\right)^{\frac{2}{\theta}-\frac{2}{q}} \ \leq \ \frac{1}{4}\,,
\end{equation}
upon choosing~$\kk \geq c_0\nra{h_+}_{L^{q}(Q_1)} + c\nra{f_+}_{L^2(Q_1)}$, i.~\!e.~$\|G_{h,0}\|_\infty\leq \kk/4=\kk_1-\widehat{\kk}_0$.

We estimate now~$\Xb_1$. For this, let~$F_T$ be the tail contribution at level~$\widehat{\kk}_0$ and let
\[
\beta_T \ := \ 2\left(\frac1{a_{\rm d}}-\frac1{\bar p}\right)
=1+\frac1{p_*}-\frac2{\bar p}\in(0,1].
\]

By H\"older's Inequality with~$1/a_{\tt  d}=1/\bar p+\beta_T/2$,~\eqref{sharp:gain:FT-pointwise} and~\eqref{sharp:bdd:initial-alpha}, since $a^\frac{1}{\beta_T} \leq a$, for~$a \in [0,1]$, we obtain
\begin{equation}\label{sharp:bdd:tail-input}
\begin{aligned}
& \left(\biint_{U_{\hat{r}_0}}
\nra{F_T}_{L^2(B_{\hat{r}_0})}^{a_{\tt d}}\dx\dt\right)^\frac{1}{a_{\tt d}}\\
&\qquad\qquad \quad \leq \ c\,
\nra{\tail(f_+;B_\frac{1}{2})}_{L^{p}(U_1)}\left({ \biint_{U_{\hat{r}_0}}}\left(\mean{B_{\hat{r}_0}}\mathbbm{1}_{\{f > \hat{\kk}_0\}}\dv\right)^{1/\beta_T}\dx\dt\right)^\frac{\beta_T}{2}\\ 
&\qquad\qquad \quad \leq \  c\,\nra{\tail(f_+;B_\frac{1}{2})}_{L^{p}(U_1)}
\left(\frac{\nra{f_+}_{L^2(Q_1)}}{\kk}\right)^{\beta_T}.
\end{aligned}
\end{equation}

All in all, we can follow the proof Theorem~\ref{thm:gain} with~$\hat{r}_0< 1$,~$\hat{\kk}_0$, with~\eqref{sharp:bdd:tail-input} in place of~\eqref{sharp:gain:FT-mixed-input}. Then, since by~\eqref{sharp:bdd:initial-source-gap} we obtain that~$\lambda_{h,\hat{\kk}_0}(\hat{r}_0) \leq \kk_1 -\hat{\kk}_0$, the energy estimate~\eqref{sharp:bdd:initial-energy} together with~\eqref{sharp:gain_scalar_source}  and~\eqref{sharp:bdd:K-tail-source} yield
\begin{equation}\label{sharp:bdd:initial-x}
\mathcal{X}_1 \ \leq \  c_\delta\left[
\frac{\nra{f_+}_{L^2(Q_1)}^2}{\kk^2}
+\frac{\nra{f_+}_{L^2(Q_1)}}{\kk}
+\left(\frac{\nra{f_+}_{L^2(Q_1)}}{\kk}\right)^{2\beta_T}
\right].
\end{equation}

Also, by H\"older's Inequality we get
\begin{equation}\label{sharp:bdd:initial-y}
\mathcal{Y}_1+\mathcal{Y}_2 \ \leq \ c\,\left(\frac{\nra{f_+}_{L^2(Q_1)}}{\kk}\right)^2.
\end{equation}

Therefore, by choosing~$\kk$ such that
\[
\kk \ :=  \ c_1 \nra{f_+}_{L^2(Q_1)}+ c_0\,\nra{h_+}_{L^{q}(Q_1)}+ \delta\,\nra{\tail(f_+;B_\frac{1}{2})}_{L^{p}(U_1)},
\]
with~$c_0 \equiv c_0(n,s,\Lambda,q)>1$ and~$c_1 \equiv c_1(\eps,\omega,M,\delta)>1$ sufficiently large so that~\eqref{eq:con-Yj},~\eqref{sharp:bdd:initial-decay} and~\eqref{sharp:bdd:initial-source-gap} hold. Let us remark that, via~\eqref{sharp:bdd:omega-choice},~\eqref{eq:useful-def} and~\eqref{sharp:bdd:epsilon-choice},~$c_1$ actually depends on~$n$,~$s$,~$\Lambda$,~$p$,~$q$ and~$\delta$.

Finally, by~\eqref{sharp:bdd:xy-decay},  since~$r_j\downarrow1/2$ and~$\kk_j\uparrow \kk$, we obtain that~$f \leq \kk$ almost everywhere on~$Q_\frac{1}{2}$. Then, the estimate in~\eqref{kinetic_special} for general~$z_0$ and~$r >0$ will follow by the change of variables of Remark~\ref{rem:blow-up}, which gives  the~$|v_0|$-dependence of the constant and  the factor~$r^{2s}$ in front of the source term.
\end{proof}

\section{Nonlocal Strong Harnack inequality}
\label{sec_harnack}

This section is devoted to the completion of the proof of the strong Harnack inequality in Theorem~\ref{thm_strong}.

 \subsection{Proof of Theorem~\ref{thm_strong}}\label{sec_proof-strong} 
  Let~$r_{0}$ and~$\zeta$ be given by the weak Harnack inequality in Theorem~\ref{thm:weak}. 
 We assume that~$\zeta \in (0,2)$. This is admissible, up to lowering~$\zeta$ via H\"older's Inequality, being~$|Q^-_{r_{0}}|$ finite.
  
  Set ~$1/2 \leq  \sigma' < \sigma \leq 1$ and~$z_{0}:=(-1+r_{0}^{2s},0,0)$ and for any~$z_1 \in {Q}_{\sigma'r_{0}}(z_{0})$ it holds that
  \[
  {Q}_{(c_*(\sigma-\sigma')r_{0})^\beta}(z_1) \subset {Q}_{\sigma r_{0}}(z_{0}).
  \]
   Indeed, by~\eqref{eq:inclusion} in Lemma~\ref{lemma:cov} applied with radii~$\sigma' < \sigma$ and center~$\delta_{1/r_{0}}(z_{0})$, upon applying the dilation~$\delta_{r_{0}}$ -- for which~$\delta_{r_{0}}(z\circ z')=\delta_{r_{0}}(z)\circ\delta_{r_{0}}(z')$ and~$\delta_{r_{0}}({Q}_{\tau}(z))={Q}_{\tau r_{0}}(\delta_{r_{0}}(z))$ -- we get~${Q}_{r_{0}(c_*(\sigma-\sigma'))^\beta}(z_1) \subset {Q}_{\sigma r_{0}}(z_{0})$, and it suffices to note that~$(c_*(\sigma-\sigma')r_{0})^\beta \leq r_{0}(c_*(\sigma-\sigma'))^\beta$, being~$\beta\geq1$ and~$r_{0}<1$.
  
   Hence, applying the boundedness estimate in~\eqref{kinetic_special} in Theorem~\ref{thm_bdd} 
   with radius $(c_*(\sigma-\sigma')r_{0})^\beta$, and observing that 
   ${Q}_{(c_*(\sigma-\sigma')r_{0})^\beta}(z_1) \subseteq {Q}_{\sigma r_{0}}(z_{0})$, we get
  \begin{equation}\label{eq:bdd-ref-1} 
  \begin{aligned}
  	f(z_1)   & \ \leq \  \frac{c(\delta)}{[(\sigma-\sigma')r_{0}]^{\beta_{0}}}\left(\iiint_{{Q}_{\sigma r_{0}}(z_{0})}f^2 \d z\right)^\frac{1}{2}\\
	& \quad +
  	\delta\left(\biint_{{U}_{(c_*(\sigma-\sigma')r_{0})^\beta}(z_1)}\tail(f;B_{\frac{1}{2}(c_*(\sigma-\sigma')r_{0})^\beta}(v_1))^{p} \dx\dt\right)^\frac{1}{p}\,,
    \end{aligned}
  \end{equation}
  for some~$\beta_{0} \equiv \beta_{0}(n,s)>0$and $z_1 \in {Q}_{\sigma'r_{0}}(z_{0})$. 
  
Now, we estimate the nonlocal term.
  For a.~\!e.~$(t,x) \in U_{(c_*(\sigma-\sigma')r_{0})^\beta}(z_1) \subset U_{\sigma r_{0}}(-1+r_{0}^{2s},0,0)$, the latter inclusion being the projection onto the drift variables of~${Q}_{(c_*(\sigma-\sigma')r_{0})^\beta}(z_1) \subset {Q}_{\sigma r_{0}}(z_{0})$,  we have that
  $$
  \begin{aligned}
  &\Big[\displaystyle\tfrac{1}{2}(c_*(\sigma-\sigma')r_{0})^\beta\Big]^{2s} \int_{\ern \setminus B_{\frac{1}{2}(c_*(\sigma-\sigma')r_{0})^\beta}(v_1)}\frac{f(t,x,w)}{|w-v_1|^{n+2s}}\dw \\
  &\quad  =  \Big[\tfrac{1}{2}(c_*(\sigma-\sigma')r_{0})^\beta\Big]^{2s} \int_{B_{\sigma r_{0}} \setminus B_{\frac{1}{2}(c_*(\sigma-\sigma')r_{0})^\beta}(v_1)}\frac{f(t,x,w)}{|w-v_1|^{n+2s}}\dw \notag\\
  & \qquad   +  \Big[\tfrac{1}{2}(c_*(\sigma-\sigma')r_{0})^\beta\Big]^{2s} \int_{\ern \setminus B_{\sigma r_{0}}}\frac{f(t,x,w)}{|w-v_1|^{n+2s}}\dw\\
  & \quad\ \leq c\,  \,\operatorname*{ess\,sup}_{{Q}_{\sigma r_{0}}(z_{0})}f + \frac{c\,r_{0}^{2s}}{(\sigma-\sigma')^{n+2s}}\int_{\ern \setminus B_\frac{r_{0}}{2}}\frac{f(t,x,w)}{|w|^{n+2s}}\dw,
  \end{aligned}
  $$
  where we have used the fact that {$z_0=(-1+r_0^{2s}, 0,0)$,}~${Q}_{(c_*(\sigma-\sigma')r_{0})^\beta}(z_1) \subset {Q}_{\sigma r_{0}}(z_{0})$ for any~$z_1 \in {Q}_{\sigma'r_{0}}(z_{0})$,~$\sigma > \frac{1}{2}$ and that
  \[
  \frac{|w|}{|w-v_1|} \ \leq \ 1 +\frac{|v_1|}{\left||w|-|v_1|\right|} \leq 1 + \frac{\sigma'}{\sigma-\sigma'} \  \leq \  \frac{1}{\sigma-\sigma'}\,,
  \]
for any~$w \in \ern \setminus B_{\sigma r_{0}}$, and that~$\beta \geq 1$ and~$r_{0}<1$.

  Thus, we can estimate the $p$-contribution of the tail in velocity as follows,
  \begin{equation}\label{eq:tail-est-ref-bdd}
  \begin{aligned}
  	&\left(\biint_{U_{(c_*(\sigma-\sigma')r_{0})^\beta}(z_1)} \tail(f; B_{\frac{1}{2}(c_*(\sigma-\sigma')r_{0})^\beta}(v_1))^{p}\dx\dt\right)^\frac{1}{p}\\
  	&\qquad\qquad  \quad  \leq \   c\, \operatorname*{ess\,sup}_{{Q}_{\sigma r_{0}}(z_{0})}f  + \frac{c}{(\sigma-\sigma')^{n+2s+\frac{\beta N_{\rm d}}{p}}  { r_0^{\frac{(\beta - 1)N_d}{p} }}}    \left(\biint_{U_{r_{0}}({z_0})}\tail(f; B_\frac{r_{0}}{2})^{p}\dx\dt\right)^\frac{1}{p}\,,
  \end{aligned}
  \end{equation}
  where we have used that~${Q}_{(c_*(\sigma-\sigma')r_{0})^\beta}(z_1) \subset {Q}_{\sigma r_{0}}(z_{0})\subset {Q}_{r_{0}}(z_{0})$, whence~$U_{(c_*(\sigma-\sigma')r_{0})^\beta}(z_1)\subset U_{r_{0}}(-1+r_{0}^{2s},0,0)$, together with
 $$
  \left(\frac{|U_{r_{0}}|}{|U_{(c_*(\sigma-\sigma')r_{0})^\beta}|}\right)^\frac{1}{p} \ = \ \left(\frac{r_{0}}{(c_*(\sigma-\sigma')r_{0})^\beta}\right)^\frac{N_{\rm d}}{p} \ \leq \ \frac{c}{(\sigma-\sigma')^{\frac{\beta N_{\rm d}}{p}}  { r_0^{\frac{(\beta - 1)N_d}{p} }}} \,,
 $$
  for some~$c\equiv c(n,s,p)>0$.

  Then, combining~\eqref{eq:tail-est-ref-bdd} with~\eqref{eq:bdd-ref-1},
  applying the interpolation inequality with parameters $\zeta$, $\infty$ and $\theta=\zeta/2$,
  we obtain 	
  \begin{equation}\label{eq:harnack-3}
        \begin{aligned}
    		\operatorname*{ess\,sup}_{{Q}_{\sigma'r_{0}}(z_{0})}f 
		& \leq \  \frac{c(\delta)}{[(\sigma-\sigma')r_{0}]^{\beta_{0}}}\left(\iiint_{{Q}_{\sigma r_{0}}(z_{0})}f^2 \d z\right)^\frac{1}{2}\\
		&\qquad  + c\delta \, \operatorname*{ess\,sup}_{{Q}_{\sigma r_{0}}(z_{0})}f  
		+ \frac{c\delta}{(\sigma-\sigma')^{n+2s+\frac{\beta N_{\rm d}}{p}}  { r_0^{\frac{(\beta - 1)N_d}{p} }}}    \left(\biint_{U_{r_{0}}(z_0)}\tail(f; B_\frac{r_{0}}{2})^{p}\dx\dt\right)^\frac{1}{p}\\
		&{ \leq \  \frac{c(\delta)}{[(\sigma-\sigma')r_{0}]^{\beta_{0}}}\| f \|_{L^\zeta({Q}_{r_{0}}(z_{0}) ) }^{\frac{\zeta}{2}} \left(  \sup_{{Q}_{\sigma r_{0}}(z_{0})}f \right)^{1-\frac{\zeta}{2}}}\\
		&\qquad { + c\delta \, \operatorname*{ess\,sup}_{{Q}_{\sigma r_{0}}(z_{0})}f  
		+ \frac{c\delta}{(\sigma-\sigma')^{n+2s+\frac{\beta N_{\rm d}}{p}}  { r_0^{\frac{(\beta - 1)N_d}{p} }}}    \left(\biint_{U_{r_{0}}(z_0)}\tail(f; B_\frac{r_{0}}{2})^{p}\dx\dt\right)^\frac{1}{p}} \\
		&\  \leq \  \frac{c(\delta, {\zeta} )\|f\|_{L^\zeta({Q}_{r_{0}}(z_{0}))} }{[(\sigma-\sigma')r_{0}]^{\beta_1}}  +\left(c\delta+\frac{2-\zeta}{2}\right)\, \operatorname*{ess\,sup}_{{Q}_{\sigma r_{0}}(z_{0})}f\\
    		& \quad   + \frac{c(\delta)}{(\sigma-\sigma')^{\beta_2}}\,\left(\biint_{U_{r_{0}}(z_0)}\tail(f; B_\frac{r_{0}}{2})^{p}\dx\dt\right)^\frac{1}{p}\,,
        \end{aligned}
    	\end{equation}
by also making use of an application of Young's Inequality (with exponents~$2/\zeta$ and~$2/(2-\zeta)$), where~$\beta_1 \equiv \beta_1(n,s,\Lambda, {\zeta})>0$ and~$\beta_2 \equiv \beta_2(n,s,p)>0$.

    Now, we choose the interpolation parameter~$\delta \in (0,1)$ such that 
       \[
    		c\delta +\frac{2-\zeta}{2}  \ =: \  \eps <1\,, 
       \]
    	which together with~\eqref{eq:harnack-3}  yields
    	\[
        \begin{aligned}
    		\operatorname*{ess\,sup}_{{Q}_{\sigma'r_{0}}(z_{0})} f
    		&  \ \leq  \ \frac{c\,\|f\|_{L^\zeta({Q}_{r_{0}}(z_{0}))} }{(\sigma-\sigma')^{\beta_1}}  +\eps\, \operatorname*{ess\,sup}_{{Q}_{\sigma r_{0}}(z_{0})}f \\
             & \quad + \frac{c\,}{(\sigma-\sigma')^{\beta_2}}\,\left(\biint_{U_{r_{0}}(z_0)}\tail(f; B_\frac{r_{0}}{2})^{p}\dx\dt\right)^\frac{1}{{p}}\notag.
    	\end{aligned}
        \]

    	Now, notice that
	\begin{equation}\label{scelta}
	\Uppsi(\tau) := \operatorname*{ess\,sup}_{Q_{\tau r_{0}}(z_{0})}f 
	\end{equation} 
    is finite for every~$\tau \in [1/2,1]$.

    Thus, 	 a final application of Lemma~\ref{lemma_giusti}, with $\Uppsi$ as in~\eqref{scelta}, $\varrho := \frac 12$, $r:=1$, $A_1:= c\,\|f\|_{L^\zeta({Q}_{r_{0}}(z_{0}))} $ and $A_2:= c\nra{\tail(f; B_\frac{r_{0}}{2})}_{L^{p}(U_{r_{0}}(z_0))}$  yields
         \[
    \operatorname*{ess\,sup}_{{Q}_{\frac{r_{0}}{2}}(z_{0})} f
    		  \ \leq \   c\,\|f\|_{L^\zeta({Q}_{r_{0}}(z_{0}))} + c\,\left(\biint_{U_{r_{0}}(z_0)}\tail(f; B_\frac{r_{0}}{2})^{p}\dx\dt\right)^\frac{1}{p}\,,
        \]
        which  together with the weak Harnack inequality in Theorem~\ref{thm:weak}
        yields the desired strong Harnack inequality in~\eqref{strong} up to relabeling $r_{0}/2$ as~$r_0$. \hfill~$\Box$

\begin{rem}
    Still in theme of Harnack-type inequalities for kinetic equations and $L^2$-$L^\infty$ estimates, it is worth mentioning~{\cite{Loh24c}} in which amongst other results the author proves a strong Harnack inequality and an $L^1$-$L^\infty$ estimate (with no tail) for kinetic integral equations. In order to achieve such results several conditions are imposed on the considered solutions and the author works under a restricted class of nonnegative functions with $L^2$-transport, conservation of mass and polynomial decay in $v$; see Definition~1 there.  The expected nonlocality issues are then partially annihilated by the peculiar framework there, so that no tail contributions do appear both in~\eqref{kinetic_special} and~\eqref{strong}. 
\end{rem}
 
\medskip
  \section{On the sharpness of the $L^p_{t,x}$ tail summability threshold}\label{sec_psharp}
This section is devoted to the construction of the counter-example in Theorem~\ref{thm_sharp}. This is close to the approach used in~\cite{KW24c} for the stationary counterpart. However, the stationary bound only involves the spatial volume~$\eps^{n(1+2s)}$ and leads to the exponent~$n(1+2s)/(2s)$. In the construction below, the time localization contributes the additional factor~$\eps^{2s}$. This gives the drift dimension~$N_{\tt d}=2s+n(1+2s)$ and the threshold~$1+n(1+2s)/(2s)$.

\subsection{The Dirichlet problem and the kinetic semigroup}\label{psharp:subsec_kernel}
Our aim is to prove some 
 properties of the fundamental solution associated to the fractional Kolmogorov equation on a domain~$\mathcal{D} \subset \mathbb{R}^{2n}$, that will be employed later on to build the family~$\{f_\eps\}$ in Theorem~\ref{thm_sharp}.
 
 \vspace{2mm}
 Let~$(\Xi,\mathcal{F},\mathbb{P})$ be a probability space carrying a rotationally symmetric~$2s$-stable process
$\Xb=\{\Xb_t\}_{t\ \geq\ 0}$ in~$\ern$.

We choose a c\`adl\`ag modification of the process, which exists
for a stochastically continuous process with independent and stationary
increments; see, e.~\!g.,~\cite[Chapter~2 and the discussion following Formula~(2.3)]{Sch14}. 
On the completed probability space, there is a single measurable null set~$N$ outside which all its paths are c\`adl\`ag and~$\Xb_0=0$. After setting
the whole path equal to zero on~$N$, for every~$\omega\in\Xi$,  we have that 
\begin{enumerate}[\it (i)]
\item $\Xb_0(\omega)=0$;
\item $t\mapsto\Xb_t(\omega)$ is right-continuous;
\item $\lim_{\tau\to t^-}\Xb_\tau(\omega)$ exists and is finite
for every~$t>0$.
\end{enumerate}

For every point~$(x,v)\in\er^{2n}$, consider now the kinetic trajectory
$(x,v)\in\er^{2n}$ as
\[
\Zs_t^{(x,v)}(\omega) \ := \  \left(x-tv-\int_0^t\Xb_\tau(\omega)\d\tau,v+\Xb_t(\omega)\right) \quad \text{for}~t\ \geq\ 0\,,
\]
and, for any~$\mathcal{D} \subset \er^{2n}$,  we define its first exit
time from~$\mathcal{D}$ by
\[
\tau_{\mathcal{D}}^{(x,v)} \ := \ \inf\left\{\tau>0: \Zs_\tau^{(x,v)}\notin \mathcal{D} \right\},
\]
with the convention~$\inf\varnothing=+\infty$.

We denote with~$P_t^{\mathcal{D}}$ the density associated to~$\Zs^{(x,v)}_t$, killed outside~$\mathcal{D}$ (see, e.~\!g.,~\cite[Section~2]{KW24c}). For every
Borel set~$\Ac \subset \mathcal{D}$,~$t>0$ and~$(x,v)\in \mathcal{D}$, it holds
\begin{equation}\label{eq:density-identity}
\mathbb{P}\left(\left\{t<\tau_{\mathcal{D}}^{(x,v)},
\Zs_t^{(x,v)}\in \Ac\right\}\right)
\ = \ \iint_{\Ac}P_t^{\mathcal{D}}(x,v;y,w)\d y\d w.
\end{equation}

Proceeding similarly as in Section~\ref{sec_preliminaries}, we therefore can define  the kinetic semigroup associated with the fractional Kolmogorov equation in~$\mathcal{D}$ as
\begin{equation}\label{def:kin-semigroup-D}
\mathscr{S}_t^{\mathcal{D}}f(x,v)
 \ := \ \iint_{\mathcal{D}}P_t^{\mathcal{D}}(x,v;y,w)g(y,w)\d y\d w.
\end{equation}

For any~$g \in L^2(\mathcal{D})$, we obtain that the function
\[
\mathscr{P}^{\mathcal{D}}g(t,x,v) = \int_{-\infty}^t [\mathscr{S}_{t-\tau}^\mathcal{D}g](x,v)\d\tau
\]
is the solution of
\[
	\begin{cases}
		\partial_t f+ v\cdot \nabla_x f + (-\Delta_v)^s f &= g \quad \text{in}~\er \times \mathcal{D}\,, \\
		\hfill f&=0 \quad \text{in}~\er \times (\mathbb{R}^{2n} \setminus \mathcal{D}).
	\end{cases}
\]

\vspace{2mm}
In order to prove Theorem~\ref{thm_sharp} we will employ several properties of~$P^\mathcal{D}_t$ which we are going to establish.

 First of all, we note that for any two bounded open sets~$\mathcal{D}\subset\mathcal{D}' \subset \mathbb{R}^{2n}$, every Borel set
$\mathcal{A}\subset\mathcal{D}$, every~$t>0$, and every
$(x,v)\in\mathcal{D}$, one has
\begin{equation}\label{eq:maximum-principle}
0 \ \leq\ \iint_{\Ac}P_t^{\mathcal{D}}(x,v;y,w)\d y\d w \ \leq\ \iint_{\Ac}P_t^{\mathcal{D}'}(x,v;y,w)\d y\d w 
\ \leq\ \iint_{\Ac}P_t(x-y-tw,v-w)\d y\d w\,,
\end{equation}
where~$P_t$ is the fundamental solution in~\eqref{eq:fundamental}; see, e.~\!g.,~\cite[Formula~(2.4)]{KW24c}.

Next, for~$r>0$, we define the phase dilation
\[
\delta_r^0(x,v):=(r^{1+2s}x,rv).
\]
Then,~$P^{\mathcal{D}}_t$ satisfies the following scaling property for~$(0,0) \in \mathcal{D}$
\begin{equation}\label{eq:fundamental-dilation}
P_{\eps^{2s}t}^{\mathcal{D}}
 (0,0;\eps^{1+2s}y,\eps w) \ = \ \eps^{-n(2+2s)}
P_t^{\delta_\frac{1}{\eps}^0(\mathcal{D})}(0,0;y,w).
\end{equation}
Here~$t,\eps>0$, the initial point~$(x,v)$ belongs to
$\delta_{1/\eps}^0(\mathcal{D})$, and the identity holds for almost every
final point~$(y,w)$ in that domain; see, e.~\!g.,~\cite[Formula~(2.5)]{KW24c}.

Now, as consequence of~\eqref{eq:maximum-principle} we also deduce an almost everywhere control of~$P^{\mathcal{D}}_t$ in terms of~$P_t$.

\begin{prop}\label{prop:ae-max-principle}
   For~$t>0$,~$(x,v)\in\mathcal{D}$, 
  \[
  P_t^{\mathcal{D}}(x,v;y,w) \ \leq\ P_t(x-y-tw,v-w) \quad \text{for a.~\!e.}~(y,w) \in \mathcal{D}.
  \]
Moreover, for each fixed~$t>0$, this inequality holds for almost every
$(x,v,y,w)\in\mathcal{D}\times\mathcal{D}$.
\end{prop}

\begin{proof}
Let us fix~$t>0$ and~$(x,v)\in\mathcal{D}$. Define
\[ 
F_t(y,w) \ := \ P_t^{\mathcal{D}}(x,v;y,w)
- P_t(x-y-tw,v-w)\,,
\]
and, for every~$m\in\mathbb{N}$, let 
\[
\mathcal{A}_m \ := \ 
\left\{(y,w)\in\mathcal{D}: F_t(y,w)\ \geq\ \frac{1}{m}\right\}.
\]

We apply~\eqref{eq:maximum-principle} on~$\mathcal{A}_m$, so that recalling the definition of~$\Ac_m$ gives
\[
\frac{|\Ac_m|}{m} \ \leq\ \iint_{\mathcal{A}_m}F_t(y,w)\d y\d w \ \leq\ 0.
\]

Then,~$|\Ac_m|=0$ for any~$m \in \en$. Hence, since
\[
|\left\{(y,w)\in\mathcal{D}:F_t(y,w)>0\right\}|
 \ \leq\ \sum_{m=1}^{+\infty}|\mathcal{A}_m|=0,
\]
we have that~$F_t(y,w)\ \leq\  0$ a.~\!e. on~$\mathcal{D}$. This gives the first assertion. By joint measurability the set where the
inequality fails is measurable in~$\mathcal{D}\times\mathcal{D}$.
Its section in~$(y,w)$ has measure zero for every fixed~$(x,v)$.
Tonelli's theorem applied to its indicator therefore gives product
measure zero, proving the last assertion.
\end{proof}

Lastly, we establish a lower estimate for~$P_t^{B_2\times B_2}$ on a sufficiently small subdomain.

\begin{lemma}\label{psharp:lem_kernel_mass}
There exist~$T \in (0,1]$, depending only on~$n$ and~$s$, such that
\begin{equation}\label{psharp:eq_kernel_mass}
\inf_{(t,x,v)\in [T/4,T]\times B_{\frac{1}{16}}\times B_{\frac{1}{16}}}\iint_{B_1\times B_1}
P_t^{B_2\times B_2}(x,v;y,w)\d y\d w \ \geq\ \frac{1}{2}.
\end{equation}
\end{lemma}

\begin{proof}
For~$T>0$, consider the event~$E_T \in \mathcal{F}$ given by
\[
E_T \ := \ 
\left\{\omega\in\Xi:
\sup_{t\in[0,T]}|\Xb_t(\omega)|<\frac14
\right\}.
\]

For any~$m \in \en$, let
\[
T_m \ := \ \frac{1}{m}.
\]
Since the trajectories are right-continuous at~$t=0$ and
$\Xb_0=0$, for every~$\omega\in\Xi$ there exists
$\delta_\omega>0$ such that
\[
|\Xb_t(\omega)|<\frac18
\qquad
\text{for every }t\in[0,\delta_\omega).
\]
Consequently, for every~$\omega\in\Xi$ there exists~$m(\omega)\in\mathbb{N}$ such that~$
\omega\in E_{T_m}$ for every~$m\ \geq\  m(\omega)$.

It follows that
\[
\lim_{m\to +\infty}\mathbbm{1}_{E_{T_m}}(\omega)
\ = \ 1 \quad \forall \omega \in \Xi\,,
\]
and, by the Dominated Convergence Theorem, we obtain
\[
\lim_{m\to+\infty}\mathbb{P}(E_{T_m})
\ = \ 1.
\]

Therefore, we choose~$m_0\in\mathbb{N}$ such that
\[
\mathbb{P}(E_{T_{m_0}})\ \geq\ \frac{1}{2}\,,
\]
and set~$T  :=  T_{m_0}  =  \frac1{m_0} \in (0,1]$.

We remark that, since the probability law of the normalized rotationally symmetric
$2s$-stable process is determined by~$n$ and~$s$, the
probabilities~$\mathbb{P}(E_{T_m})$, and therefore the choice of
$m_0$ and~$T$, depend only on~$n$ and~$s$.

Let now~$(x,v)\in B_{\frac{1}{16}}\times B_{\frac{1}{16}}$ 
and assume that the event~$E_T$ occurs. For every~$t\in[0,T]$, since~$T\ \leq\ 1$, we have 
\[
|v+\Xb_t(\omega)| 
\ < \ \frac{5}{16} \quad \text{and} \quad 
\left|x-tv-\int_0^t\Xb_\tau(\omega)\d\tau\right| \ <\ \frac{3}{8}. 
\]

Thus, on~$E_T$, we have that~$Z_t^{(x,v)}\in B_{\frac{3}{8}}\times B_{\frac{5}{16}}$, for every~$t \in [0,T]$, and in particular~$Z_t^{(x,v)}$ does not leave~$B_2\times B_2$ before time~$T$. 

Let us verify now that the trajectory does not exit~$B_2\times B_2$ at time~$T$ either. Since the map~$t \mapsto Z_t^{(x,v)}(\omega)$ is right-continuous and~$Z_T^{(x,v)}(\omega) \in B_{\frac38}\times B_{\frac5{16}}$, there exists~$\varepsilon_\omega>0$ such that
\[
Z_t^{(x,v)}(\omega)\in B_2\times B_2
\quad \text{for every}~t\in[T,T+\varepsilon_\omega)\,,
\]
i.~\!e.~$\tau_{B_2\times B_2}^{(x,v)}>T$.

Therefore, for every~$t\in[T/4,T]$,
\begin{equation}\label{eq:ET-inclusion}
E_T \ \subset \ \left\{t<\tau_{B_2\times B_2}^{(x,v)}, Z_t^{(x,v)}\in B_1\times B_1 \right\}.
\end{equation}

By~\eqref{eq:density-identity} and~\eqref{eq:ET-inclusion}, we conclude
that
\[ 
\begin{aligned}
\iint_{B_1\times B_1}P_t^{B_2\times B_2}(x,v;y,w)\d y\d w & \ = \ \mathbb{P}\left(\left\{t<\tau_{B_2\times B_2}^{(x,v)},\ Z_t^{(x,v)}\in B_1\times B_1\right\}\right) \\
& \ \geq\ \mathbb{P}(E_T) \ \geq\ \frac{1}{2}. 
\end{aligned}
\]
for every~$(t,x,v)\in[T/4,T]\times B_{\frac{1}{16}}\times B_{\frac{1}{16}}$. Thus, the desired~\eqref{psharp:eq_kernel_mass} follows.
\end{proof}

The following properties of the semigroup~$\mathscr{S}_t^{\mathcal{D}}$ hold true.

\begin{prop}\label{prop:L1-Linf-contraction}
The following properties hold true:
\begin{enumerate}[\it (i)]
     \item For any~$g \in L^\infty(\mathcal{D})$, we have that
     \[
      \|\mathscr{S}_t^{\mathcal{D}}g\|_{L^\infty(\mathcal{D})}
      \ \leq\ \|g\|_{L^\infty(\mathcal{D})}.
     \]
    \item For any~$g \in L^1(\mathcal{D})$, we have that
     \[
      \|\mathscr{S}_t^{\mathcal{D}}g\|_{L^1(\mathcal{D})}
      \ \leq\ \|g\|_{L^1(\mathcal{D})}.
     \]
\end{enumerate}
\end{prop}

\begin{proof}
The case~$t=0$ follows by defining~$\mathscr{S}_0^{\mathcal{D}}$ as the
identity. Fix~$t>0$. Nonnegativity of the kernel gives
\[
|\mathscr{S}_t^{\mathcal{D}}g|
\ \leq\ \mathscr{S}_t^{\mathcal{D}}|g|.
\]
For~$g\in L^\infty(\mathcal{D})$, by~\eqref{eq:maximum-principle},
\[
\begin{aligned}
|\mathscr{S}_t^{\mathcal{D}}g(x,v)|
&\ \leq\ \|g\|_{L^\infty(\mathcal{D})}
 \iint_{\mathcal{D}}P_t^{\mathcal{D}}(x,v;y,w)\d y\dw
\\
&\ \leq\ \|g\|_{L^\infty(\mathcal{D})}
 \iint_{\mathbb{R}^{2n}}P_t(x-y-tw,v-w)\d y\dw
=\|g\|_{L^\infty(\mathcal{D})}.
\end{aligned}
\]
where in the last equality we used the change of variables~$X=x-y-tw$ and~$V=v-w$ and  Proposition~\ref{prop:fund}~\textit{(iii)}.
This proves~\textit{(i)}.

For~\textit{(ii)}, the joint almost everywhere comparison in
Proposition~\ref{prop:ae-max-principle} gives, for almost every~$(y,w)$,
\[
\begin{aligned}
\iint_{\mathcal{D}}P_t^{\mathcal{D}}(x,v;y,w)\dx\dv
&\ \leq\ \iint_{\mathbb{R}^{2n}}P_t(x-y-tw,v-w)\dx\dv
=1.
\end{aligned}
\]
By Fubini-Tonelli's Theorem, we this obtain
\[
\begin{aligned}
\|\mathscr{S}_t^{\mathcal{D}}g\|_{L^1(\mathcal{D})}
&\ \leq\ \iint_{\mathcal{D}}|g(y,w)|
 \left(\iint_{\mathcal{D}}P_t^{\mathcal{D}}(x,v;y,w)\dx\dv\right)
 \d y\dw
\\
&\ \leq\ \iint_{\mathcal{D}}|g(y,w)|\d y\dw.
\end{aligned}
\]
\end{proof}

\subsection{Construction of the family~$\{f_\eps\}$}\label{psharp:subsec_solution}
In this section we provide the explicit construction of the function~$f_\eps$.
Fix once and for all two cut-off functions
\[
\eta\in C_c^\infty(B_2)
\quad\text{and}\quad
\psi\in C_c^\infty(B_\frac{1}{2}(3e_n)),
\]
such that
\[
\eta=1\quad\text{in }B_1,
\qquad 0\ \leq\ \eta,\psi\ \leq\ 1,
\qquad \int_{\ern}\psi(w)\dw>0.
\]
For every~$\eps\in(0,1)$, define
\[
\eta_\eps(x) \ := \ \eta\left(\frac{x}{\eps^{1+2s}}\right),
\quad \text{and} \quad 
\beta_\eps(x,v) \ := \ \eta_\eps(x)\psi(v).
\]
Since~$\operatorname{supp}\psi\Subset B_\frac{1}{2}(3e_n)$, one has
\[
\operatorname{dist}(\operatorname{supp}\psi,\overline{B_2})>0,
\qquad
\operatorname{supp}\psi\subset\ern\setminus\overline{B_2}.
\]
Thus~$\beta_\eps$ vanishes in a neighborhood of every point of
$B_2\times\overline{B_2}$. 

For~$(x,v)\in B_2\times B_2$, set
\begin{equation}\label{def:g-eps}
g_\eps(x,v) \ := \ -(-\Delta_v)^s\beta_\eps(x,v).
\end{equation}

Since~$\psi$ vanishes in a neighborhood of~$\overline{B_2}$, the preceding formula gives
\[
g_\eps(x,v) \ = \ 
 c_{n,s}\eta_\eps(x)
\int_{\mathbb{R}^n}
\frac{\psi(w)}{|v-w|^{n+2s}}\dw \quad 
\text{for}~(x,v)\in B_2\times B_2.
\]
In particular, there exists two constants~$0<c_1\ \leq\  c_2$, depending only on~$n$,~$s$ and~$\psi$, such that
\begin{equation}\label{psharp:eq_gKW_comp}
 c_1\mathbbm{1}_{B_{\eps^{1+2s}}\times B_2} \ \leq\  \
 g_\eps \ \leq\ c_2\mathbbm{1}_{B_{2\eps^{1+2s}}\times B_2}
\quad \text{in}~B_2 \times B_2\,,
\end{equation}
and, consequently,
\begin{equation}\label{psharp:eq_gKW_bounds}
\left\|g_\eps\right\|_{L^\infty(B_2 \times B_2)}
\ \leq\ c \quad \text{and} \quad 
\left\|g_\eps\right\|_{L^1(B_2\times B_2)} \ \leq\ c\eps^{n(1+2s)}.
\end{equation}

Now, let~$T\in(0,1]$ be the constant given by Lemma~\ref{psharp:lem_kernel_mass}.  Consider
\[
a\in C_c^\infty((-2T,-T/4))\,,
\]
such that
\[
a \ = \ 1
\quad
\text{on}~[-T,-T/2] \quad \text{and} \quad 0 \ \leq\ a \ \leq\  1.
\]

For~$\eps\in(0,1)$ define
\[
a_\eps(t) \ := \ a\left(\frac{t}{\eps^{2s}}\right).
\]

Define
\begin{equation}\label{def:b-eps}
b_\eps(t,x,v) \ := \ a_\eps(t)\beta_\eps(x,v)
\quad
\text{in}~\er^{1+2n}\,,
\end{equation}
and, for~$z=(x,v)\in B_2\times B_2$, let
\begin{equation}\label{psharp:eq_correction}
w_\eps^{\rm corr}(t,x,v) \ := \ \int_{-\infty}^ta_\eps(\tau)\left[\mathscr{S}_{t-\tau}^{B_2\times B_2}g_\eps\right](z)\d\tau\,,
\end{equation}
extended to zero in~$\er\times\big(\er^{2n}\setminus(B_2\times B_2)\big)$.
Let
\begin{equation}\label{def:f-eps}
f_\eps \ := \ b_\eps+w_\eps^{\rm corr} \quad \text{in}~\er^{1+2n}.
\end{equation}

\begin{prop}\label{prop:solution}
Let~$\mathcal{Q}:= (-2,1/2) \times B_2 \times B_2$. Then, there exists~$\eps_0 \in (0,1)$, depending only on~$n$ and~$s$, such that, for every~$\eps\in(0,\eps_0)$, the function~$f_\eps$  given in~\eqref{def:f-eps} is a globally nonnegative weak solution to
\[
(\partial_t+v\cdot\nabla_x)f_\eps+(-\Delta_v)^sf_\eps=0
\quad\text{in } \mathcal{Q}\,,
\]
according to  Definition~{\rm\ref{def_weak_sol}}.

Moreover, for every~$t\in\er$, it holds
\begin{equation}\label{psharp:eq_correction_bounds}
\begin{aligned}
\|w_\eps^{\rm corr}(t)\|_{L^\infty(B_2\times B_2)} & \ \leq\ c\eps^{2s},\\
\|w_\eps^{\rm corr}(t)\|_{L^1(B_2\times B_2)} &\ \leq\ c\,\eps^{2s+n(1+2s)}\,,
\end{aligned}
\end{equation}
for some~$c>0$ depending only on~$n$,~$s$ and the fixed functions~$\eta$,~$\psi$ and~$a$, but not on~$\eps$.
\end{prop}

\begin{proof}
We first justify the equation and the energy regularity of the correction.
By Proposition~\ref{prop:L1-Linf-contraction} and Riesz--Thorin's Theorem (with
$1/2=(1-1/2)/1+(1/2)/\infty$), we have that~$\mathscr{S}_t^{B_2\times B_2}$ is a contruction on~$L^2(B_2\times B_2)$. Moreover, as in~\cite[Section~2]{KW24c},~$w_\eps^{\rm corr}$ satisfies (in the distributional sense)
the following conditions
\[
\begin{cases}
(\partial_t+v\cdot\nabla_x)w_\eps^{\rm corr}
 +(-\Delta_v)^sw_\eps^{\rm corr}=a_\eps(t)g_\eps
 &\text{in }\er\times B_2\times B_2,\\
w_\eps^{\rm corr}(t,x,v)=0
 &\text{in }\er
 \times (\er^{2n}\setminus(B_2\times B_2)).
\end{cases}
\]

We now prove that~$w_\eps^{\rm corr} \in \Wc(\mathcal{Q})$. 
Fix~$t\in\er$,~$\delta>0$ and split the integral in~\eqref{psharp:eq_correction} via the semigroup property of~$\mathscr{S}^{B_2\times B_2}_\tau$
\[
w_\eps^{\rm corr}(t+\delta)
\ = \ \mathscr{S}_{\delta}^{B_2\times B_2}w_\eps^{\rm corr}(t)
 +\int_t^{t+\delta}a_\eps(\tau)
 \mathscr{S}_{t+\delta-\tau}^{B_2\times B_2}g_\eps\d\tau.
\]
Then, we obtain that
\begin{equation}\label{eq:correct-quotient}
\frac{\mathscr{S}_{\delta}^{B_2\times B_2}w_\eps^{\rm corr}(t)
             -w_\eps^{\rm corr}(t)}\delta
=
 \frac{w_\eps^{\rm corr}(t+\delta)-w_\eps^{\rm corr}(t)}\delta
 -\bint_0^{\delta} a_\eps(t+\delta-\sigma)
 \mathscr{S}_{\sigma}^{B_2\times B_2}g_\eps\,\d \sigma.
\end{equation}

Let us begin noticing that changing~$\sigma=t-\tau$ in~\eqref{psharp:eq_correction} yields
\[
w_\eps^{\rm corr}(t)
\ =\ \int_0^\infty
 a_\eps(t-\sigma)\mathscr{S}_{\sigma}^{B_2\times B_2}g_\eps\,\d \sigma\,,
\]
so that, by letting
\[
\partial_tw_\eps^{\rm corr}(t)
\ :=\ \int_0^\infty
 a_\eps'(t-\sigma)\mathscr{S}_{\sigma}^{B_2\times B_2}g_\eps\d \sigma\,,
\]
we obtain from the Fundamental Theorem of Calculus, Fubini's Theorem and  the $L^2$-contraction property
\[
\begin{aligned}
&\left\|\frac{w_\eps^{\rm corr}(t+\delta)-w_\eps^{\rm corr}(t)}\delta-\partial_tw^{\rm corr}_\eps(t)\right\|_{L^2(B_2\times B_2)}\\
& \quad \leq \int_0^\infty \left| \frac{a_\eps(t+\delta-\sigma)-a_\eps(t-\sigma)}{\delta}-a_\eps'(t-\sigma)\right| \|\mathscr{S}_{\sigma}^{B_2\times B_2}g_\eps\|_{L^2(B_2\times B_2)}\d\sigma\\
& \quad\ \leq\ \|g_\eps\|_{L^2(B_2\times B_2)}\int_0^1
 \|a_\eps'(\cdot+\theta \delta)-a_\eps'\|_{L^1(\er)}\,\d\theta\\
& \quad\ \leq\ \frac{|\delta|}{2}\|g_\eps\|_{L^2(B_2\times B_2)}
 \|a_\eps''\|_{L^1(\er)}\ \to\ 0.
\end{aligned}
\]
Thus the first term in~\eqref{eq:correct-quotient} converges strongly to~$\partial_tw_\eps^{\rm corr}(t)$. 

Moreover, by combining once again the $L^2$-contraction with the Fundamental Theorem of Calculus we have
\[
\begin{aligned}
\|\partial_tw_\eps^{\rm corr}(t+\delta)-\partial_tw_\eps^{\rm corr}(t)\|_{L^2(B_2\times B_2)}
&\ \leq\ \|g_\eps\|_{L^2(B_2\times B_2)}
 \|a_\eps'(\cdot+\delta)-a_\eps'\|_{L^1(\er)}\\
&\ \leq\ |\delta|\|g_\eps\|_{L^2(B_2\times B_2)}\|a_\eps''\|_{L^1(\er)} \to 0\,, \quad \text{as}~\delta\to 0^+.
\end{aligned}
\]
Consequently,
\begin{equation}\label{eq:time-reg-corrector}
w_\eps^{\rm corr}\in C^1(\er;L^2(B_2\times B_2)).
\end{equation}

Now, the second term in~\eqref{eq:correct-quotient} satisfies
\[
\begin{aligned}
&\left\|\bint_0^{\delta} a_\eps(t+\delta-\sigma)
 \mathscr{S}_{\sigma}^{B_2\times B_2}g_\eps\,\d \sigma-a_\eps(t)g_\eps\right\|_{L^2(B_2\times B_2)}\\
&\quad\ \leq\ \delta\|a_\eps'\|_{L^\infty(\er)}\|g_\eps\|_{L^2(B_2\times B_2)}
 +|a_\eps(t)|\sup_{0\leq \sigma\leq \delta}
 \|\mathscr{S}_{\sigma}^{B_2\times B_2}g_\eps-g_\eps\|_{L^2(B_2\times B_2)}
 \ \to\ 0\,,
\end{aligned}
\]
where in the last line we have used the $L^2$-contraction and right-continuity of the killed process together with the Dominated Convergence Theorem.

All in all, it follows that, as~$\delta\to 0^+$, we have
\begin{equation}\label{eq:corrector-quotient-limit}
\frac{\mathscr{S}_{\delta}^{B_2\times B_2}w_\eps^{\rm corr}(t)
             -w_\eps^{\rm corr}(t)}\delta
\ \to\
\partial_tw_\eps^{\rm corr}(t)-a_\eps(t)g_\eps
\quad\text{in}~L^2(B_2\times B_2).
\end{equation}

The convergence~\eqref{eq:corrector-quotient-limit} yields 
\begin{multline}\label{eq:corrector-limit}
\lim_{\delta\to 0^+}\frac{\|w_\eps^{\rm corr}(t)\|_{L^2(B_2\times B_2)}^2
 -\|\mathscr{S}_{\delta}^{B_2\times B_2}w_\eps^{\rm corr}(t)\|_{L^2(B_2\times B_2)}^2}{2\delta} \\
 = - \langle \partial_t w_\eps^{\rm corr}(t)-a_\eps(t)g_\eps, w^{\rm corr}_\eps(t)\rangle_{L^2(B_2\times B_2)}.
\end{multline}

Now, we can establish the energy regularity of~$w_\eps^{\rm corr}$. For this, let us recall that, for~$g\in L^2(\er^{2n})$, the Fourier formula for the kinetic semigroup and Plancherel's Identity give
\[
\|\mathscr{S}_tg\|_{L^2(\er^{2n})}^2 \ = \ \frac1{(2\pi)^{2n}}\iint_{\ern\times\ern}
 e^{-2\int_0^t|\eta-\tau\xi|^{2s}\d \tau}|\widehat g(\xi,\eta)|^2\d\xi\d\eta.
\]
Moreover, since as~$ t\to 0^+$, we have
\[
\frac{1-e^{-2\int_0^t|\eta-\tau\xi|^{2s}\d \tau}}{2t}
\to  |\eta|^{2s},
\]
Moreover, by Fatou's lemma and a combination of Fubini-Tonelli's Theorem with Plancherel's Identity, we obtain that
\[
\begin{aligned}
& \liminf_{t \to 0^+}\frac{\|g\|_{L^2(\er^{2n})}^2 -\|\mathscr{S}_t g\|_{L^2(\er^{2n})}^2}{2t}\\
& \quad  \ \geq\ \frac1{(2\pi)^{2n}}
 \iint_{\ern\times\ern}|\eta|^{2s}|\widehat{g}(\xi,\eta)|^2\d\xi\d\eta
 \ = \ \frac{c_{n,s}}2\int_{\ern}[g]^2_{H^s(\ern)}\dx.
\end{aligned}
\]

Fix~$t\in \er$, upon extending~$w_\eps^{\rm corr}(t,\cdot,\cdot)$ to zero outside its support, by Proposition~\ref{prop:ae-max-principle} we obtain
\[
\|\mathscr{S}_\delta^{B_2\times B_2}w_\eps^{\rm corr}(t,\cdot,\cdot)\|_{L^2(B_2\times B_2)}^2
\ \leq\ \|\mathscr{S}_\delta w_\eps^{\rm corr}(t,\cdot,\cdot)\|_{L^2(\er^{2n})}^2.
\]
Using~\eqref{eq:corrector-limit} above, we obtain
\[
\begin{aligned}
&-\langle\partial_tw_\eps^{\rm corr}(t)-a_\eps(t)g_\eps,
                  w_\eps^{\rm corr}(t)\rangle_{L^2(B_2\times B_2)}
\\
&\qquad=
\lim_{\delta \to 0^+}
\frac{\|w_\eps^{\rm corr}(t)\|_2^2
 -\|\mathscr{S}_\delta^{B_2\times B_2}w_\eps^{\rm corr}(t)\|_{L^2(B_2\times B_2)}^2}{2\delta}
\\
&\qquad\ \geq\ \frac{c_{n,s}}2\int_{B_2} \iint_{\ern\times\ern}
 \frac{|w_\eps^{\rm corr}(t,x,v)-w_\eps^{\rm corr}(t,x,w)|^2}
 {|v-w|^{n+2s}}\dw\dv\dx.
\end{aligned}
\]

Since~\eqref{eq:time-reg-corrector} gives~$\langle\partial_tw_\eps^{\rm corr},w_\eps^{\rm corr}\rangle_{L^2(B_2\times B_2)}
=\frac12\partial_t\|w_\eps^{\rm corr}\|_{L^2(B_2\times B_2)}^2$, we can test the equation with~$w_\eps^{\rm corr}(t)\ \geq\ 0$. Hence, integrating from~$t_1> - T\eps^{2s}/2$ gives
\[
\begin{aligned}
&\frac12\|w_\eps^{\rm corr}(t_1)\|_{L^2(B_2\times B_2)}^2
+\frac{c_{n,s}}2
 \int_{-\infty}^{t_1}\int_{B_2}
 \iint_{\ern\times\ern}
 \frac{|w_\eps^{\rm corr}(t,x,v)-w_\eps^{\rm corr}(t,x,w)|^2}
 {|v-w|^{n+2s}}\dw\dz
\\
&\quad\ \ \leq\ \left(\int_{\er}a_\eps(t)\|g_\eps\|_{L^2(B_2\times B_2)}\dt\right)^2.
\end{aligned}
\]
by the $L^2$-contraction properties of~$\mathscr{S}^{B_2\times B_2}_t$. 

Thus,
\[
w_\eps^{\rm corr}\in
 C^1(\er;L^2(B_2\times B_2)) \cap L^2_{\rm loc}(\er\times B_2;H^s(\ern)).
\]

Moreover, for every compact~$K\Subset\mathbb{R}\times B_2$, H\"older's Inequality gives
\[
\int_K\int_{\ern}
\frac{|w_\eps^{\rm corr}(t,x,v)|}{(1+|v|)^{n+2s}}\dv\dx\dt
\  \leq\ |B_2|^{1/2}|K|^{1/2}
\|w_\eps^{\rm corr}\|_{L^2(K\times B_2)}<+\infty.
\]

Now, we are left with computing the transport derivative on~$w_\eps^{\rm corr}$. For this, let us choose~$\eps_0\in(0,1)$, depending only on~$n,s$, such that
\begin{equation}\label{eq:corr-eps-condition}
2\eps_0^{1+2s}+4T\eps_0^{2s}<1,
\end{equation}
and fix~$\eps \in (0,\eps_0)$.

Extend~$g_\eps$ to~$\ern\times B_2$ by
its defining formula, which is smooth and supported in~$B_{2\eps^{1+2s}}$ in~$x$. Let us consider the exit time in~$B_2$:
\[
\tau_{B_2}^{v}:=\inf\{\tau>0:v+\Xb_{\tau}\notin B_2\}.
\]

For~$\sigma\geq0$, if~$a_\eps(t-\sigma)\neq0$, then~$\sigma \in (T\eps^{2s}/4+t,t+2T\eps^{2s})$. Moreover, let us note that on~$\{\sigma<\tau_{B_2}^{v}\}$, all velocities up to~$\sigma$ have modulus less than two. Also, if
\[
g_\eps\left(x-\sigma{}v-\int_0^{\sigma}\Xb_\tau\,\d\tau,
             v+\Xb_{\sigma}\right)\neq0,
\]
then, for every~$0\leq\ell\leq \sigma$, since~$\eta_\eps$ is supported on~$2\eps^{1+2s}$, by~\eqref{eq:corr-eps-condition} we obtain
\[
\begin{aligned}
\left|x-\ell v-\int_0^\ell\Xb_\tau\,\d\tau\right|
&\ \leq\
\left|x-\sigma{}v-\int_0^{\sigma}\Xb_\tau\,\d\tau\right|
 +\int_\ell^{\sigma}|v+\Xb_\tau|\,\d\tau\\
&\ <\ 2\eps^{1+2s}+2\sigma
\ <\ 2\eps^{1+2s}+4T\eps^{2s} + 2t<2\,,
\end{aligned}
\]
for any~$t \in (-2,\frac{1}{2})$.

Thus, for~$(x,v)\in B_2\times B_2$, the probabilistic representation becomes
\[
w_\eps^{\rm corr}(t,x,v)
=\int_0^{+\infty}a_\eps(t-\sigma)\,
\mathbb{E}\!\left[
 g_\eps\left(x-\sigma{}v-\int_0^{\sigma}\Xb_\tau\,\d\tau,
             v+\Xb_{\sigma}\right)
 \mathbbm{1}_{\{\sigma<\tau_{B_2}^{v}\}}
\right]\d \sigma,
\]
where integrand is defined to be zero when~$\sigma\geq\tau_{B_2}^v$.
Its time support is contained in
$[0,\infty)\cap(t+T\eps^{2s}/4,t+2T\eps^{2s})$.

Since the indicator is independent of~$x$ and
\[
\|\nabla_xg_\eps\|_{L^\infty(\ern\times B_2)}
\ \leq\ c\eps^{-1-2s},
\]
up to differentiation under the integral, by Dominate Convergence Theorem, we obtain from the $L^\infty$-contraction in Proposition~\ref{prop:L1-Linf-contraction} that
\[
\|\nabla_xw_\eps^{\rm corr}(t)\|_{L^\infty(B_2\times B_2)}
\ \leq\ \|\nabla_xg_\eps\|_{L^\infty(\ern\times B_2)}
 \|a_\eps\|_{L^1(\er)}\ \leq\ c\eps^{-1}\,,
\]
for~$t \in (-2,1/2)$.

This, together with~\eqref{eq:time-reg-corrector}, yields
\[
\begin{aligned}
(\partial_t+v\cdot\nabla_x)w_\eps^{\rm corr}
&\in L^2_{\rm loc}((-2,1/2)\times B_2;L^2(\ern))\\
&\subset L^2_{\rm loc}((-2,1/2)\times B_2;H^{-s}(\ern)).
\end{aligned}
\]
Therefore~$w_\eps^{\rm corr}\in\Wc(\mathcal{Q})$.

Now we compute the equation satisfied by~$b_\eps$ in~$\er\times B_2\times B_2$. Note that, since~$\psi$ vanishes in a neighborhood of~$\overline{B_2}$, one has
\begin{equation}\label{eq:trasport-beps}
b_\eps \ = \  \partial_tb_\eps \ = \ v\cdot\nabla_xb_\eps \ = \ 0
\quad
\text{in}~\er \times B_2\times B_2.
\end{equation}
Moreover, by the definition of~$g_\eps$ in~\eqref{def:g-eps}, we have that
\begin{equation}\label{eq:fraclap-beps}
(-\Delta_v)^sb_\eps(t,x,v) \ = \  a_\eps(t)(-\Delta_v)^s\beta_\eps(x,v) \ = \ 
- a_\eps(t)g_\eps(x,v) \quad \text{in}~\er \times B_2 \times B_2.
\end{equation}

Therefore, by~\eqref{def:f-eps} we obtain that, in the distributional sense, it holds
\[
(\partial_t+v\cdot\nabla_x)f_\eps+(-\Delta_v)^sf_\eps \ = \ 0
\quad
\text{in}~(-2,1/2)\times B_2\times B_2\,,
\]
and, by~\eqref{eq:trasport-beps},~\eqref{eq:fraclap-beps} and since~$w_\eps^{\it corr} \in \Wc$, we also have 
that~$f_\eps\in\Wc(\mathcal{Q})$.. 

Lastly, since the distributional identity in~$B_2\times B_2$
extends to the tests of Definition~\ref{def_weak_sol} by density and
continuity of the~$H^{-s}$-$H^s$ pairing and of the fractional bilinear
form, we have that~$f_\eps$ is a weak solution to~\eqref{problema} in~$\mathcal{Q}$ according to Definition~\ref{def_weak_sol}.

We now prove~\eqref{psharp:eq_correction_bounds}.  By Proposition~\ref{prop:L1-Linf-contraction}~{\it (i)} and~\eqref{psharp:eq_gKW_bounds} we obtain
\[
0\ \leq\ w_\eps^{\rm corr}(t,x,v)
\ \leq\ \|g_\eps\|_{L^\infty(B_2\times B_2)}
 \int_{-\infty}^ta_\eps(\tau)\d\tau \ \leq\ c\eps^{2s}.
\]
Thus~$\|w_\eps^{\rm corr}(t)\|_{L^\infty(B_2\times B_2)}\ \leq\  c\eps^{2s}$.

 For the~$L^1$-bound, we apply Fubini-Tonelli's Theorem, Proposition~\ref{prop:L1-Linf-contraction}~{\it (ii)} and~\eqref{psharp:eq_gKW_bounds}, so that 
\[
\begin{aligned}
\|w_\eps^{\rm corr}(t)\|_{L^1(B_2\times B_2)}
& \ = \ \int_{-\infty}^ta_\eps(\tau)
 \|\mathscr{S}_{t-\tau}^{B_2\times B_2}g_\eps\|_{L^1(B_2\times B_2)}\d\tau
\\
&\ \leq\ \|g_\eps\|_{L^1(B_2\times B_2)}
 \int_{\mathbb{R}}a_\eps(\tau)\d\tau
\\
&\ \leq\ c\eps^{2s+n(1+2s)}.
\end{aligned}
\]
This concludes the proof. 
\end{proof}

\subsection{Proof of Theorem~\ref{thm_sharp}}\label{psharp:sec_sharpi}
We now divide the proof depending on the range of the parameter~$p$.

\vspace{2mm} 
{\it Case 1: Assume~$p \in [1,p^\ast)$}.
Consider the function~$f_\eps$ given in~\eqref{def:f-eps}. By Proposition~\ref{prop:solution}, we obtain that it is a weak solution to~\eqref{problema} in~$\mathcal{Q}$. 
Now, we show  that there exist~$c,C>0$ depending only on~$n$ and~$s$, such that, for any~$\eps \in (0,\eps_0)$, it holds
\begin{align}
f_\eps(0) & \ \geq \ c\,\eps^{2s}\,, \label{psharp:eq_lower_set}\\
\nra{\tail(f_\eps;B_\frac{1}{2})}_{L^p(U_1)}
& \ \leq  \ C\eps^\frac{N_{\tt d}}{p} \quad \forall p \in (1,\infty)\,, \label{psharp:eq_tail_estimate}\\
\text{and} \quad 
\nra{f_\eps}_{L^2(Q_1)}
& \ \leq \ C\eps^\frac{N_{\tt d}+2s}{2}.\label{psharp:eq_L2_estimate}
\end{align}

Let us begin noticing that by definition,~$b_\eps \equiv 0$ in~$\mathcal{Q}$. Therefore, by the formula~\eqref{psharp:eq_correction},  together with the lower bound in~\eqref{psharp:eq_gKW_comp} , the definition of~$a_\eps$  and~\eqref{psharp:eq_gKW_comp} we obtain
\[
\begin{aligned}
f_\eps(0)
& \ \geq \ \int_{-\infty}^0a_\eps(\tau)\left(\iint_{B_2\times B_2}P^{B_2\times B_2}_{-\tau}(0,0;y,w)g_\eps(y,w)\d y \dw\right) \d\tau\\
& \ \geq  \ c_1
\int_{\frac{T\eps^{2s}}{2}}^{T\eps^{2s}}
\int_{B_{\eps^{1+2s}}}\int_{B_\eps}
 P_{\tau}^{B_2\times B_2}(0,0;y,w) \d y \dw \d\tau.
\end{aligned}
\]

Now, letting~$\tau := \eps^{2s}\tau'$,~$y:= \eps^{1+2s}y'$,~$w := \eps w$, and recalling the scaling property in~\eqref{eq:fundamental-dilation} we obtain
\[
\begin{aligned}
f_\eps(0)
 \ \geq  \ c_1
\eps^{2s}\int_{\frac{T}{2}}^{T}
\iint_{B_1\times B_1}
 P_{\tau}^{B_\frac{2}{\eps^{1+2s}}\times B_\frac{2}{\eps}}(0,0;y,w) \d y \dw \d\tau.
\end{aligned}
\]

Hence,~\eqref{eq:maximum-principle} together with Lemma~\ref{psharp:eq_kernel_mass} gives
\[
f_\eps(0)
 \ \geq  \ \frac{c_1T
\eps^{2s}}{2}\inf_{t \in [T/4,T]\times B_{\frac{1}{16}}\times B_{\frac{1}{16}}}\iint_{B_1\times B_1}
P_t^{B_2\times B_2}(x,v;y,w)\d y\d w \ \geq \  \frac{c_1T
\eps^{2s}}{4}\,,
\]
namely,~\eqref{psharp:eq_lower_set} with~$c:=c_1T/4$.

We now prove~\eqref{psharp:eq_tail_estimate}. First of all, let us note that by~\eqref{def:f-eps} we obtain
\begin{equation}\label{psharp:eq_tail_Linfty}
\|\tail(f_\eps;B_\frac{1}{2})(t,\cdot)\|_{L^\infty(\er^n)}
\ \leq  \  C\,,
\end{equation}
for every~$t\in\er$ and for some~$C \equiv C(n,s)>0$. 

For the exterior contribution, by the choice of the cut-off~$\eta_\eps$ and~$\psi$, we have
\begin{equation}\label{psharp:eq_tail_boundary_L1}
\begin{aligned}
\int_{\er^n}\tail(b_\eps;B_\frac{1}{2})\dx
& \ = \ 2^{-2s}a_\eps(t)\int_{\er^n}\eta_\eps(x)\dx
\int_{\er^n}\frac{\psi(v)}{|v|^{n+2s}}\dv\\
& \ \leq  \ C a_\eps(t)\eps^{n(1+2s)}\,,
\end{aligned}
\end{equation}
for~$C \equiv C(n,s)>0$.

Also, the correction is supported in~$B_{2}$ in the velocity variables, which by~\eqref{psharp:eq_correction_bounds} gives
\begin{equation}\label{psharp:eq_tail_correction_L1}
\int_{\er^n}\tail(w_\eps^{\rm corr};B_\frac{1}{2})(t,x)\dx
\ \leq \ C\norma{w_\eps^{\rm corr}(t)}_{L^1(B_2\times B_2)}
 \ \leq \ C\eps^{N_{\tt d}}.
\end{equation}

Combining~\eqref{psharp:eq_tail_boundary_L1} and~\eqref{psharp:eq_tail_correction_L1}, for~$\eps \in (0,1)$, gives
\begin{equation}\label{psharp:eq_tail_L1}
\|\tail(f_\eps;B_\frac{1}{2})(t,\cdot)\|_{L^1(\er^n)}
\ \leq \ C\,\eps^{n(1+2s)}\,,
\end{equation}
for every~$t \in (-2,1/2)$. 

Interpolating~\eqref{psharp:eq_tail_Linfty} and~\eqref{psharp:eq_tail_L1} yields
\[
\|\tail(f_\eps;B_\frac{1}{2})(t,\cdot)\|_{L^p(\er^n)}^p
\ \leq \ C\eps^{n(1+2s)}.
\]
which, after integration in time,
\[
\iint_{U_1}
\tail(f_\eps;B_\frac{1}{2})^p\dx\dt
\ \leq \  C\,\eps^{N_{\tt d}}\,,
\]
i.~\!e.~\eqref{psharp:eq_tail_estimate} follows.

Lastly,~\eqref{psharp:eq_L2_estimate} follows since ~$b_\eps=0$ in~$Q_{1}$, so~$f_\eps=w_\eps^{\rm corr}$ there and~\eqref{psharp:eq_correction_bounds},
gives
\[
\begin{aligned}
\iiint_{Q_1}f_\eps^2\dz
& \ \leq \ \int_{-1}^{0}
\norma{w_\eps^{\rm corr}(t)}_{L^\infty(B_2\times B_2)}
\norma{w_\eps^{\rm corr}(t)}_{L^1(B_2\times B_2)}\dt\\
& \ \leq \  C\eps^{N_{\tt d}+2s}.
\end{aligned}
\]
Now, since~$p<N_{\tt d}/(2s)$ and~$N_{\tt d}>2s$, we have that~$ 2s < \min\left\{N_{\tt d}/p,(N_{\tt d}+2s)/2\right\}$. Therefore,~\eqref{psharp:eq_lower_set},~\eqref{psharp:eq_tail_estimate} and~\eqref{psharp:eq_L2_estimate} yields
\[
\frac{f_\eps(0)}{\nra{f_\eps}_{L^2(Q_1)}+\nra{\tail(f_\eps;B_\frac{1}{2})}_{L^{p}(U_1)}} \ \gtrsim \ \frac{1}{\eps^{\min\left\{\frac{N_{\tt d}}{p},\frac{N_{\tt d}+2s}{2}\right\}-2s}}.
\]
Letting now,~$\eps:= \eps_1/m$, for~$\eps_1 \in (0,\eps_0)$ and $m \in \en$, we obtain the desired property letting~$m\to +\infty$.

\vspace{2mm} 
{\it Case 2: Assume~$p = p^\ast$}.
Fix~$\eps_1\in(0,\eps_0)$, where~$\eps_0$ is given by
Proposition~\ref{prop:solution}, and, for every~$j\in\en$, set
\[
\eps_j:=\frac{\eps_1}{16^{(j-1)/(2s)}}.
\]
so that the intervals
\[
\left(-2T\eps_j^{2s},-\frac{T}{4}\eps_j^{2s}\right),
\qquad j\in\en,
\]
are pairwise disjoint.

Fix~$\gamma\in(1/p^\ast,1)$ and define for $m \in \mathbb{N}$
\[
f_m(t,x,v)
:=
\frac{1}{m^\gamma}
\sum_{j=1}^m
\frac{
b_{\eps_j}(t,x,v)
+w_{\eps_j}^{\mathrm{corr}}(t,x,v)}{\eps_j^{2s}},
\]
where~$b_\eps$ and~$w_\eps^{\mathrm{corr}}$ are given by
\eqref{def:b-eps} and~\eqref{psharp:eq_correction}, respectively.

As in the subcritical case, linearity gives
$f_m\in\Wc(\mathcal{Q})$, and~$f_m$ is a globally nonnegative
weak solution to~\eqref{problema} in~$\mathcal{Q}$ with~$h=0$.
Moreover, for every~$\alpha>0$, let us note that
\begin{equation}\label{eq:counter-sum}
\sum_{j=1}^m\eps_j^\alpha
\leq C(\alpha,\eps_1).
\end{equation}
uniformly in~$m$.

Using the estimates established in
\eqref{psharp:eq_lower_set}--\eqref{psharp:eq_L2_estimate},
together with~\eqref{eq:counter-sum}, we obtain
\[
\begin{aligned}
f_m(0)
&\geq c\,m^{1-\gamma},
\\
\nra{\tail(f_m;B_{\frac{1}{2}})}_{L^{p^\ast}(U_1)}
& \leq C\left(m^{1/p^\ast-\gamma}+m^{-\gamma}\right),
\\
\nra{f_m}_{L^2(Q_1)}
&
\leq C m^{-\gamma}.
\end{aligned}
\]
We only justify the tail estimate, since the others simply follow by triangle inequality of the $L^{p^\ast}$-norm. Indeed, since the time supports are pairwise disjoint, by~\eqref{psharp:eq_tail_boundary_L1} we obtain
\[
\int_{\ern}\tail(b_{\eps_j};B_\frac{1}{2})\dx \leq Ca_{\eps_j}(t)\eps^{n(1+2s)}_j.
\]
On the other hand,
\[
\|\tail(b_{\eps_j};B_\frac{1}{2})\|_{L^\infty(\ern)} \leq C\,,
\]
for some~$C \equiv C(n,s)>0$.

Therefore, by interpolation, we obtain that
\[
\begin{aligned}
\frac{1}{m^{\gamma p^*}}\sum_{j=1}^m\eps_j^{-2sp^*}\int_{-1}^{0}\|
\tail(b_{\eps_j};B_{\frac{1}{2}})
\|_{L^{p^\ast}(B_1)}^{p^\ast}\dt & \ \leq \  \frac{C}{m^{\gamma p^*}}\sum_{j=1}^m\eps_j^{-2sp^*}\int_{-1}^{0}\|
\tail(b_{\eps_j};B_{\frac{1}{2}})
\|_{L^{1}(B_1)}\dt\\
& \ \leq  \
\frac{C}{m^{\gamma p^\ast}}
\sum_{j=1}^m\eps_j^{-2sp^\ast+n(1+2s)}
\int_{-2T\eps_j^{2s}}^{-\frac{T\eps^{2s}_j}{4}}a_{\eps_j}(t)\dt\\
&  \ \leq \ 
\frac{C}{m^{\gamma p^\ast}}
\sum_{j=1}^m1
=C m^{1-\gamma p^\ast}.
\end{aligned}
\]

On the other hand, for the corrections we simply interpolate between $L^1$ and~$L^\infty$ and use the triangle inequality so that, by~\eqref{eq:counter-sum}, we obtain
\[
\frac{1}{m^\gamma}
\sum_{j=1}^m\eps_j^{-2s}
\nra{\tail(w_{\eps_j}^{\mathrm{corr}};B_{\frac{1}{2}})}
_{L^{p^\ast}(U_1)}
 \ \leq \ 
\frac{C}{m^\gamma}\sum_{j=1}^m\eps_j^{2s}
\leq \frac{C}{m^\gamma}.
\]

Since~$1/p^\ast<\gamma<1$, both norms in the denominator tend
to zero, whereas~$f_m(0)$ diverges. More precisely,
\[
\frac{f_m(0)}{\nra{f_m}_{L^2(Q_1)}+\nra{\tail(f_m;B_{\frac{1}{2}})}_{L^{p^\ast}(U_1)}}
\gtrsim  m^{1-1/p^\ast}\,,
\]
which blows up as~$m \to +\infty$.
\hfill $\square$

\vspace{3mm}


\vspace{3mm}




\begin{thebibliography}{99}


\vs \bibitem{AN26}
{P. Auscher, L. Niebel}:
Kinetic Sobolev Spaces.
{\it Preprint}~(2026).
Available at~\href{https://arxiv.org/abs/2603.17491}{\tt arXiv:2603.17491}



\vs \bibitem{AP20}
{F. Anceschi, S. Polidoro}:
A survey on the classical theory for Kolmogorov equation.
{\it Le Matematiche} {\bf LXXV} (2020), no.~1, 221--258.
\href{http://dx.doi.org/10.4418/2020.75.1.11}{\tt DOI:10.4418/2020.75.1.11}
\vs


\bibitem{BCD11}
H. Bahouri, J.-Y. Chemin, R. Danchin:
\emph{Fourier Analysis and Nonlinear Partial Differential Equations},
Grundlehren der mathematischen Wissenschaften, Vol.~343,
Springer, Berlin, Heidelberg, 2011.
\vs

\bibitem{Bou02}
{F. Bouchut}:
Hypoelliptic regularity in kinetic equations.
{\it J. Math. Pures Appl.} {\bf 81} (2002), no.~11, 1135--1159.
\href{http://dx.doi.org/10.1016/S0021-7824(02)01264-3}{\tt DOI:10.1016/S0021-7824}

 \vs \bibitem{CS18}
{L. A. Caffarelli, Y. Sire}:
{\it On some pointwise inequalities involving nonlocal operators}.
Appl. Numer. Harmon. Anal., Birkh\"auser/Springer, Cham, 2017, 1--18.
\href{http://dx.doi.org/10.1007/978-3-319-52742-0_1}{\tt DOI:10.1007/978-3-319-52742-0\_1} 
 
\vs \bibitem{Coz17}
{M. Cozzi}:
Regularity results and Harnack inequalities for minimizers and solutions of nonlocal problems: A unified approach via fractional De Giorgi classes.
{\it J. Funct. Anal.} {\bf 272} (2017), no.~11, 4762--4837.
\href{http://dx.doi.org/10.1016/j.jfa.2017.02.016}{\tt DOI:10.1016/j.jfa.2017.02.016} 
 
 \vs \bibitem{DKP14}
{A. Di Castro, T. Kuusi, G. Palatucci}:
Nonlocal Harnack inequalities.
{\it J. Funct. Anal.} {\bf 267} (2014), no.~6, 1807--1836.
\href{http://dx.doi.org/10.1016/j.jfa.2014.05.023}{\tt DOI:10.1016/j.jfa.2014.05.023} 

\vs \bibitem{DKP16}
{A. Di Castro, T. Kuusi, G. Palatucci}:
Local behavior of fractional $p$-minimizers.
{\it Ann. Inst. H. Poincar\'e Anal. Non Lin\'eaire} {\bf 33} (2016), 1279--1299.
\href{http://dx.doi.org/10.1016/j.anihpc.2015.04.003}{\tt DOI:10.1016/j.anihpc.2015.04.003} 

\vs \bibitem{DH22}
{H. Dieter, J. Hirsch}:
Regularity for rough hypoelliptic equations.
{\it Preprint}~(2022).
\href{https://arxiv.org/abs/2209.08077}{\tt arXiv:2209.08077} 

\vs \bibitem{DPV12}
{E. Di Nezza, G. Palatucci, E. Valdinoci}:
Hitchhiker's guide to the fractional Sobolev spaces.
{\it Bull. Sci. Math.} {\bf 136} (2012), 521--573.
\href{http://dx.doi.org/10.1016/j.bulsci.2011.12.004}{\tt DOI:10.1016/j.bulsci.2011.12.004} 


\vs \bibitem{FRW24}
{X. Fern\'andez-Real, X.~Ros-Oton, M. Weidner}:
Regularity for the Boltzmann equation conditional to pressure and moment bounds.
{\it Commun. Math. Phys.} {\bf 406} (2025), Art.~175.
\href{http://dx.doi.org/10.1007/s00220-025-05356-9}{\tt DOI:10.1007/s00220-025-05356-9}

\vs \bibitem{Fol75}
{G.~\!B.~Folland}:
Subelliptic estimates and function spaces on nilpotent Lie groups.
{\it Ark. Math.} {\bf 13} (1975), 161--207.
\href{http://dx.doi.org/10.1007/BF02386204}{\tt  DOI:10.1007/BF02386204}



\vs \bibitem{GIMV19}
{F.~Golse, C.~Imbert, C.~Mouhot, A.~!F.~Vasseur}:
Harnack inequality for kinetic Fokker-Planck equations with rough coefficients and application to the Landau equation.
{\it Ann. Sc. Norm. Super. Pisa, Cl. Sci. (5)} {\bf 19} (2019), no.~1, 253--295.
\href{http://dx.doi.org/10.2422/2036-2145.201702_001}{\tt DOI:10.2422/2036-2145.201702\_001}

\vs \bibitem{Gru24}
{F. Grube}:
Pointwise estimates of the fundamental solution to the fractional Kolmogorov equation.
{\it Preprint}~(2024).
\href{https://arxiv.org/abs/2411.00687}{\tt  arXiv:2411.00687}

\vs \bibitem{GI23}
{J.~Guerand, C.~Imbert}:
Log-transform and the weak {H}arnack inequality for kinetic {F}okker-{P}lanck equations.
{\it J. Inst. Math. Jussieu} {\bf 22} (2023), no.~6, 2749--2774.
\href{http://dx.doi.org/10.1017/S1474748022000160}{\tt   DOI:10.1017/S1474748022000160}



\vs \bibitem{GM22}
{J. Guerand, C. Mouhot}:
Quantitative De Giorgi methods in kinetic theory.
{\it J. \'Ec. polytech. Math.} {\bf 9} (2022), 1159--1181.
\href{http://dx.doi.org/10.5802/jep.203}{\tt  DOI:10.5802/jep.203}


\vs \bibitem{HPZ21}
{Z. Hao, X. Peng, X. Zhang}:
H\"ormander's Hypoelliptic Theorem for Nonlocal Operators.
{\it J. Theor. Probability} {\bf 34} (2021), 1870--1916.
\href{http://dx.doi.org/10.1007/s10959-020-01020-1}{\tt   DOI:10.1007/s10959-020-01020-1}

{
\vs \bibitem{Hor67}
{L. H{\"o}rmander}: Hypoelliptic second order differential equations.
{\em Acta Math.} {\bf 119} (1967), 147--171.
\href{http://dx.doi.org/10.1007/BF02392081}{\tt  DOI:10.1007/BF02392081}
}

\vs \bibitem{HZ24}
{H. Hou, X.~Zhang}:
Heat kernel estimates for nonlocal kinetic operators.
{\it Prob. Theory  Related Fields}~(2026).
\href{http://dx.doi.org/10.1007/s00440-026-01523-8}{\tt  DOI:10.1007/s00440-026-01523-8}

\vs \bibitem{IMS20}
{C. Imbert, C. Mouhot, L. Silvestre}:
Decay estimates for large velocities in the Boltzmann equation without cutoff.
{\it J. \'Ec. polytech. Math.} {\bf 7} (2020), 143--184.
\href{http://dx.doi.org/10.5802/jep.113}{\tt  DOI:10.5802/jep.113}

\vs \bibitem{IS20a}
{C.~Imbert, L.~Silvestre}:
Regularity for the Boltzmann equation conditional to macroscopic bounds.
{\it EMS Surv. Math. Sci.} {\bf 7} (2020), no.~1, 117--172.
\href{http://dx.doi.org/10.4171/EMSS/37}{\tt  DOI:10.4171/EMSS/37}

\vs \bibitem{IS20}
{C.~Imbert, L.~Silvestre}:
The weak Harnack inequality for the Boltzmann equation without cut-off.
{\it J. Eur. Math. Soc. (JEMS)} {\bf 22} (2020), no.~2, 507--592.
\href{http://dx.doi.org/10.4171/JEMS/928}{\tt   DOI:10.4171/JEMS/928}

\vs \bibitem{IS22}
{C.~Imbert, L.~Silvestre}:
Global regularity estimates for the Boltzmann equation without cut-off.
{\it J. Amer. Math. Soc. (JAMS)} {\bf 35} (2022), no.~3, 625--703.
\href{http://dx.doi.org/10.1090/jams/986}{\tt  DOI:10.1090/jams/986}



\vs \bibitem{Kas07}
{M.~Kassmann}:
The classical Harnack inequality fails for nonlocal operators.
{\it SFB 611-preprint} {\bf 360} (2007).
\href{https://citeseerx.ist.psu.edu/viewdoc/download?doi=10.1.1.454.223andrep=rep1andtype=pdf}{\tt  https://citeseerx.ist.psu.edu/viewdoc/download?doi=10.1.1.454.223})

\vs \bibitem{Kas09}
{M.~Kassmann}:
A priori estimates for integro-differential operators with measurable kernels.
{\it Calc. Var. Partial Differential Equations} {\bf 34} (2009), 11--21.
\href{http://dx.doi.org/10.1007/s00526-008-0173-6}{\tt  DOI:10.1007/s00526-008-0173-6}

\vs \bibitem{Kas11}
{M.~Kassmann}:
Harnack inequalities and H\"older regularity estimates for nonlocal operator revisited.
{\it SFB 11015-preprint} (2011).
\href{https://sfb701.math.uni-bielefeld.de/preprints/sfb11015.pdf}{\tt https://sfb701.math.uni-bielefeld.de/preprints/sfb11015.pdf})

\vs \bibitem{KW24}
{M. Kassmann, M. Weidner}:
Nonlocal operators related to nonsymmetric forms II: Harnack inequalities.
{\it Anal. PDE} {\bf 17} (2024), no.~9, 3189--3249.
\href{http://dx.doi.org/10.2140/apde.2024.17.3189}{\tt  DOI:10.2140/apde.2024.17.3189}

\vs \bibitem{KW23}
{M. Kassmann, M. Weidner}:
The parabolic Harnack inequality for nonlocal equations.
{\it Duke Math. J.} {\bf 173} (2024), no.~17, 3413--3451.
\href{http://dx.doi.org/10.1215/00127094-2024-0008}{\tt  DOI:10.1215/00127094-2024-0008}

\vs \bibitem{KW24c}
{M.~Kassmann, M. Weidner}:
The Harnack inequality fails for nonlocal kinetic equations.
{\it Adv. Math.} {\bf 459} (2024), Art.~110030.
\href{http://dx.doi.org/10.1016/j.aim.2024.110030}{\tt  DOI:10.1016/j.aim.2024.110030}



{
\vs \bibitem{Kol34}
{A. Kolmogorov}:
Zuf{\"a}llige {Bewegungen}. ({Zur} {Theorie} der {Brownschen} {Bewegung}.).
{\it Ann. Math.} {\bf 35} (2) (1934), 116--117.
\href{http://dx.doi.org/10.2307/1968123}{\tt DOI:10.2307/1968123}
}

\vs \bibitem{KKP16}
{J. Korvenp\"a\"a, T. Kuusi, G. Palatucci}:
The obstacle problem for nonlinear integro-differential operators.
{\it Calc. Var. Partial Differential Equations} {\bf 55} (2016), no.~3, Art.~63.
\href{http://dx.doi.org/10.1007/s00526-016-0999-2}{\tt   DOI:10.1007/s00526-016-0999-2}

\vs \bibitem{KKP17}
{J. Korvenp\"a\"a, T. Kuusi, G. Palatucci}:
Fractional superharmonic functions and the Perron method for nonlinear integro-differential equations.
{\it Math. Ann.} {\bf 369} (2017), no.~3--4, 1443--1489.
\href{http://dx.doi.org/10.1007/s00208-016-1495-x}{\tt   DOI:10.1007/s00208-016-1495-x}


\vs \bibitem{Loh22}
{A. Loher}:
Quantitative De Giorgi methods in kinetic theory for non-local operators.
{\it J. Funct. Anal.} {\bf 286} (2024), no.~6, Art.~110312.
\href{http://dx.doi.org/10.1016/j.jfa.2023.110312}{\tt   DOI:10.1016/j.jfa.2023.110312}

\vs \bibitem{Loh24c} {A. Loher}: Semi-local behaviour of non-local hypoelliptic equations: divergence form. \href{https://arxiv.org/abs/2404.05612}{\tt arXiv:2404.05612v4} (2026).

\vs \bibitem{Mou18}
{C. Mouhot}:
De Giorgi-Nash-Moser and H\"ormander theories: new interplays.
In {\it Proceedings of the International Congress of Mathematicians--Rio de Janeiro 2018}, {\bf Vol.~III} (2018), 2467--2493.



\vs \bibitem{OS23}
{Z. Ouyang, L. Silvestre}:
Conditional $L^\infty$ estimates for the non-cutoff Boltzmann equations in a bounded domain.
{\it Arch. Rational Mech. Anal.} {\bf 248} (2024), Art.~59.
\href{http://dx.doi.org/10.1007/s00205-024-02002-x}{\tt  DOI:10.1007/s00205-024-02002-x}

\vs \bibitem{PasBook}
{A. Pascucci}:
{\it PDE and Martingale Methods in Option Pricing.}
Bocconi \& Springer Series, ISBN 978-88-470-1780-1 (2010) 

\vs \bibitem{PP04}
{A. Pascucci, S. Polidoro}:
The Moser's iterative method for a class of ultraparabolic equations.
{\it Comm. Cont. Math.} {\bf 6} (2004), no.~3, 395--417.
\href{http://dx.doi.org/10.1142/S0219199704001355}{\tt   DOI:10.1142/S0219199704001355}


\vs {\bibitem{Sch14} {R.~\!L. Schilling}: \emph{An Introduction to Lévy and Feller Processes. Advanced Courses in Mathematics - CRM Barcelona} 2014. \href{https://arxiv.org/abs/1603.00251}{\tt   arXiv1603.00251}
}



\vs \bibitem{Sil16}
{L. Silvestre}:
A New Regularization Mechanism for the Boltzmann Equation Without Cut-Off.
{\it Commun. Math. Phys.} {\bf 348} (2016), 69--100.
\href{http://dx.doi.org/10.1007/s00220-016-2757-x}{\tt   DOI:10.1007/s00220-016-2757-x}

\vs \bibitem{Sto19}
{L.~\!F. Stokols}:
H\"older continuity for a family of nonlocal hypoelliptic kinetic equations.
{\it SIAM J. Math. Anal.} {\bf 51} (2019), no.~6, 4815--4847.
\href{http://dx.doi.org/10.1137/18M1234953}{\tt  DOI:10.1137/18M1234953}

\vs \bibitem{WZ11}
{W. Wang, L. Zhang}:
The $C^\alpha$ regularity of weak solutions of ultraparabolic equations.
{\it Discrete Contin. Dyn. Syst.} {\bf 29} (2011), no.~3, 1261--1275.
\href{http://dx.doi.org/10.3934/dcds.2011.29.1261}{\tt  DOI:10.3934/dcds.2011.29.1261}


\end{thebibliography}
\end{document}